\documentclass[12pt]{article}
\usepackage{graphics}
\usepackage{epstopdf}
\usepackage{float}

\usepackage{fullpage}
\usepackage{amsfonts}
\usepackage{amsmath}
\usepackage{amssymb}
\usepackage{latexsym}
\usepackage{indentfirst}
\usepackage{fancyhdr}
\usepackage{stmaryrd}
\usepackage{mathrsfs}
\usepackage{cite}
\usepackage[all]{xy}
\usepackage{amscd,amsthm,graphicx, epsfig}
\usepackage{color}
\usepackage[numbers,sort&compress]{natbib}

\numberwithin{equation}{section}

\newtheorem{thm}{Theorem}[section]
\newtheorem{defn}[thm]{Definition}
\newtheorem{prop}[thm]{Proposition}

\newtheorem{lemma}[thm]{Lemma}
\newtheorem{rema}[thm]{Remark}

\newcommand{\halmos}{\rule{1ex}{1.4ex}}

\newcommand{\ki}{{\rm Ker}\; \mathbf{I}}
\newcommand{\kp}{{\rm Ker}\; \mathbf{P}}

 \newcommand{\res}{\mbox{\rm Res}}
 \newcommand{\re}{\mbox{\rm Re}}
 \newcommand{\im}{\mbox{\rm Im}}
 
\renewcommand{\hom}{\mbox{\rm Hom}}

 \newcommand{\pf}{{\it Proof.}\hspace{2ex}}
 \newcommand{\epf}{\hspace*{\fill}\mbox{$\halmos$}}
 \newcommand{\epfv}{\hspace*{\fill}\mbox{$\halmos$}\vspace{1em}}

\newcommand{\lbar}{\bigg\vert}

\newcommand{\A}{\mathcal{A}}
\newcommand{\Y}{\mathcal{Y}}
\newcommand{\C}{\mathbb{C}}
\newcommand{\Z}{\mathbb{Z}}

\newcommand{\F}{\mathcal{F}}

\newcommand{\V}{\mathcal{V}}

\newcommand{\nno}{\nonumber}
\newcommand{\bea}{\begin{eqnarray}}
\newcommand{\eea}{\end{eqnarray}}

\newcommand{\be}{\begin {equation}}

\newcommand{\ee}{\end{equation}}

\newcommand{\mbar}{\vert}

\makeatletter

\newcommand{\Rmnum}[1]{\expandafter\@slowromancap\romannumeral#1@}
\makeatother

\begin{document}

\bibliographystyle{unsrt}
\baselineskip=16pt

 \title{{\bf An $S_3$-symmetry of the Jacobi Identity for Intertwining Operator Algebras}}
 \author{Ling Chen}
    \date{}
    \maketitle

\begin{abstract}
We prove an $S_{3}$-symmetry of the Jacobi identity for intertwining operator algebras. Since this Jacobi identity involves the braiding and fusing isomorphisms satisfying the genus-zero Moore-Seiberg equations, our proof uses not only the basic properties of intertwining operators, but also the properties of braiding and fusing isomorphisms and the genus-zero Moore-Seiberg equations. Our proof depends heavily on the theory of multivalued analytic functions of several variables, especially the theory of analytic extensions.

$\vspace{0.1cm}$



\end{abstract}

\section{Introduction}

Intertwining operator algebras were introduced and studied by Huang in \cite{H1,H2}.
In \cite{C}, the author studied intertwining operator algebras in a setting more general than \cite{H2}. In particular, the duality properties, Jacobi identity, Moore-Seiberg equations, locality and some other properties of intertwining operator algebras were studied. For the background on  intertwining operator algebras, we refer the reader to \cite{H1,H2,C}.

For vertex operator algebras, the Jacobi identity has an $S_3$-symmetry which corresponds to the obvious $S_3$-symmetry of the Jacobi identity for Lie algebras \cite{FHL}. For abelian intertwining operator algebras (see \cite{DL1,DL2}), Guo \cite{G} proved that the Jacobi identity for these algebras also has an $S_3$-symmetry. In this paper, we prove an $S_3$-symmetry of the Jacobi identity for intertwining operator algebras introduced by Huang \cite{H2} and studied by the author \cite{C}. See Theorem \ref{t1} for the statement of this $S_3$-symmetry. The $S_3$-symmetry in this general case is much more complicated but is also much more interesting and much deeper. Note that the Jacobi identity for general intertwining operator algebras in \cite{H2} and \cite{C} involves the braiding and fusing isomorphisms satisfying the genus-zero Moore-Seiberg equations. The proof of the $S_3$-symmetry in the present paper uses not only the properties of the intertwining operators (for example, the skew-symmetry) but also the properties of braiding and fusing isomorphisms and the genus-zero Moore-Seiberg equations. In particular, our proof depends heavily on the theory of multivalued analytic functions of several variables, especially the theory of analytic extensions.

This paper is organized as follows. In section 2, we review some preliminaries concerning the theory of intertwining operator algebras which we need to formulate and prove the main result of this paper. In section 3, we prove an $S_3$-symmetry of the Jacobi identity for intertwining operator algebras.

\paragraph{Acknowledgments}
The author is
very grateful to Professor Yi-Zhi Huang for his support,
encouragement, many discussions on the paper
and help with the exposition of the paper. The author is supported by NSFC
grant 11401559.

\section{Preliminaries}

In this section, we first recall some notations and facts in formal calculus
and complex analysis (see \cite{FLM2,FHL,H2} for more details), then we review some definitions and properties in the theory of
intertwining operator algebras in \cite{H2,C}. These are necessary preliminaries for formulating and proving the main result of this paper.

In this paper, as in \cite{FHL,H2,C}, $x$, $x_{0}, \dots$ are independent
commuting formal variables. And
for a vector space $W$ and a formal variable $x$, as in \cite{FHL,H2,C}, we
shall use $W[x]$, $W[x,x^{-1}]$, $W[[x]]$, $W[[x,x^{-1}]]$, $W((x))$ and $W\{x\}$ to denote the spaces of all polynomials in $x$, all Laurent polynomials in $x$, all formal power series in $x$, all formal Laurent series in $x$, all formal Laurent series in $x$ with finitely many negative powers and all formal series with arbitrary powers of $x$ in $\mathbb{C}$, respectively.
For series with more than one formal variables, we shall use similar notations.
We shall use $\res_{x}f(x)$ to denote the
coefficient of $x^{-1}$ in $f(x)$ for any $f(x)\in W\{x\}$.  As in \cite{FHL,H2,C}, $z,
z_{0}, \dots,$ are complex numbers, {\it not} formal
variables.

Let
\begin{equation}
\delta(x)=\sum_{n\in
\mathbb{Z}}x^{n}.
\end{equation}
It has the following important property:
For any $f(x)\in \mathbb{C}[x, x^{-1}]$,
\begin{equation}
f(x)\delta(x)=f(1)\delta(x).
\end{equation}
Following \cite{FHL,H2,C}, we use the convention that negative powers of a
binomial are to be expanded in nonnegative powers of the second
summand so that, for example,
\begin{equation}
x_{0}^{-1}\delta\left(\frac{x_{1}-x_{2}}{x_{0}}\right)=\sum_{n\in \mathbb{Z}}
\frac{(x_{1}-x_{2})^{n}}{x_{0}^{n+1}}=\sum_{m\in \mathbb{N},\; n\in \mathbb{Z}}
(-1)^{m}{{n}\choose {m}} x_{0}^{-n-1}x_{1}^{n-m}x_{2}^{m}.
\end{equation}
The following identities are often very useful:
\begin{equation}
 x_{1}^{-1}\delta\left(\frac{x_{2}+x_{0}}{x_{1}}\right)=x_{2}^{-1}\delta\left(
\frac{x_{1}-x_{0}}{x_{2}}\right),
\end{equation}
\begin{equation}
 x_{0}^{-1}\delta\left(\frac{x_{1}-x_{2}}{x_{0}}\right)-
x_{0}^{-1}\delta\left(\frac{x_{2}-x_{1}}{-x_{0}}\right)=
x_{2}^{-1}\delta\left(\frac{x_{1}-x_{0}}{x_{2}}\right).
\end{equation}

As in \cite{FHL,H2,C}, $\C[x_1,x_2]_S$ is the ring of rational functions obtained by
inverting the products of (zero or more) elements of the set $S$
of nonzero homogenous linear polynomials in $x_1$ and $x_2$.
Also, $\iota_{12}$ is the operation of expanding an element
of $\C[x_1,x_2]_S$, that is, a polynomial in $x_1$ and $x_2$ divided
by a product of homogenous linear polynomials in $x_1$ and $x_2$,
as a formal series containing at most finitely many negative powers
of $x_2$ (using binomial expansions for negative powers of linear
polynomials involving both $x_1$ and $x_2$); similarly for $\iota_{21}$, and so on.
We need the following fact from \cite{FHL}.

\begin{prop}\label{p3.1.1}
Consider a rational function of the form
\begin{equation}
f(x_0,x_1,x_2)=\frac{g(x_0,x_1,x_2)}{x_0^r x_1^s x_2^t},
\end{equation}
where $g$ is a polynomial and $r,s,t\in\Z$. Then
\begin{equation}
 x_{1}^{-1}\delta\left(\frac{x_{2}+x_{0}}{x_{1}}\right)\iota_{20}(f\mbar_{x_1=x_0+x_2})
 =x_{2}^{-1}\delta\left(
\frac{x_{1}-x_{0}}{x_{2}}\right)\iota_{10}(f\mbar_{x_2=x_1-x_0})
\end{equation}
and
\begin{eqnarray}
& x_{0}^{-1}\delta\left(\displaystyle\frac{x_{1}-x_{2}}{x_{0}}\right)
\iota_{12}(f\mbar_{x_0=x_1-x_2})-
x_{0}^{-1}\delta\left(\displaystyle\frac{x_{2}-x_{1}}{-x_{0}}\right)
\iota_{21}(f\mbar_{x_0=x_1-x_2})&\nno\\
&=
x_{2}^{-1}\delta\left(\displaystyle\frac{x_{1}-x_{0}}{x_{2}}\right)
\iota_{10}(f\mbar_{x_2=x_1-x_0}).&
\end{eqnarray}

\end{prop}
\vspace{0.2cm}

As in \cite{FHL,H2,C}, the graded dual of a $\Z$-graded, or more generally, $\C$-graded, vector space
$W=\coprod_n W_{(n)}$ is denoted by
\begin{equation}
W'=\coprod_n W_{(n)}^*.
\end{equation}

For any $z\in \mathbb{C}$, we use $\log z$ to denote the value $\log |z|+i\arg z$ with $0\le\arg z<2\pi$ of logarithm of $z$. For two multivalued functions $f_{1}$ and $f_{2}$ on a region, $f_{1}$ and $f_{2}$ are equal if on any simply
connected open subset of the region, any single-valued branch of
$f_{1}$ is equal to a single-valued branch of $f_{2}$, and vice versa.
\vspace{0.2cm}

Now we recall some basic notions and results in the theory of intertwining
operator algebras. For the details of the definitions and properties of vertex operator algebras, their modules
and intertwining operators,
the reader is referred to \cite{FHL,FLM2,H2}. And for more details of the properties of intertwining operator algebras, the reader is referred to \cite{H2,C}.

Let $(V,Y,{\bf 1},\omega )$ be a vertex operator algebra, and let $W_1,W_2,W_3$ be
modules of $V$.
The space of all intertwining operators of
type ${W_{3}}\choose {W_{1}\ W_{2}}$ is denoted by $\bar{\V}^{W_{3}}_{W_{1}W_{2}}$
instead of $\V^{W_{3}}_{W_{1}W_{2}}$, for as in \cite{C}, the latter shall be used to denote
a subspace of $\bar{\V}^{W_{3}}_{W_{1}W_{2}}$ in the definition of intertwining operator algebra.
The dimension of this vector space
is denoted by
$\bar{\mathcal{N}}^{W_{3}}_{W_{1}W_{2}}$. It is the so-called fusion rule of the same type.
Let $\mathcal{Y}$ be an intertwining operator of type
${W_{3}}\choose {W_{1}W_{2}}$.
Given any $r\in \mathbb{Z}$, as in \cite{HL,H2,C},  we define
\begin{equation}
\Omega_{r}(\mathcal{Y}):W_2\otimes W_1 \rightarrow  W_3\{ x\}
\end{equation}
by
\begin{equation}
\Omega_{r}(\mathcal{Y})(w_{(2)},x)w_{(1)} = e^{xL(-1)}
\mathcal{Y}(w_{(1)},e^{ (2r+1)\pi i}x)w_{(2)}
\end{equation}
for  $w_{(1)}\in W_{1}$, $w_{(2)}\in W_{2}$. We have the following result
proved in \cite{HL}:

\begin{prop}
For any $\Y\in\bar{\V}^{W_{3}}_{W_{1}W_{2}}$, $r\in \mathbb{Z}$, we have  $\Omega_{r}(\mathcal{Y})\in\bar{\V}^{W_{3}}_{W_{2}W_{1}}$. Moreover,
\begin{equation}
\Omega_{-r-1}(\Omega_{r}(\mathcal{Y}))=\Omega_{r}(\Omega_{-r-1}(\mathcal{Y}))
 = \mathcal{Y}.
\end{equation}
In particular, the correspondence  $\mathcal{Y} \mapsto
\Omega_{r}(\mathcal{Y})$  defines a linear isomorphism {from}
$\bar{\V}^{W_{3}}_{W_{1}W_{2}}$ to
$\bar{\V}^{W_{3}}_{W_{2}W_{1}}$,  and we have
\begin{equation}
\bar{\mathcal{N}}^{W_{3}}_{W_{1}W_{2}} = \bar{\mathcal{N}}^{W_{3}}_{W_{2}W_{1}}.
\end{equation}
\end{prop}
\vspace{0.2cm}

Now we recall the first definition of intertwining operator
algebras in \cite{H2}:

\begin{defn}[{\bf Intertwining operator algebra}]\label{i1}
{\rm An {\it intertwining operator algebra of central charge
$c\in \mathbb{ C}$} consists of the following data:

\begin{enumerate}

\item A vector space
\begin{equation}
W=\coprod_{a\in
\mathcal{A}}W^{a}
\end{equation}
graded
by a finite set $\mathcal{ A}$
containing a special element $e$
(graded  by {\it color}).

\item A vertex operator algebra
structure of central charge $c$
on $W^{e}$, and a $W^{e}$-module structure on $W^{a}$ for
each $a\in \mathcal{ A}$.

\item A subspace $\mathcal{ V}_{a_{1}a_{2}}^{a_{3}}$ of
the space of all intertwining operators of type
${W^{a_{3}}\choose W^{a_{1}}W^{a_{2}}}$ for  each triple
$a_{1}, a_{2}, a_{3}\in \mathcal{ A}$, with its dimension denoted by
$\mathcal{ N}_{a_{1}a_{2}}^{a_3}$.

\end{enumerate}

\noindent These data satisfy the
following axioms for any $a_{1}, a_{2},
a_{3}, a_{4}, a_{5},
a_{6}\in \mathcal{ A}$,
$w_{(a_{i})}\in W^{a_{i}}$, $i=1, 2, 3$, and $w_{(a_{4})}'
\in (W^{a_{4}})'$:

\begin{enumerate}

\item The $W^{e}$-module structure on $W^{e}$ is the adjoint module
structure. For any $a\in \mathcal{ A}$, the space $\mathcal{ V}_{ea}^{a}$ is the
one-dimensional vector space spanned by the vertex operator for the
$W^{e}$-module $W^{a}$.
For any $a_{1}, a_{2}\in \mathcal{ A}$ such that
$a_{1}\ne a_{2}$, $\mathcal{ V}_{ea_{1}}^{a_{2}}=0$.

\item {\it Weight condition}: For any $a\in \mathcal{ A}$ and the
corresponding module $W^a=\coprod_{n\in \C}W_{(n)}^a$ graded by the action of
$L(0)$, there exists
$h_{a}\in \mathbb{ R}$ such that $W^{a}_{(n)}=0$ for $n\not
\in h_{a}+\mathbb{ Z}$.

\item {\it Convergence properties}: For any $m\in \mathbb{ Z}_{+}$,
$a_{i}, b_{j}\in \mathcal{ A}$, $w_{(a_{i})}
\in W^{a_{i}}$, $\mathcal{ Y}_{i}\in \mathcal{
V}_{a_{i}\;b_{i+1}}^{b_{i}}$, $i=1, \dots, m$, $j=1, \dots, m+1$, $w_{(b_{1})}'
\in (W^{b_{1}})'$ and
$w_{(b_{m+1})}\in W^{b_{m+1}}$, the series
\begin{equation}\label{conv-pr}
\langle w_{(b_{1})}', \mathcal{ Y}_{1}(w_{(a_{1})}, x_{1})
\cdots\mathcal{ Y}_{m}(w_{(a_{m})},
x_{m})w_{(b_{m+1})}\rangle_{W^{b_{1}}}\mbar_{x^{n}_{i}=e^{n\log z_{i}},\;
i=1, \dots, m,\; n\in \mathbb{ R}}
\end{equation}
is absolutely convergent when $|z_{1}|>\cdots >|z_{m}|>0$
and its sum can be analytically extended to a multivalued analytic function on the region given by $z_i\not= 0$, $i=1, \dots, m$, $z_i\not=z_j$, $i\not=j$, such that
for any set of possible singular points with either $z_i=0$, $z_i=\infty$ or $z_i=z_j$ for $i\not=j$, this multivalued analytic function can be expanded near the singularity as a series having the same form as the expansion near the singular points of a solution of a system of differential equations with regular singular points.
 For
any $\mathcal{ Y}_{1}\in \mathcal{
V}_{a_{1}a_{2}}^{a_{5}}$ and $\mathcal{ Y}_{2}\in \mathcal{
V}_{a_{5}a_{3}}^{a_{4}}$, the series
\begin{equation}\label{conv-it}
\langle w_{(a_{4})}', \mathcal{ Y}_{2}(\mathcal{ Y}_{1}(w_{(a_{1})},
x_{0})w_{(a_{2})},
x_{2})w_{(a_{3})}\rangle_{W^{a_{4}}}
\mbar_{x^{n}_{0}=e^{n\log (z_{1}-z_{2})},\;
x^{n}_{2}=e^{n\log z_{2}},\; n\in \mathbb{ R}}
\end{equation}
is absolutely convergent when
$|z_{2}|>|z_{1}-z_{2}|>0$.

\item {\it Associativity}: For any $\mathcal{ Y}_{1}\in \mathcal{
V}_{a_{1}a_{5}}^{a_{4}}$ and $\mathcal{ Y}_{2}\in\mathcal{
V}_{a_{2}a_{3}}^{a_{5}}$, there exist $\mathcal{ Y}^{a}_{3, i}
\in \mathcal{
V}_{a_{1}a_{2}}^{a}$ and $\mathcal{ Y}^{a}_{4, i}\in \mathcal{
V}_{aa_{3}}^{a_{4}}$ for $i=1, \dots,
\mathcal{ N}_{a_{1}a_{2}}^{a}
\mathcal{N}_{aa_{3}}^{a_{4}}$ and $a\in \mathcal{ A}$,
such that the (multivalued) analytic function
\begin{equation}\label{prod}
\langle w_{(a_{4})}',
\mathcal{ Y}_{1}(w_{(a_{1})}, x_{1})\mathcal{ Y}_{2}(w_{(a_{2})},
x_{2})w_{(a_{3})}\rangle_{W^{a_{4}}}\mbar_{x_{1}=z_{1},
x_{2}=z_{2}}
\end{equation}
defined on the region
$|z_{1}|>|z_{2}|>0$
and the (multivalued) analytic function
\begin{equation}\label{iter}
\sum_{a\in \mathcal{ A}}\sum_{i=1}^{\mathcal{ N}_{a_{1}a_{2}}^{a}\mathcal{
N}_{aa_{3}}^{a_{4}}}\langle w_{(a_{4})}', \mathcal{ Y}^{a}_{4, i}
(\mathcal{ Y}^{a}_{3, i}(w_{(a_{1})},
x_{0})w_{(a_{2})}, x_{2})w_{(a_{3})}\rangle_{W^{a_{4}}}
\lbar_{x_{0}=z_{1}-z_{2},
x_{2}=z_{2}}
\end{equation}
defined on the region
$|z_{2}|>|z_{1}-z_{2}|>0$ are equal on the intersection
$|z_{1}|> |z_{2}|>|z_{1}-z_{2}|>0$.

\item {\it Skew-symmetry}: The restriction of  $\Omega_{-1}$ to
$\mathcal{ V}_{a_{1}a_{2}}^{a_{3}}$ is an isomorphism from
$\mathcal{ V}_{a_{1}a_{2}}^{a_{3}}$ to
$\mathcal{ V}_{a_{2}a_{1}}^{a_{3}}$.

\end{enumerate}}
\end{defn}

\begin{rema}
The skew-symmetry isomorphisms
$\Omega_{-1}(a_{1}, a_{2}; a_{3})$
for all $a_{1}, a_{2}, a_{3}\in \mathcal{ A}$ give an isomorphism
\begin{equation}
\Omega_{-1}: \coprod_{a_{1}, a_{2}, a_{3}\in \mathcal{ A}}
\mathcal{ V}_{a_{1}a_{2}}^{a_{3}}\to \coprod_{a_{1}, a_{2}, a_{3}\in \mathcal{ A}}
\mathcal{ V}_{a_{1}a_{2}}^{a_{3}},
\end{equation}
which, as in \cite{H2,C}, is still called the {\it skew-symmetry isomorphism}.
In this paper, as in \cite{H2,C}, we shall omit subscript $-1$ in
$\Omega_{-1}$ for simplicity and denote it
by $\Omega$.
\end{rema}

We denote the intertwining operator algebra just defined by
$(W,
\mathcal{ A}, \{\mathcal{ V}_{a_{1}a_{2}}^{a_{3}}\}, {\bf 1},
\omega)$
 or simply by $W$.
\vspace{0.3cm}

Next, as in \cite{H2,C}, we give the two linear maps corresponding to the multiplication and
iterates of
intertwining operators, respectively. Let
\begin{eqnarray}
{\bf P}: \coprod_{a_{1}, a_{2}, a_{3}, a_{4},
a_{5}\in \mathcal{ A}}\mathcal{ V}_{a_{1}a_{5}}^{a_{4}}\otimes
\mathcal{ V}_{a_{2}a_{3}}^{a_{5}}&\to &
(\hom(W\otimes W\otimes W, W))\{x_{1}, x_{2}\}\nno\\
\mathcal{ Z}&\mapsto& {\bf P}(\mathcal{ Z})
\end{eqnarray}
be the linear map defined using products of intertwining operators
as follows: For
\begin{equation}
\mathcal{ Z}\in
\coprod_{a_{1}, a_{2}, a_{3}, a_{4},
a_{5}\in \mathcal{ A}}\mathcal{ V}_{a_{1}a_{5}}^{a_{4}}\otimes
\mathcal{ V}_{a_{2}a_{3}}^{a_{5}},
\end{equation}
the element ${\bf P}(\mathcal{ Z})$ to be defined
is a linear map from $W\otimes W \otimes W$ to
$W\{x_{1}, x_{2}\}$.
We denote the image of $w_{1}\otimes
w_{2}\otimes w_{3}$ under this map by
$({\bf P}(\mathcal{ Z}))(w_{1},  w_{2}, w_{3}; x_{1}, x_{2})$ for any $w_{1}, w_{2}, w_{3}\in W$.
Then we define ${\bf P}$ by linearity and by
\begin{eqnarray}
\lefteqn{({\bf P}(\mathcal{ Y}_{1}\otimes
\mathcal{ Y}_{2}))(w_{(a_{6})},  w_{(a_{7})}, w_{(a_{8})}; x_{1}, x_{2})}\nno\\
&&=\left\{\begin{array}{ll}
\mathcal{ Y}_{1}(w_{(a_{6})}, x_{1})\mathcal{ Y}_{2}(w_{(a_{7})}, x_{2})w_{(a_{8})},&
a_{6}=a_{1}, a_{7}=a_{2}, a_{8}=a_{3},\\
0,&\mbox{\rm otherwise}
\end{array}\right.
\end{eqnarray}
for $a_{1}, \dots, a_{8}\in \mathcal{ A}$, $\mathcal{ Y}_{1}\in \mathcal{
V}_{a_{1}a_{5}}^{a_{4}}$, $\mathcal{ Y}_{2}\in \mathcal{
V}_{a_{2}a_{3}}^{a_{5}}$, and $w_{(a_{6})}\in W^{a_{6}}$,
$w_{(a_{7})}\in W^{a_{7}}$, $w_{(a_{8})}\in W^{a_{8}}$. So we have an isomorphism
\begin{equation}
\tilde{\bf P}:\ \ \frac{ \displaystyle \coprod_{a_{1}, a_{2}, a_{3}, a_{4},
a_{5}\in \mathcal{ A}}\mathcal{ V}_{a_{1}a_{5}}^{a_{4}}\otimes
\mathcal{ V}_{a_{2}a_{3}}^{a_{5}}}{\kp} \longrightarrow
{\bf P}\left(\coprod_{a_{1}, a_{2}, a_{3}, a_{4},
a_{5}\in \mathcal{ A}}\mathcal{ V}_{a_{1}a_{5}}^{a_{4}}\otimes
\mathcal{ V}_{a_{2}a_{3}}^{a_{5}}\right)
\end{equation}
which makes the following diagram commute:
\begin{equation}
\xymatrix{
                \displaystyle  \coprod_{a_{1}, a_{2}, a_{3}, a_{4},
a_{5}\in \mathcal{ A}}\mathcal{ V}_{a_{1}a_{5}}^{a_{4}}\otimes
\mathcal{ V}_{a_{2}a_{3}}^{a_{5}}  \ar[d]_{ \pi_P } \ar[r]^-{ {\bf P} }
&   {\bf P}\left(\displaystyle \coprod_{a_{1}, a_{2}, a_{3}, a_{4},
a_{5}\in \mathcal{ A}}\mathcal{ V}_{a_{1}a_{5}}^{a_{4}}\otimes
\mathcal{ V}_{a_{2}a_{3}}^{a_{5}}\right)      \\
   \frac{ \displaystyle \coprod_{a_{1}, a_{2}, a_{3}, a_{4},
a_{5}\in \mathcal{ A}}\mathcal{ V}_{a_{1}a_{5}}^{a_{4}}\otimes
\mathcal{ V}_{a_{2}a_{3}}^{a_{5}}}{\displaystyle \kp}  \ar[ur]^{ {\bf\tilde{ P}}}},
\end{equation}
where $\pi_P$ is the corresponding canonical projective map. As in \cite{C}, we also denote $\pi_P(\mathcal{ Z})$ by
$[\mathcal{ Z}]_P$ or $\mathcal{ Z}+\kp$ for
$\mathcal{ Z}\in
\coprod_{a_{1}, a_{2}, a_{3}, a_{4},
a_{5}\in \mathcal{ A}}\mathcal{ V}_{a_{1}a_{5}}^{a_{4}}\otimes
\mathcal{ V}_{a_{2}a_{3}}^{a_{5}}$ when
there is no ambiguity.
The second linear map is
\begin{eqnarray}
{\bf I}: \coprod_{a_{1}, a_{2}, a_{3}, a_{4},
a_{5}\in \mathcal{ A}}\mathcal{ V}_{a_{1}a_{2}}^{a_{5}}\otimes \mathcal{
V}_{a_{5}a_{3}}^{a_{4}}&\to&
(\hom(W\otimes W\otimes W, W))\{x_{0}, x_{2}\}\nno\\
\mathcal{ Z}&\mapsto& {\bf I}(\mathcal{ Z})
\end{eqnarray}
defined using iterates of intertwining operators as follows:
For
\begin{equation}
\mathcal{ Z}\in \coprod_{a_{1}, a_{2}, a_{3}, a_{4},
a_{5}\in \mathcal{ A}}\mathcal{ V}_{a_{1}a_{2}}^{a_{5}}\otimes \mathcal{
V}_{a_{5}a_{3}}^{a_{4}},
\end{equation}
the element ${\bf I}(\mathcal{ Z})$ to be defined
is a linear map from $W\otimes W \otimes W$ to
$W\{x_{0}, x_{2}\}$. We denote the image of $w_{1}\otimes
w_{2}\otimes w_{3}$ under this map by
$({\bf I}(\mathcal{ Z}))(w_{1}, w_{2}, w_{3};
x_{0}, x_{2})$ for any $w_{1}, w_{2}, w_{3}\in W$. Then we define ${\bf I}$ by linearity and by
\begin{eqnarray}
\lefteqn{({\bf I}(\mathcal{ Y}_{1}\otimes
\mathcal{ Y}_{2}))(w_{(a_{6})}, w_{(a_{7})}, w_{(a_{8})}; x_{0}, x_{2})}\nno\\
&&=\left\{\begin{array}{ll}
\mathcal{ Y}_{2}(\mathcal{ Y}_{1}(w_{(a_6)}, x_{0})w_{(a_7)}, x_{2})w_{(a_{8})},
&a_{6}=a_{1}, a_{7}=a_{2}, a_{8}=a_{3},\\
0,&\mbox{\rm otherwise}
\end{array}\right.
\end{eqnarray}
for $a_{1}, \dots, a_{8}\in \mathcal{ A}$, $\mathcal{ Y}_{1}\in \mathcal{
V}_{a_{1}a_{2}}^{a_{5}}$, $\mathcal{ Y}_{2}\in \mathcal{
V}_{a_{5}a_{3}}^{a_{4}}$, and $w_{(a_{6})}\in W^{a_{6}}$,
$w_{(a_{7})}\in W^{a_{7}}$, $w_{(a_{8})}\in W^{a_{8}}$. Therefore we have an isomorphism
\begin{equation}
\tilde{\bf I}: \ \ \ \frac{ \displaystyle \coprod_{a_{1}, a_{2}, a_{3}, a_{4},
a_{5}\in \mathcal{ A}}\mathcal{ V}_{a_{1}a_{2}}^{a_{5}}\otimes
\mathcal{ V}_{a_{5}a_{3}}^{a_{4}}}{\ki}  \longrightarrow  {\bf I}
\left( \displaystyle\coprod_{a_{1}, a_{2}, a_{3}, a_{4},
a_{5}\in \mathcal{ A}}\mathcal{ V}_{a_{1}a_{2}}^{a_{5}}\otimes
\mathcal{ V}_{a_{5}a_{3}}^{a_{4}} \right)
\end{equation}
which makes the following diagram commute:
\begin{equation}
\xymatrix{    \displaystyle  \coprod_{a_{1}, a_{2}, a_{3}, a_{4},
a_{5}\in \mathcal{ A}}\mathcal{ V}_{a_{1}a_{2}}^{a_{5}}\otimes
\mathcal{ V}_{a_{5}a_{3}}^{a_{4}}  \ar[d]_{ \pi_I } \ar[r]^-{ {\bf I} }
&   {\bf I}\left(\displaystyle \coprod_{a_{1}, a_{2}, a_{3}, a_{4},
a_{5}\in \mathcal{ A}}\mathcal{ V}_{a_{1}a_{2}}^{a_{5}}\otimes
\mathcal{ V}_{a_{5}a_{3}}^{a_{4}}\right)\\
   \frac{\displaystyle \coprod_{a_{1}, a_{2}, a_{3}, a_{4},
a_{5}\in \mathcal{ A}}\mathcal{ V}_{a_{1}a_{2}}^{a_{5}}\otimes
\mathcal{ V}_{a_{5}a_{3}}^{a_{4}}}{\displaystyle \ki}  \ar[ur]^{ {\bf\tilde{ I}}}},
\end{equation}
where $\pi_I$ is the corresponding canonical projective map. As in \cite{C}, we also denote $\pi_I(\mathcal{ Z})$ by
$[\mathcal{ Z}]_I$ or $\mathcal{ Z}+\ki$ for
$\mathcal{ Z}\in
\coprod_{a_{1}, a_{2}, a_{3}, a_{4},
a_{5}\in \mathcal{ A}}\mathcal{ V}_{a_{1}a_{2}}^{a_{5}}\otimes \mathcal{
V}_{a_{5}a_{3}}^{a_{4}}$ when
there is no ambiguity.

The two linear maps
${\bf P}$ and ${\bf I}$ are called the {\it multiplication
of intertwining operators} and
the {\it iterates of intertwining operators}, respectively. \vspace{0.3cm}

Moreover, in \cite{H2,C}, Huang and the author obtained isomorphisms from the associativity of intertwining operator algebras and from the skew-symmetry isomorphism $\Omega$.
The {\it fusing isomorphism} which we obtained from the associativity of intertwining operator algebras
is a map
\begin{equation}
\F: \quad \frac{\displaystyle \coprod_{a_{1}, a_{2}, a_{3}, a_{4},
a_{5}\in \mathcal{ A}}\mathcal{ V}_{a_{1}a_{5}}^{a_{4}}\otimes
\mathcal{ V}_{a_{2}a_{3}}^{a_{5}}}{\displaystyle \kp}
 \longrightarrow  \frac{\displaystyle \coprod_{a_{1}, a_{2}, a_{3}, a_{4},
a_{5}\in \mathcal{ A}}\mathcal{ V}_{a_{1}a_{2}}^{a_{5}}\otimes
\mathcal{ V}_{a_{5}a_{3}}^{a_{4}}}{\ki}
\end{equation}
determined by linearity and by
\begin{equation}
\F(\Y_1\otimes\Y_2+\kp)=\sum_{a\in \mathcal{ A}}\sum_{i=1}^{\mathcal{ N}_{a_{1}a_{2}}^{a}
\mathcal{N}_{aa_{3}}^{a_{4}}}\mathcal{ Y}^{a}_{3, i}\otimes
\mathcal{ Y}^{a}_{4,i}+\ki
\end{equation}
for $a_1,\cdots,a_5\in\A$, $\Y_1\in \V_{a_1a_5}^{a_4}$ and
$\Y_2\in \V_{a_2a_3}^{a_5}$, where
\begin{equation}
\{\mathcal{ Y}^{a}_{3, i}\in\V_{a_1a_2}^{a},\mathcal{ Y}^{a}_{4,i}\in\V_{aa_3}^{a_4}\mid i=1,\cdots, \mathcal{ N}_{a_{1}a_{2}}^{a}\mathcal{N}_{aa_{3}}^{a_{4}}, a\in\A\}
\end{equation}
is a set of intertwining operators satisfying that for any $w_1,w_2,w_3\in W$ and $w'\in W'$,
\begin{equation}\label{e2:2}
\left.\sum_{a\in \mathcal{ A}}\sum_{i=1}^{\mathcal{ N}_{a_{1}a_{2}}^{a}\mathcal{
N}_{aa_{3}}^{a_{4}}}\langle w', ({\bf I}(\mathcal{ Y}^{a}_{3,i}
\otimes\mathcal{ Y}^{a}_{4, i}))(w_1, w_2, w_3;
x_0,x_2)\rangle_{W}\right|_{x_{0}^{n}=e^{n \log (z_{1}-z_{2})},
x_{2}^{n}=e^{n\log z_{2}}}
\end{equation}
is equal to
\begin{equation}\label{e2:1}
\langle w',
({\bf P}(\mathcal{ Y}_{1}\otimes\mathcal{ Y}_{2}))(w_1, w_2, w_3;
x_1,x_{2})\rangle_{W}\mbar_{x_{1}^{n}=e^{n \log z_{1}},
x_{2}^{n}=e^{n\log z_{2}}}
\end{equation}
on the region
\begin{equation*}
S_1=\{(z_1,z_2)\in \mathbb{C}^2 \mid \re z_1>\re z_2 > \re (z_1-z_2) >0,\ \im z_1 >  \im z_2 > \im (z_1-z_2) > 0\}.
\end{equation*}
It was also proved that
\begin{eqnarray}
\F(a_1,a_2,a_3,a_4): \quad \pi_P\left(\coprod_{a_{5}\in \mathcal{ A}}
\mathcal{ V}_{a_{1}a_{5}}^{a_{4}}\otimes
\mathcal{ V}_{a_{2}a_{3}}^{a_{5}}\right)  &\longrightarrow&
\pi_I\left(\coprod_{a_{5}\in \mathcal{ A}} \mathcal{ V}_{a_{1}a_{2}}^{a_{5}}\otimes
\mathcal{ V}_{a_{5}a_{3}}^{a_{4}}\right)\qquad\nno\\
\Y_1\otimes\Y_2+\kp &\longmapsto& \F(\Y_1\otimes\Y_2+\kp)\qquad
\end{eqnarray}
is an isomorphism for any $a_1,\cdots,a_4\in\A$, where $\mathcal{ Y}_{1}\in \mathcal{
V}_{a_{1}a_{5}}^{a_{4}}$, $\mathcal{ Y}_{2}\in \mathcal{
V}_{a_{2}a_{3}}^{a_{5}}$. These isomorphisms are also called {\it fusing isomorphisms}.
The isomorphisms we obtained from $\Omega$ and its inverse are linear isomorphic maps:
\begin{equation}\label{e3:1}
\tilde{\Omega}^{(1)}, (\widetilde{\Omega^{-1}})^{(1)}:\
\frac{ \displaystyle\coprod_{a_{1}, a_{2},
a_{3}, a_{4}, a_{5}\in \mathcal{ A}}
\mathcal{ V}_{a_{1}a_{2}}^{a_{5}}\otimes
\mathcal{ V}_{a_{5}a_{3}}^{a_{4}}}{\displaystyle\ki}
\longrightarrow \frac{ \displaystyle
\coprod_{a_{1}, a_{2}, a_{3}, a_{4}, a_{5}\in \mathcal{ A}}
\mathcal{ V}_{a_{2}a_{1}}^{a_{5}}\otimes
\mathcal{ V}_{a_{5}a_{3}}^{a_{4}}}{\displaystyle\ki}
\end{equation}
defined by linearity and by
\begin{equation}\label{e3:2}
\tilde{\Omega}^{(1)}(\Y_1\otimes\Y_2+\ki)= \Omega(\Y_1)\otimes\Y_2+\ki,
\end{equation}
\begin{equation}\label{e3:3}
(\widetilde{\Omega^{-1}})^{(1)}(\Y_1\otimes\Y_2+\ki)=\Omega^{-1}(\Y_1)\otimes\Y_2+\ki
 \end{equation}
 for $a_{1}, \dots, a_{5}\in \mathcal{ A}$, $\mathcal{ Y}_{1}\in \mathcal{
V}_{a_{1}a_{2}}^{a_{5}}$, $\mathcal{ Y}_{2}\in \mathcal{
V}_{a_{5}a_{3}}^{a_{4}}$;
\begin{equation}
\tilde{\Omega}^{(2)}, (\widetilde{\Omega^{-1}})^{(2)}: \
\frac{ \displaystyle \coprod_{a_{1}, a_{2}, a_{3}, a_{4}, a_{5}\in \mathcal{ A}}
\mathcal{ V}_{a_{1}a_{5}}^{a_{4}}\otimes \mathcal{ V}_{a_{2}a_{3}}^{a_{5}}}{\displaystyle\kp}
\longrightarrow \frac{ \displaystyle \coprod_{a_{1}, a_{2}, a_{3}, a_{4}, a_{5}\in \mathcal{ A}}
\mathcal{ V}_{a_{2}a_{3}}^{a_{5}}\otimes \mathcal{ V}_{a_{5}a_{1}}^{a_{4}}}{\displaystyle\ki}
\end{equation}
defined by linearity and by
\begin{equation}
\tilde{\Omega}^{(2)}(\Y_1\otimes\Y_2+\kp)=\Y_2\otimes \Omega(\Y_1)+\ki,
\end{equation}
\begin{equation}
(\widetilde{\Omega^{-1}})^{(2)}(\Y_1\otimes\Y_2+\kp)=\Y_2\otimes \Omega^{-1}(\Y_1)+\ki
\end{equation}
 for $a_{1}, \dots, a_{5}\in \mathcal{ A}$, $\mathcal{ Y}_{1}\in \mathcal{
V}_{a_{1}a_{5}}^{a_{4}}$, $\mathcal{ Y}_{2}\in \mathcal{
V}_{a_{2}a_{3}}^{a_{5}}$;
\begin{equation}
\tilde{\Omega}^{(3)}, (\widetilde{\Omega^{-1}})^{(3)}:\
\frac{\displaystyle \coprod_{a_{1}, a_{2},
a_{3}, a_{4}, a_{5}\in \mathcal{ A}}
\mathcal{ V}_{a_{1}a_{2}}^{a_{5}}\otimes \mathcal{ V}_{a_{5}a_{3}}^{a_{4}}}{\displaystyle\ki}
\longrightarrow \frac{\displaystyle \coprod_{a_{1}, a_{2}, a_{3}, a_{4}, a_{5}\in \mathcal{ A}}
\mathcal{ V}_{a_{3}a_{5}}^{a_{4}}\otimes \mathcal{ V}_{a_{1}a_{2}}^{a_{5}}}{\displaystyle\kp}
\end{equation}
defined by linearity and by
\begin{equation}
\tilde{\Omega}^{(3)}(\Y_1\otimes\Y_2+\ki)= \Omega(\Y_2)\otimes\Y_1+\kp,
\end{equation}
\begin{equation}
(\widetilde{\Omega^{-1}})^{(3)}(\Y_1\otimes\Y_2+\ki)=\Omega^{-1}(\Y_2)\otimes\Y_1+\kp
 \end{equation}
 for $a_{1}, \dots, a_{5}\in \mathcal{ A}$, $\mathcal{ Y}_{1}\in \mathcal{
V}_{a_{1}a_{2}}^{a_{5}}$, $\mathcal{ Y}_{2}\in \mathcal{
V}_{a_{5}a_{3}}^{a_{4}}$;
\begin{equation}
\tilde{\Omega}^{(4)}, (\widetilde{\Omega^{-1}})^{(4)}:\
\frac{\displaystyle\coprod_{a_{1}, a_{2}, a_{3}, a_{4},
a_{5}\in \mathcal{ A}}
\mathcal{ V}_{a_{1}a_{5}}^{a_{4}}\otimes \mathcal{ V}_{a_{2}a_{3}}^{a_{5}}}{\displaystyle\kp}
\longrightarrow \frac{\displaystyle \coprod_{a_{1}, a_{2}, a_{3}, a_{4}, a_{5}\in \mathcal{ A}}
\mathcal{ V}_{a_{1}a_{5}}^{a_{4}}\otimes \mathcal{ V}_{a_{3}a_{2}}^{a_{5}}}{\displaystyle\kp}
\end{equation}
defined by linearity and by
\begin{equation}
\tilde{\Omega}^{(4)}(\Y_1\otimes\Y_2+\kp)=\Y_1\otimes \Omega(\Y_2)+\kp,
\end{equation}
\begin{equation}
(\widetilde{\Omega^{-1}})^{(4)}(\Y_1\otimes\Y_2+\kp)=\Y_1\otimes \Omega^{-1}(\Y_2)+\kp
 \end{equation}
 for $a_{1}, \dots, a_{5}\in \mathcal{ A}$, $\mathcal{ Y}_{1}\in \mathcal{
V}_{a_{1}a_{5}}^{a_{4}}$, $\mathcal{ Y}_{2}\in \mathcal{
V}_{a_{2}a_{3}}^{a_{5}}$.
And these isomorphisms have relations:
\begin{equation}\label{e2:21}
(\tilde{\Omega}^{(2)})^{-1}=(\widetilde{\Omega^{-1}})^{(3)},\quad ((\widetilde{\Omega^{-1}})^{(2)})^{-1}=\tilde{\Omega}^{(3)},
\end{equation}
\begin{equation}\label{e2:22}
(\tilde{\Omega}^{(1)})^{-1}=(\widetilde{\Omega^{-1}})^{(1)}, \quad (\tilde{\Omega}^{(4)})^{-1}=(\widetilde{\Omega^{-1}})^{(4)}.
\end{equation}

The above isomorphisms are not independent, we proved the following relations in \cite{C}:
\begin{thm}
The above isomorphisms satisfy the following {\it genus-zero Moore-Seiberg
equations}:
\begin{eqnarray}
\F\circ \tilde{\Omega}^{(3)} \circ \F&=&\tilde{\Omega}^{(1)}\circ
\F\circ \tilde{\Omega}^{(4)},\label{hexagon1}\\
\F\circ (\widetilde{\Omega^{-1}})^{(3)} \circ \F&=&(\widetilde{\Omega^{-1}})^{(1)}\circ
\F\circ (\widetilde{\Omega^{-1}})^{(4)}.\label{hexagon2}
\end{eqnarray}
\end{thm}

Using the fusing isomorphism and the isomorphism $\tilde{\Omega}^{(1)}$,
we deduced a {\it braiding isomorphism}
\begin{eqnarray}
\mathcal{B}=\F^{-1}\circ \tilde{\Omega}^{(1)}\circ \F:\qquad
\qquad\qquad \qquad\qquad \qquad\qquad \qquad \qquad \qquad \qquad \qquad&&\nno\\
 \qquad \qquad\frac{\displaystyle \coprod_{a_{1}, a_{2}, a_{3}, a_{4},
a_{5}\in \mathcal{ A}}\mathcal{ V}_{a_{1}a_{5}}^{a_{4}}\otimes
\mathcal{ V}_{a_{2}a_{3}}^{a_{5}}}{\displaystyle\kp} \longrightarrow
\frac{\displaystyle\coprod_{a_{1}, a_{2}, a_{3}, a_{4},
a_{5}\in \mathcal{ A}}\mathcal{ V}_{a_{2}a_{5}}^{a_{4}}\otimes
\mathcal{ V}_{a_{1}a_{3}}^{a_{5}}}{\displaystyle\kp}.&&
\end{eqnarray}
Moreover, we get an isomorphism
\begin{eqnarray}
\mathcal{B}(a_1,a_2,a_3,a_4): \quad \pi_P\left(\coprod_{a_{5}\in
\mathcal{ A}}\mathcal{ V}_{a_{1}a_{5}}^{a_{4}}\otimes
\mathcal{ V}_{a_{2}a_{3}}^{a_{5}}\right)  &\longrightarrow&
\pi_P\left(\coprod_{a_{5}\in \mathcal{ A}} \mathcal{ V}_{a_{2}a_{5}}^{a_{4}}\otimes
\mathcal{ V}_{a_{1}a_{3}}^{a_{5}}\right)\qquad\nno\\
\Y_1\otimes\Y_2+\kp &\longmapsto& \mathcal{B}(\Y_1\otimes\Y_2+\kp)\qquad
\end{eqnarray}
for any $a_1,\cdots,a_4\in\A$, where $\mathcal{ Y}_{1}\in \mathcal{
V}_{a_{1}a_{5}}^{a_{4}}$, $\mathcal{ Y}_{2}\in \mathcal{
V}_{a_{2}a_{3}}^{a_{5}}$. We also call these isomorphisms {\it braiding isomorphisms}.
\vspace{0.2cm}

Before formulating the Jacobi identity for intertwining operator algebras,
we need to recall the specifics of one more property, which is about certain special multivalued analytic
functions, and were discussed in \cite{H2,C}.

First we consider some simply connected regions in $\mathbb{ C}^{2}$.
Cutting the regions $|z_{1}|>|z_{2}|>0$ and
$|z_{2}|>|z_{1}|>0$ along the
intersections of these regions with $$\{(z_{1}, z_{2})\in \C^{2}\;|\;z_{1}\in [0,+\infty)\}\cup \{(z_{1}, z_{2})\in \C^{2}\;|\;z_{2}\in [0, +\infty)\},$$
we obtain two simply connected regions, which, as in \cite{H2,C}, are denoted by $R_{1}$ and $R_{2}$,
respectively. Also, let $R_{3}$ be the simply connected region obtained
by cutting the region
$|z_{2}|>|z_{1}-z_{2}|>0$ along the
intersection of this region with
$$\{(z_{1}, z_{2})\in \C^{2}\;|\;z_{2}\in [0,+\infty)\}\cup \{(z_{1}, z_{2})\in \C^{2}\;|\;z_{1}-z_{2}\in [0,+\infty)\},$$
and let $R_{4}$ be the simply connected region obtained
by cutting the region
$|z_{1}|>|z_{1}-z_{2}|>0$ along the
intersection of this region with $$\{(z_{1}, z_{2})\in \C^{2}\;|\;z_{1}\in [0,+\infty)\}\cup \{(z_{1}, z_{2})\in \C^{2}\;|\;z_{2}-z_{1}\in [0,+\infty)\}.$$
Then we consider some special multivalued analytic functions on
\begin{equation}
M^{2}=\{(z_{1}, z_{2})\in \mathbb{ C}^{2}\;|\; z_{1}, z_{2}\ne 0,
z_{1}\ne z_{2}\}.
\end{equation}
For $a_{1}, a_{2}, a_{3}, a_{4}\in \mathcal{ A}$, as in \cite{H2,C}, we let $\mathbb{ G}^{a_{1},
a_{2}, a_{3}, a_{4}}$ be the set of multivalued analytic functions on
$M^{2}$
with a choice of a single-valued branch on the region $R_{1}$ satisfying the following property:
Any
branch of $f(z_{1}, z_{2})\in \mathbb{ G}^{a_{1}, a_{2}, a_{3}, a_{4}}$
on the regions $|z_{1}|>|z_{2}|>0$,
$|z_{2}|>|z_{1}|>0$ and $|z_{2}|>|z_{1}-z_{2}|>0$, respectively, can be expanded as
\begin{eqnarray}
&{\displaystyle \sum_{a\in \mathcal{ A}}z_{1}^{h_{a_{4}}-h_{a_{1}}-h_{a}}
z_{2}^{h_{a}-h_{a_{2}}-h_{a_{3}}}F_{a}(z_{1}, z_{2})},&\\
&{\displaystyle \sum_{a\in \mathcal{ A}}z_{2}^{h_{a_{4}}-h_{a_{2}}-h_{a}}
z_{1}^{h_{a}-h_{a_{1}}-h_{a_{3}}}G_{a}(z_{1}, z_{2})}&
\end{eqnarray}
and
\begin{equation}
\sum_{a\in \mathcal{ A}}z_{2}^{h_{a_{4}}-h_{a}-h_{a_{3}}}
(z_{1}-z_{2})^{h_{a}-h_{a_{1}}-h_{a_{2}}}H_{a}(z_{1}, z_{2}),
\end{equation}
respectively, where for $a\in \mathcal{ A}$,
\begin{equation}
F_{a}(z_{1}, z_{2})\in \mathbb{ C}[[z_{2}/z_{1}]][z_{1}, z_{1}^{-1}, z_{2}, z_{2}^{-1}],
\end{equation}
\begin{equation}
G_{a}(z_1,z_2)\in \mathbb{ C}
[[z_{1}/z_{2}]][z_{1}, z_{1}^{-1}, z_{2}, z_{2}^{-1}]
\end{equation}
and
\begin{equation}
H_{a}(z_1,z_2)\in \mathbb{ C}
[[(z_{1}-z_{2})/z_{2}]][z_{2}, z_{2}^{-1}, z_{1}-z_{2},
(z_{1}-z_{2})^{-1}].
\end{equation}
The chosen single-valued branch on $R_{1}$ of an element
of $\mathbb{ G}^{a_{1}, a_{2}, a_{3}, a_{4}}$ is called
the {\it preferred branch on $R_{1}$}. As in \cite{C}, we use the nonempty simply connected regions
\begin{equation*}
S_1=\{(z_1,z_2)\in \mathbb{C}^2 \mid \re z_1>\re z_2 > \re (z_1-z_2) >0,\ \im z_1 >  \im z_2 > \im (z_1-z_2) > 0\}
\end{equation*}
and
\begin{equation*}
S_2=\{(z_1,z_2)\in \mathbb{C}^2 \mid \re z_2>\re z_1 > \re (z_2-z_1) >0,\ \im z_2 >  \im z_1 > \im (z_2-z_1) > 0\}
\end{equation*}
to determine other special branches of an element
of $\mathbb{ G}^{a_{1}, a_{2}, a_{3}, a_{4}}$ on $R_{2}$, $R_{3}$ and $R_{4}$ related to the preferred branch on $R_{1}$. Firstly,
the restriction of the preferred branch on $R_{1}$ of an element of
$\mathbb{ G}^{a_{1}, a_{2}, a_{3}, a_{4}}$
to the region $S_1 \subset R_{1}\cap R_{3}$ gives a single-valued
branch of the element on $R_{3}$, which is then called the
{\it preferred branch on $R_{3}$}. Secondly, the restriction of the preferred branch on $R_{1}$ to the region $S_1 \subset R_{1}\cap R_{4}$ also gives a
single-valued branch of the element on $R_{4}$, which is then called the {\it preferred branch on $R_{4}$}. Moreover, the restriction of the preferred branch on $R_{4}$ to the region $S_2 \subset R_{4}\cap R_{2}$ then gives a
single-valued branch of the element on $R_{2}$ and we
call it the {\it preferred branch on $R_{2}$}.
It was verified in \cite{H2,C} that
$\mathbb{ G}^{a_{1}, a_{2}, a_{3}, a_{4}}$ is a
vector space.

For any element of
$\mathbb{ G}^{a_{1}, a_{2}, a_{3}, a_{4}}$,
the preferred branches of this function on $R_{1}$, $R_{2}$ and $R_{3}$
give  formal series in
\begin{equation}
\coprod_{a\in \mathcal{ A}}x_{1}^{h_{a_{4}}-h_{a_{1}}-h_{a}}
x_{2}^{h_{a}-h_{a_{2}}-h_{a_{3}}}
\mathbb{ C}[[x_{2}/x_{1}]][x_{1}, x_{1}^{-1}, x_{2}, x_{2}^{-1}],
\end{equation}
\begin{equation}
\coprod_{a\in \mathcal{ A}}x_{2}^{h_{a_{4}}-h_{a_{2}}-h_{a}}
x_{1}^{h_{a}-h_{a_{1}}-h_{a_{3}}}
\mathbb{ C}[[x_{1}/x_{2}]][x_{1}, x_{1}^{-1}, x_{2}, x_{2}^{-1}]
\end{equation}
and
\begin{equation}
\coprod_{a\in \mathcal{ A}}x_{2}^{h_{a_{4}}-h_{a}-h_{a_{3}}}
x_{0}^{h_{a}-h_{a_{1}}-h_{a_{2}}}\mathbb{ C}[[x_{0}/x_{2}]][
x_{0}, x_{0}^{-1}, x_{2}, x_{2}^{-1}],
\end{equation}
respectively, which induce linear maps
\begin{eqnarray}
\iota_{12}: \mathbb{ G}^{a_{1}, a_{2}, a_{3}, a_{4}}&\to&
\coprod_{a\in \mathcal{ A}}x_{1}^{h_{a_{4}}-h_{a_{1}}-h_{a}}
x_{2}^{h_{a}-h_{a_{2}}-h_{a_{3}}}\mathbb{ C}[[x_{2}/x_{1}]][x_{1}, x_{1}^{-1}, x_{2},
x_{2}^{-1}],\qquad\\
\iota_{21}: \mathbb{ G}^{a_{1}, a_{2}, a_{3}, a_{4}}&\to&
\coprod_{a\in \mathcal{ A}}x_{2}^{h_{a_{4}}-h_{a_{2}}-h_{a}}
x_{1}^{h_{a}-h_{a_{1}}-h_{a_{3}}}\mathbb{ C}[[x_{1}/x_{2}]][x_{1}, x_{1}^{-1}, x_{2},
x_{2}^{-1}],\qquad\\
\iota_{20}: \mathbb{ G}^{a_{1}, a_{2}, a_{3}, a_{4}}&\to&
\coprod_{a\in \mathcal{ A}}x_{2}^{h_{a_{4}}-h_{a}-h_{a_{3}}}
x_{0}^{h_{a}-h_{a_{1}}-h_{a_{2}}}\mathbb{ C}[[x_{0}/x_{2}]][ x_{0},
x_{0}^{-1}, x_{2}, x_{2}^{-1}]\qquad
\end{eqnarray}
generalizing $\iota_{12}$, $\iota_{21}$ and $\iota_{20}$ discussed
at the beginning of this section. These maps are
injective because analytic extensions are unique.

For $a_{1}, a_{2}, a_{3}, a_{4}\in \mathcal{ A}$,
$\mathbb{ G}^{a_{1}, a_{2}, a_{3}, a_{4}}$ is a module over the ring
$\mathbb{ C}[x_{1}, x_{1}^{-1}, x_{2}, x_{2}^{-1}, (x_{1}-x_{2})^{-1}]$.
Huang \cite{H2} proved the following lemma:

\begin{lemma}\label{free}
For any $a_{1}, a_{2}, a_{3}, a_{4}\in \mathcal{ A}$, the module
$\mathbb{ G}^{a_{1}, a_{2}, a_{3}, a_{4}}$ is free.
\end{lemma}

\begin{rema}
In the following theorem for Jacobi identity and for the rest of the paper,
we fix a basis
$\{e^{a_{1}, a_{2}, a_{3}, a_{4}}_{\alpha}\}_{\alpha\in
\mathbb{ A}(a_{1}, a_{2}, a_{3}, a_{4})}$ of the free module
$\mathbb{ G}^{a_{1}, a_{2}, a_{3}, a_{4}}$ over the ring
$\mathbb{ C}[x_{1}, x_{1}^{-1}, x_{2}, x_{2}^{-1}, (x_{1}-x_{2})^{-1}]$ for any
$a_{1}, a_{2}, a_{3}, a_{4}\in \mathcal{ A}$, where $\mathbb{ A}(a_{1}, a_{2}, a_{3}, a_{4})$ is the index set of the basis.
\end{rema}

Now we give the Jacobi identity derived in \cite{H2}:

\begin{thm}[{\bf Jacobi identity}]\label{jacobi}
For any $a_{1}, a_{2}, a_{3}, a_{4}\in \mathcal{ A}$, there exist
linear maps
\begin{eqnarray}
\lefteqn{f^{a_{1}, a_{2}, a_{3},
a_{4}}_{\alpha}: W^{a_{1}}\otimes W^{a_{2}}\otimes W^{a_{3}}\otimes
\pi_P(\coprod_{a_{5}\in \mathcal{ A}}
\mathcal{ V}_{a_{1}a_{5}}^{a_{4}}\otimes
\mathcal{ V}_{a_{2}a_{3}}^{a_{5}})}\nno\\
&&\hspace{8em}\to  W^{a_{4}}[[x_{2}/x_{1}]][x_{1}, x_{1}^{-1}, x_{2},
x_{2}^{-1}]\nno\\
&&\hspace{3em} w_{(a_{1})}\otimes w_{(a_{2})}\otimes w_{(a_{3})}\otimes [\mathcal{Z}]_P \nno\\
&&\hspace{8em} \mapsto f^{a_{1}, a_{2}, a_{3},
a_{4}}_{\alpha}(w_{(a_{1})},
w_{(a_{2})}, w_{(a_{3})}, [\mathcal{Z}]_P; x_{1}, x_{2}),\label{e1:6}
\end{eqnarray}
\begin{eqnarray}
\lefteqn{g^{a_{1}, a_{2}, a_{3},
a_{4}}_{\alpha}: W^{a_{1}}\otimes W^{a_{2}}\otimes W^{a_{3}}\otimes
\pi_P(\coprod_{a_{5}\in \mathcal{ A}}
\mathcal{ V}_{a_{2}a_{5}}^{a_{4}}\otimes
\mathcal{ V}_{a_{1}a_{3}}^{a_{5}})}\nno\\
&&\hspace{8em} \to  W^{a_{4}}[[x_{1}/x_{2}]][x_{1}, x_{1}^{-1}, x_{2},
x_{2}^{-1}]\nno\\
&&\hspace{3em} w_{(a_{1})}\otimes w_{(a_{2})}\otimes w_{(a_{3})}\otimes [\mathcal{Z}]_P \nno\\
&&\hspace{8em}\mapsto g^{a_{1}, a_{2}, a_{3},
a_{4}}_{\alpha}(w_{(a_{1})},
w_{(a_{2})}, w_{(a_{3})}, [\mathcal{Z}]_P; x_{1}, x_{2})\label{e1:7}
\end{eqnarray}
and
\begin{eqnarray}\label{1.-3}
\lefteqn{h^{a_{1}, a_{2}, a_{3},
a_{4}}_{\alpha}: W^{a_{1}}\otimes W^{a_{2}}\otimes W^{a_{3}}\otimes
\pi_I(\coprod_{a_{5}\in \mathcal{ A}}
\mathcal{ V}_{a_{1}a_{2}}^{a_{5}}\otimes
\mathcal{ V}_{a_{5}a_{3}}^{a_{4}})}\nno\\
&&\hspace{8em}\to W^{a_{4}}[[x_{0}/x_{2}]][x_{0},
x_{0}^{-1}, x_{2}, x_{2}^{-1}]\nno\\
&&\hspace{3em} w_{(a_{1})}\otimes w_{(a_{2})}\otimes w_{(a_{3})}\otimes [\mathcal{Z}]_I \nno\\
&&\hspace{8em} \mapsto h^{a_{1}, a_{2}, a_{3},
a_{4}}_{\alpha}(w_{(a_{1})},
w_{(a_{2})}, w_{(a_{3})}, [\mathcal{Z}]_I; x_{0}, x_{2})\label{e1:8}
\end{eqnarray}
for $\alpha\in \mathbb{ A}(a_{1}, a_{2}, a_{3}, a_{4})$,
such that  for
any
$w_{(a_{1})}\in W^{a_{1}}$, $w_{(a_{2})}\in W^{a_{2}}$, $w_{(a_{3})}\in
W^{a_{3}}$, and any
\begin{equation}
\mathcal{ Z}\in \coprod_{a_{5}\in \mathcal{ A}}
\mathcal{ V}_{a_{1}a_{5}}^{a_{4}}\otimes
\mathcal{ V}_{a_{2}a_{3}}^{a_{5}}\subset \coprod_{a_{1}, a_{2}, a_{3}, a_{4},
a_{5}\in \mathcal{ A}}\mathcal{ V}_{a_{1}a_{5}}^{a_{4}}\otimes
\mathcal{ V}_{a_{2}a_{3}}^{a_{5}},
\end{equation}
only finitely many of
\begin{equation}\label{e4:1}
f^{a_{1}, a_{2}, a_{3},
a_{4}}_{\alpha}(w_{(a_{1})},
w_{(a_{2})}, w_{(a_{3})}, [\mathcal{Z}]_P; x_{1}, x_{2}),
\end{equation}
\begin{equation}\label{e4:2}
g^{a_{1}, a_{2}, a_{3},
a_{4}}_{\alpha}(w_{(a_{1})},
w_{(a_{2})}, w_{(a_{3})}, \mathcal{B}([\mathcal{Z}]_P); x_{1}, x_{2}),
\end{equation}
and
\begin{equation}\label{e4:3}
h^{a_{1}, a_{2}, a_{3},
a_{4}}_{\alpha}(w_{(a_{1})},
w_{(a_{2})}, w_{(a_{3})}, \mathcal{ F}([\mathcal{Z}]_P); x_{0}, x_{2}),
\end{equation}
$\alpha\in \mathbb{ A}(a_{1}, a_{2}, a_{3}, a_{4})$, are nonzero,
\begin{eqnarray}
\lefteqn{(\tilde{\bf P}([\mathcal{Z}]_P))(w_{(a_{1})},
w_{(a_{2})}, w_{(a_{3})}; x_{1}, x_{2})}\nno\\
&&=\sum_{\alpha\in \mathbb{ A}(a_{1}, a_{2}, a_{3}, a_{4})}f^{a_{1}, a_{2}, a_{3},
a_{4}}_{\alpha}(w_{(a_{1})},
w_{(a_{2})}, w_{(a_{3})}, [\mathcal{Z}]_P; x_{1}, x_{2}) \iota_{12}\left(e^{a_{1},
a_{2}, a_{3}, a_{4}}_{\alpha}\right),\label{e1:9}
\end{eqnarray}
\begin{eqnarray}
\lefteqn{(\tilde{\bf P}(\mathcal{B}([\mathcal{Z}]_P)))(w_{(a_{2})},
w_{(a_{1})}, w_{(a_{3})}; x_{2}, x_{1})}\nno\\
&&=\sum_{\alpha\in \mathbb{ A}(a_{1}, a_{2}, a_{3}, a_{4})}g^{a_{1}, a_{2}, a_{3},
a_{4}}_{\alpha}(w_{(a_{1})},
w_{(a_{2})}, w_{(a_{3})}, \mathcal{B}([\mathcal{Z}]_P); x_{1}, x_{2}) \iota_{21}\left(e^{a_{1},
a_{2}, a_{3}, a_{4}}_{\alpha}\right),\quad \label{e1:10}
\end{eqnarray}
\begin{eqnarray}
\lefteqn{
(\tilde{\bf I}(\F([\mathcal{Z}]_P)))(w_{(a_{1})}, w_{(a_{2})},
w_{(a_{3})}; x_{0}, x_{2})}\nno\\
&&=\sum_{\alpha\in \mathbb{ A}(a_{1}, a_{2}, a_{3}, a_{4})}h^{a_{1}, a_{2}, a_{3},
a_{4}}_{\alpha}(w_{(a_{1})},
w_{(a_{2})}, w_{(a_{3})}, \F([\mathcal{Z}]_P); x_{0}, x_{2}) \iota_{20}
(e^{a_{1},
a_{2}, a_{3}, a_{4}}_{\alpha}),\quad \label{e1:11}
\end{eqnarray}
and the following {\it Jacobi identity} holds:
\begin{eqnarray}
\lefteqn{x_{0}^{-1}
\delta\left(\frac{x_{1}-x_{2}}{x_{0}}\right)
f^{a_{1}, a_{2}, a_{3},
a_{4}}_{\alpha}(w_{(a_{1})},
w_{(a_{2})}, w_{(a_{3})}, [\mathcal{Z}]_P; x_{1}, x_{2})}\nno\\
&&\quad -x_{0}^{-1}\delta\left(\frac{x_{2}-x_{1}}{-x_{0}}\right)
g^{a_1, a_2, a_{3},
a_{4}}_{\alpha}(w_{(a_1)},
w_{(a_2)}, w_{(a_{3})}, \mathcal{ B}([\mathcal{Z}]_P); x_1, x_2)
\nno\\
&&=x_{2}^{-1}\delta\left(\frac{x_{1}-x_{0}}{x_{2}}\right)
h^{a_{1}, a_{2}, a_{3},
a_{4}}_{\alpha}(w_{(a_{1})},
w_{(a_{2})}, w_{(a_{3})}, \mathcal{ F}([\mathcal{Z}]_P); x_{0}, x_{2})\label{e1:12}
\end{eqnarray}
for $\alpha\in \mathbb{ A}(a_{1}, a_{2}, a_{3}, a_{4})$.
\end{thm}

\section{$S_3$-symmetry of the Jacobi identity}

In this section, we formulate and prove a symmetric property of the Jacobi identity for intertwining operator algebras under the symmetric group $S_3$, which is the main result of this paper. Here is the precise statement of the main result:

\begin{thm}\label{t1}
In the presence of the axioms for an intertwining operator algebra except for the associativity property, we assume that there exists an isomorphism
\begin{equation}\label{s5}
\F: \quad \frac{\coprod_{a_{1}, a_{2}, a_{3}, a_{4},
a_{5}\in \mathcal{ A}}\mathcal{ V}_{a_{1}a_{5}}^{a_{4}}\otimes
\mathcal{ V}_{a_{2}a_{3}}^{a_{5}}}{\kp}  \longrightarrow  \frac{\coprod_{a_{1}, a_{2}, a_{3}, a_{4},
a_{5}\in \mathcal{ A}}\mathcal{ V}_{a_{1}a_{2}}^{a_{5}}\otimes
\mathcal{ V}_{a_{5}a_{3}}^{a_{4}}}{\ki}
\end{equation}
satisfying
\begin{equation}\label{a18}
\F\left(\pi_P\left(\coprod_{a_{5}\in \mathcal{ A}}\mathcal{ V}_{a_{1}a_{5}}^{a_{4}}\otimes
\mathcal{ V}_{a_{2}a_{3}}^{a_{5}}\right)\right)=\pi_I \left(\coprod_{a_{5}\in \mathcal{ A}} \mathcal{ V}_{a_{1}a_{2}}^{a_{5}}\otimes
\mathcal{ V}_{a_{5}a_{3}}^{a_{4}}\right)
\end{equation}
for any $a_1,\cdots,a_4\in\A$, and that the Moore-Seiberg equations (\ref{hexagon1}) and
(\ref{hexagon2}) hold,
then the Jacobi identity for the ordered triple
$(\tilde{w}_{(a_{1})}, \tilde{w}_{(a_{2})}, \tilde{w}_{(a_{3})})\in \coprod_{i=1}^{3}W^{a_{i}}$
implies the Jacobi identity for
the triple
$(\tilde{w}_{(a_{\tau(1)})}, \tilde{w}_{(a_{\tau(2)})}, \tilde{w}_{(a_{\tau(3)})})
\in \coprod_{i=1}^{3}W^{a_{\tau(i)}}$ for any $\tau\in S_{3}$.

\end{thm}

\begin{rema}\label{r1}
In the above theorem, since there's no associativity in the assumptions, we have no fusing isomorphism. The assumption that the Moore-Seiberg equations (\ref{hexagon1}) and
(\ref{hexagon2}) hold in fact means that they hold with the fusing isomorphism replaced by the given isomorphism $\F$ in (\ref{s5}). Moreover, from $\F$ in (\ref{s5}) and the isomorphism $\tilde{\Omega}^{(1)}$ we obtain an isomorphism
\begin{eqnarray}
\mathcal{B}=\F^{-1}\circ \tilde{\Omega}^{(1)}\circ \F:\qquad
\qquad\qquad \qquad\qquad \qquad\qquad \qquad \qquad \qquad \qquad \qquad&&\nno\\
 \qquad \qquad\frac{\displaystyle \coprod_{a_{1}, a_{2}, a_{3}, a_{4},
a_{5}\in \mathcal{ A}}\mathcal{ V}_{a_{1}a_{5}}^{a_{4}}\otimes
\mathcal{ V}_{a_{2}a_{3}}^{a_{5}}}{\displaystyle\kp} \longrightarrow
\frac{\displaystyle\coprod_{a_{1}, a_{2}, a_{3}, a_{4},
a_{5}\in \mathcal{ A}}\mathcal{ V}_{a_{2}a_{5}}^{a_{4}}\otimes
\mathcal{ V}_{a_{1}a_{3}}^{a_{5}}}{\displaystyle\kp}.&&\label{s6}
\end{eqnarray}
And we also assume that the isomorphisms $\F$ and $\mathcal{B}$ involved in the formulation of the Jacobi identity are replaced by the given isomorphism in (\ref{s5}) and the deduced isomorphism in  (\ref{s6}), respectively.
\end{rema}
\vspace{0.1cm}

The rest of this section is devoted to proving the $S_3$-symmetry of the Jacobi identity.
We achieve this goal by establishing three results that lead to Theorem \ref{t1}.

\begin{prop}\label{l1}
In the presence of the axioms for an intertwining operator algebra except for the associativity property, we assume that there exists an isomorphism
\begin{equation}
\F: \quad \frac{\coprod_{a_{1}, a_{2}, a_{3}, a_{4},
a_{5}\in \mathcal{ A}}\mathcal{ V}_{a_{1}a_{5}}^{a_{4}}\otimes
\mathcal{ V}_{a_{2}a_{3}}^{a_{5}}}{\kp}  \longrightarrow  \frac{\coprod_{a_{1}, a_{2}, a_{3}, a_{4},
a_{5}\in \mathcal{ A}}\mathcal{ V}_{a_{1}a_{2}}^{a_{5}}\otimes
\mathcal{ V}_{a_{5}a_{3}}^{a_{4}}}{\ki}
\end{equation}
satisfying
\begin{equation}
\F\left(\pi_P\left(\coprod_{a_{5}\in \mathcal{ A}}\mathcal{ V}_{a_{1}a_{5}}^{a_{4}}\otimes
\mathcal{ V}_{a_{2}a_{3}}^{a_{5}}\right)\right)=\pi_I \left(\coprod_{a_{5}\in \mathcal{ A}} \mathcal{ V}_{a_{1}a_{2}}^{a_{5}}\otimes
\mathcal{ V}_{a_{5}a_{3}}^{a_{4}}\right)
\end{equation}
for any $a_1,\cdots,a_4\in\A$, and that the Jacobi identity for the ordered triple
$(\tilde{w}_{(a_{1})}, \tilde{w}_{(a_{2})}, \tilde{w}_{(a_{3})})\in \coprod_{i=1}^{3}W^{a_{i}}$ holds,
then
for any $a_{4}\in \mathcal{A}$, $w_{(a_4)}'\in (W^{a_4})'$ and $\mathcal{ Z}\in \coprod_{ a_{5}\in \mathcal{ A}}
\mathcal{ V}_{a_{1}a_{5}}^{a_{4}}\otimes
\mathcal{ V}_{a_{2}a_{3}}^{a_{5}}$, there exists a multivalued analytic function
\begin{equation}\label{s7}
\Phi(w'_{(a_{4})}, \tilde{w}_{(a_{1})},
\tilde{w}_{(a_{2})}, \tilde{w}_{(a_{3})}, [\mathcal{Z}]_P;z_1,z_2)\in \mathbb{ G}^{a_{1}, a_{2}, a_{3}, a_{4}}
\end{equation}
such that
\begin{equation}\label{s8}
\langle w'_{(a_{4})}, (\tilde{\bf P}([\mathcal{Z}]_P))(\tilde{w}_{(a_{1})},
\tilde{w}_{(a_{2})}, \tilde{w}_{(a_{3})}; x_{1}, x_{2})\rangle_{W^{a_4}}
\mbar_{x_{1}^{n}=e^{n \log z_1},
x_{2}^{n}=e^{n\log z_2}},
\end{equation}
\begin{equation}\label{s9}
\langle w_{(a_4)}', (\tilde{\bf P}(\mathcal{B}([\mathcal{Z}]_P)))(\tilde{w}_{(a_{2})},
\tilde{w}_{(a_{1})}, \tilde{w}_{(a_{3})}; x_{2}, x_{1})\rangle_{W^{a_4}}\mbar_{x_{1}^{n}=e^{n \log z_1},
x_{2}^{n}=e^{n\log z_2}},
\end{equation}
\begin{equation}\label{s10}
\langle w_{(a_4)}', (\tilde{\bf I}(\mathcal{ F}([\mathcal{Z}]_P)))(\tilde{w}_{(a_{1})}, \tilde{w}_{(a_{2})},
\tilde{w}_{(a_{3})}; x_{0}, x_{2})\rangle_{W^{a_4}}\mbar_{x_{0}^{n}=e^{n \log (z_1-z_2)},
x_{2}^{n}=e^{n\log z_2}}
\end{equation}
and
\begin{equation}\label{s11}
\langle w_{(a_4)}', (\tilde{\bf I}(\tilde{\Omega}^{(1)}(\mathcal{F}([\mathcal{Z}]_P)))) (\tilde{w}_{(a_{2})}, \tilde{w}_{(a_{1})},
\tilde{w}_{(a_{3})}; x_{0}, x_1)\rangle_{W^{a_4}}\mbar_{x_{0}^{n}=e^{n \log (z_2-z_1)},
x_{1}^{n}=e^{n\log z_1}}
\end{equation}
are its preferred branches on $R_1$, $R_2$, $R_3$ and $R_4$, respectively. Moreover,
\begin{eqnarray}
\lefteqn{ \langle w_{(a_4)}', (\tilde{\bf I}(\mathcal{ F}([\mathcal{Z}]_P)))(\tilde{w}_{(a_{1})}, \tilde{w}_{(a_{2})},
\tilde{w}_{(a_{3})}; x_{0}, x_{2})\rangle_{W^{a_4}}\mbar_{x_{0}^{n}=e^{n \log (z_1-z_2)},
x_{2}^{n}=e^{n\log z_2}}}\nno\\
&&= \langle w'_{(a_{4})}, (\tilde{\bf P}([\mathcal{Z}]_P))(\tilde{w}_{(a_{1})},
\tilde{w}_{(a_{2})}, \tilde{w}_{(a_{3})}; x_{1}, x_{2})\rangle_{W^{a_4}}
\mbar_{x_{1}^{n}=e^{n \log z_1},
x_{2}^{n}=e^{n\log z_2}} \label{s12}
\end{eqnarray}
on the region
\begin{equation*}
S_1=\{(z_1,z_2)\in \mathbb{C}^2 \mid \re z_1>\re z_2 > \re (z_1-z_2) >0,\ \im z_1 >  \im z_2 > \im (z_1-z_2) > 0\},
\end{equation*}
and
\begin{eqnarray}
\lefteqn{\langle w_{(a_4)}', (\tilde{\bf I}(\tilde{\Omega}^{(1)}(\mathcal{F}([\mathcal{Z}]_P)))) (\tilde{w}_{(a_{2})}, \tilde{w}_{(a_{1})},
\tilde{w}_{(a_{3})}; x_{0}, x_1)\rangle_{W^{a_4}}\mbar_{x_{0}^{n}=e^{n \log (z_2-z_1)},
x_{1}^{n}=e^{n\log z_1}}}\nno\\
&& =\langle w_{(a_4)}', (\tilde{\bf P}(\mathcal{B}([\mathcal{Z}]_P)))(\tilde{w}_{(a_{2})},
\tilde{w}_{(a_{1})}, \tilde{w}_{(a_{3})}; x_{2}, x_{1})\rangle_{W^{a_4}}\mbar_{x_{1}^{n}=e^{n \log z_1},
x_{2}^{n}=e^{n\log z_2}}\label{s13}
\end{eqnarray}
on the region
\begin{equation*}
S_2=\{(z_1,z_2)\in \mathbb{C}^2 \mid \re z_2>\re z_1 > \re (z_2-z_1) >0,\ \im z_2 >  \im z_1 > \im (z_2-z_1) > 0\}.
\end{equation*}

\end{prop}
\vspace{0.1cm}

\begin{rema}
In the above proposition, the formulations involving the isomorphisms $\F$ and $\mathcal{B}$ have the same assumptions as we discussed in Remark \ref{r1} below Theorem \ref{t1}.
\end{rema}

\begin{rema}
In \cite{C}, we proved that for intertwining operator algebras, the generalized rationality, commutativity and associativity follow from the Jacobi identity. And its proof involves only one ordered triple $(\tilde{w}_{(a_{1})}, \tilde{w}_{(a_{2})}, \tilde{w}_{(a_{3})})\in \coprod_{i=1}^{3}W^{a_{i}}$ for both commutativity, associativity and the Jacobi identity. Minus the skew-symmetry condition, the above proposition becomes a one-ordered-triple version of Theorem 3.3 in \cite{C}. However, for the sake of proving Theorem \ref{t1}, we add the extra skew-symmetry condition in the above proposition to obtain the preferred branch (\ref{s11}) on $R_4$ and the analytic extension relation (\ref{s13}).
\end{rema}
\vspace{0.1cm}

{\it Proof of Proposition \ref{l1}.}\hspace{2ex}
Since the Jacobi identity holds for the ordered triple \\
$(\tilde{w}_{(a_{1})}, \tilde{w}_{(a_{2})}, \tilde{w}_{(a_{3})})\in \coprod_{i=1}^{3}W^{a_{i}}$, then
for any $a_{4}\in \mathcal{A}$, there exist
linear maps
\begin{eqnarray}
\lefteqn{f^{a_{1}, a_{2}, a_{3},
a_{4}}_{\alpha}: W^{a_{1}}\otimes W^{a_{2}}\otimes W^{a_{3}}\otimes
\pi_P(\coprod_{a_{5}\in \mathcal{ A}}
\mathcal{ V}_{a_{1}a_{5}}^{a_{4}}\otimes
\mathcal{ V}_{a_{2}a_{3}}^{a_{5}})}\nno\\
&&\hspace{8em}\to  W^{a_{4}}[[x_{2}/x_{1}]][x_{1}, x_{1}^{-1}, x_{2},
x_{2}^{-1}]\nno\\
&&\hspace{3em} w_{(a_{1})}\otimes w_{(a_{2})}\otimes w_{(a_{3})}\otimes [\mathcal{Z}]_P \nno\\
&&\hspace{8em} \mapsto f^{a_{1}, a_{2}, a_{3},
a_{4}}_{\alpha}(w_{(a_{1})},
w_{(a_{2})}, w_{(a_{3})}, [\mathcal{Z}]_P; x_{1}, x_{2}),\label{s1}
\end{eqnarray}
\begin{eqnarray}
\lefteqn{g^{a_{1}, a_{2}, a_{3},
a_{4}}_{\alpha}: W^{a_{1}}\otimes W^{a_{2}}\otimes W^{a_{3}}\otimes
\pi_P(\coprod_{a_{5}\in \mathcal{ A}}
\mathcal{ V}_{a_{2}a_{5}}^{a_{4}}\otimes
\mathcal{ V}_{a_{1}a_{3}}^{a_{5}})}\nno\\
&&\hspace{8em} \to  W^{a_{4}}[[x_{1}/x_{2}]][x_{1}, x_{1}^{-1}, x_{2},
x_{2}^{-1}]\nno\\
&&\hspace{3em} w_{(a_{1})}\otimes w_{(a_{2})}\otimes w_{(a_{3})}\otimes [\mathcal{Z}]_P \nno\\
&&\hspace{8em}\mapsto g^{a_{1}, a_{2}, a_{3},
a_{4}}_{\alpha}(w_{(a_{1})},
w_{(a_{2})}, w_{(a_{3})}, [\mathcal{Z}]_P; x_{1}, x_{2})\label{s2}
\end{eqnarray}
and
\begin{eqnarray}
\lefteqn{h^{a_{1}, a_{2}, a_{3},
a_{4}}_{\alpha}: W^{a_{1}}\otimes W^{a_{2}}\otimes W^{a_{3}}\otimes
\pi_I(\coprod_{a_{5}\in \mathcal{ A}}
\mathcal{ V}_{a_{1}a_{2}}^{a_{5}}\otimes
\mathcal{ V}_{a_{5}a_{3}}^{a_{4}})}\nno\\
&&\hspace{8em}\to W^{a_{4}}[[x_{0}/x_{2}]][x_{0},
x_{0}^{-1}, x_{2}, x_{2}^{-1}]\nno\\
&&\hspace{3em} w_{(a_{1})}\otimes w_{(a_{2})}\otimes w_{(a_{3})}\otimes [\mathcal{Z}]_I \nno\\
&&\hspace{8em} \mapsto h^{a_{1}, a_{2}, a_{3},
a_{4}}_{\alpha}(w_{(a_{1})},
w_{(a_{2})}, w_{(a_{3})}, [\mathcal{Z}]_I; x_{0}, x_{2})\label{s3}
\end{eqnarray}
for $\alpha\in \mathbb{ A}(a_{1}, a_{2}, a_{3}, a_{4})$,
such that  for
any
\begin{equation}
\mathcal{ Z}\in \coprod_{ a_{5}\in \mathcal{ A}}
\mathcal{ V}_{a_{1}a_{5}}^{a_{4}}\otimes
\mathcal{ V}_{a_{2}a_{3}}^{a_{5}}\subset \coprod_{a_{1}, a_{2}, a_{3}, a_{4},
a_{5}\in \mathcal{ A}}\mathcal{ V}_{a_{1}a_{5}}^{a_{4}}\otimes
\mathcal{ V}_{a_{2}a_{3}}^{a_{5}},
\end{equation}
only finitely many of
\begin{equation}\label{a10}
f^{a_{1}, a_{2}, a_{3},
a_{4}}_{\alpha}(\tilde{w}_{(a_{1})},
\tilde{w}_{(a_{2})}, \tilde{w}_{(a_{3})}, [\mathcal{Z}]_P; x_{1}, x_{2}),
\end{equation}
\begin{equation}
g^{a_{1}, a_{2}, a_{3},
a_{4}}_{\alpha}(\tilde{w}_{(a_{1})},
\tilde{w}_{(a_{2})}, \tilde{w}_{(a_{3})}, \mathcal{B}([\mathcal{Z}]_P); x_{1}, x_{2}),
\end{equation}
and
\begin{equation}
h^{a_{1}, a_{2}, a_{3},
a_{4}}_{\alpha}(\tilde{w}_{(a_{1})},
\tilde{w}_{(a_{2})}, \tilde{w}_{(a_{3})}, \mathcal{ F}([\mathcal{Z}]_P); x_{0}, x_{2}),
\end{equation}
$\alpha\in \mathbb{ A}(a_{1}, a_{2}, a_{3}, a_{4})$, are nonzero,
\begin{eqnarray}
\lefteqn{(\tilde{\bf P}([\mathcal{Z}]_P))(\tilde{w}_{(a_{1})},
\tilde{w}_{(a_{2})}, \tilde{w}_{(a_{3})}; x_{1}, x_{2})}\nno\\
&&=\sum_{\alpha\in \mathbb{ A}(a_{1}, a_{2}, a_{3}, a_{4})}f^{a_{1}, a_{2}, a_{3},
a_{4}}_{\alpha}(\tilde{w}_{(a_{1})},
\tilde{w}_{(a_{2})}, \tilde{w}_{(a_{3})}, [\mathcal{Z}]_P; x_{1}, x_{2}) \iota_{12}\left(e^{a_{1},
a_{2}, a_{3}, a_{4}}_{\alpha}\right),
\end{eqnarray}
\begin{eqnarray}
\lefteqn{(\tilde{\bf P}(\mathcal{B}([\mathcal{Z}]_P)))(\tilde{w}_{(a_{2})},
\tilde{w}_{(a_{1})}, \tilde{w}_{(a_{3})}; x_{2}, x_{1})}\nno\\
&&=\sum_{\alpha\in \mathbb{ A}(a_{1}, a_{2}, a_{3}, a_{4})}g^{a_{1}, a_{2}, a_{3},
a_{4}}_{\alpha}(\tilde{w}_{(a_{1})},
\tilde{w}_{(a_{2})}, \tilde{w}_{(a_{3})}, \mathcal{B}([\mathcal{Z}]_P); x_{1}, x_{2}) \iota_{21}\left(e^{a_{1},
a_{2}, a_{3}, a_{4}}_{\alpha}\right),\quad
\end{eqnarray}
\begin{eqnarray}
\lefteqn{(\tilde{\bf I}(\mathcal{ F}([\mathcal{Z}]_P)))(\tilde{w}_{(a_{1})}, \tilde{w}_{(a_{2})},
\tilde{w}_{(a_{3})}; x_{0}, x_{2})}\nno\\
&&=\sum_{\alpha\in \mathbb{ A}(a_{1}, a_{2}, a_{3}, a_{4})}h^{a_{1}, a_{2}, a_{3},
a_{4}}_{\alpha}(\tilde{w}_{(a_{1})},
\tilde{w}_{(a_{2})}, \tilde{w}_{(a_{3})}, \mathcal{ F}([\mathcal{Z}]_P); x_{0}, x_{2}) \iota_{20}
(e^{a_{1},
a_{2}, a_{3}, a_{4}}_{\alpha}),\quad
\end{eqnarray}
and the following {\it Jacobi identity} holds:
\begin{eqnarray}
\lefteqn{x_{0}^{-1}
\delta\left(\frac{x_{1}-x_{2}}{x_{0}}\right)
f^{a_{1}, a_{2}, a_{3},
a_{4}}_{\alpha}(\tilde{w}_{(a_{1})},
\tilde{w}_{(a_{2})}, \tilde{w}_{(a_{3})}, [\mathcal{Z}]_P; x_{1}, x_{2})}\nno\\
&&\quad -x_{0}^{-1}\delta\left(\frac{x_{2}-x_{1}}{-x_{0}}\right)
g^{a_1, a_2, a_{3},
a_{4}}_{\alpha}(\tilde{w}_{(a_1)},
\tilde{w}_{(a_2)}, \tilde{w}_{(a_{3})}, \mathcal{ B}([\mathcal{Z}]_P); x_1, x_2)
\nno\\
&&=x_{2}^{-1}\delta\left(\frac{x_{1}-x_{0}}{x_{2}}\right)
h^{a_{1}, a_{2}, a_{3},
a_{4}}_{\alpha}(\tilde{w}_{(a_{1})},
\tilde{w}_{(a_{2})}, \tilde{w}_{(a_{3})}, \mathcal{ F}([\mathcal{Z}]_P); x_{0}, x_{2})\label{s4}
\end{eqnarray}
for $\alpha\in \mathbb{ A}(a_{1}, a_{2}, a_{3}, a_{4})$.

In analogy with the proof of Theorem 3.3 in \cite{C},
we can obtain that, for any $a_{4}\in \mathcal{A}$, $\alpha\in \mathbb{ A}(a_{1}, a_{2}, a_{3}, a_{4})$, there exists linear map
\begin{eqnarray}
\lefteqn{F_{\alpha}: \quad (W^{a_4})'\otimes W^{a_{1}}\otimes W^{a_{2}}
\otimes W^{a_{3}}\otimes
\pi_P(\coprod_{a_{5}\in \mathcal{ A}}
\mathcal{ V}_{a_{1}a_{5}}^{a_{4}}\otimes
\mathcal{ V}_{a_{2}a_{3}}^{a_{5}})}\nno\\
&&\hspace{8em} \longrightarrow  \mathbb{ C}[x_{1}, x_{1}^{-1}, x_{2},
x_{2}^{-1}, (x_{1}-x_{2})^{-1}]
\end{eqnarray}
with $F_{\alpha}(w'_{(a_{4})}, \tilde{w}_{(a_{1})},
\tilde{w}_{(a_{2})}, \tilde{w}_{(a_{3})}, [\mathcal{Z}]_P; x_{1}, x_{2})$ denoted by the image of
$$w'_{(a_{4})}\otimes \tilde{w}_{(a_{1})}\otimes
\tilde{w}_{(a_{2})}\otimes  \tilde{w}_{(a_{3})}\otimes [\mathcal{Z}]_P$$
under $F_{\alpha}$, such that
\begin{eqnarray}
\lefteqn{\langle w'_{(a_{4})}, f^{a_{1}, a_{2}, a_{3},a_{4}}_{\alpha}(\tilde{w}_{(a_{1})},
\tilde{w}_{(a_{2})}, \tilde{w}_{(a_{3})}, [\mathcal{Z}]_P; x_{1}, x_{2})\rangle_{W^{a_4}}}\nno\\
&& =\iota_{12}F_{\alpha}(w'_{(a_{4})}, \tilde{w}_{(a_{1})},
\tilde{w}_{(a_{2})}, \tilde{w}_{(a_{3})},  [\mathcal{Z}]_P; x_{1}, x_{2}),
\end{eqnarray}
\begin{eqnarray}
\lefteqn{\langle w_{(a_4)}', g^{a_{1}, a_{2}, a_{3},a_{4}}_{\alpha}(\tilde{w}_{(a_{1})},
\tilde{w}_{(a_{2})}, \tilde{w}_{(a_{3})}, \mathcal{ B}([\mathcal{Z}]_P); x_1, x_2) \rangle_{W^{a_4}}}\nno\\
&& =\iota_{21}F_{\alpha}(w'_{(a_{4})}, \tilde{w}_{(a_{1})},
\tilde{w}_{(a_{2})}, \tilde{w}_{(a_{3})}, [\mathcal{Z}]_P; x_{1}, x_{2})
\end{eqnarray}
and
\begin{eqnarray}
\lefteqn{ \langle w_{(a_4)}', h^{a_{1}, a_{2}, a_{3},a_{4}}_{\alpha}(\tilde{w}_{(a_{1})},
\tilde{w}_{(a_{2})}, \tilde{w}_{(a_{3})}, \mathcal{F}([\mathcal{Z}]_P); x_0, x_2) \rangle_{W^{a_4}}}\nno\\
&& =\iota_{20}F_{\alpha}(w'_{(a_{4})}, \tilde{w}_{(a_{1})},
\tilde{w}_{(a_{2})}, \tilde{w}_{(a_{3})},  [\mathcal{Z}]_P; x_2+x_0, x_{2})
\end{eqnarray}
for $\alpha\in \mathbb{ A}(a_{1}, a_{2}, a_{3}, a_{4})$.
Moreover, since
$\mathbb{ G}^{a_{1}, a_{2}, a_{3}, a_{4}}$ is a free module over the ring
\begin{equation}
\mathbb{ C}[x_{1}, x_{1}^{-1}, x_{2}, x_{2}^{-1}, (x_{1}-x_{2})^{-1}]
 \end{equation}
 with a basis $\{e^{a_{1}, a_{2}, a_{3}, a_{4}}_{\alpha}\}_{\alpha\in
\mathbb{ A}(a_{1}, a_{2}, a_{3}, a_{4})}$, we have
\begin{eqnarray}
\mathbb{ G}^{a_{1}, a_{2}, a_{3}, a_{4}} &\ni&  \Phi(w'_{(a_{4})}, \tilde{w}_{(a_{1})},
\tilde{w}_{(a_{2})}, \tilde{w}_{(a_{3})}, [\mathcal{Z}]_P;z_1,z_2)\nno\\
&& =\sum_{\alpha\in \mathbb{ A}(a_{1}, a_{2}, a_{3}, a_{4})}F_{\alpha}(w'_{(a_{4})}, \tilde{w}_{(a_{1})},
\tilde{w}_{(a_{2})}, \tilde{w}_{(a_{3})}, [\mathcal{Z}]_P; x_{1}, x_{2})e^{a_{1}, a_{2}, a_{3}, a_{4}}_{\alpha},\quad \label{e1:34}
\end{eqnarray}
 and
\begin{eqnarray}
\lefteqn{\langle w_{(a_4)}', (\tilde{\bf P}([\mathcal{Z}]_P))(\tilde{w}_{(a_{1})},
\tilde{w}_{(a_{2})}, \tilde{w}_{(a_{3})}; x_{1}, x_{2}) \rangle_{W^{a_4}}}\nno\\
&& =\iota_{12}\Phi(w'_{(a_{4})}, \tilde{w}_{(a_{1})},
\tilde{w}_{(a_{2})}, \tilde{w}_{(a_{3})}, [\mathcal{Z}]_P;z_1,z_2),\label{e1:35}
\end{eqnarray}
\begin{eqnarray}
\lefteqn{\langle w_{(a_4)}', (\tilde{\bf P}(\mathcal{B}([\mathcal{Z}]_P)))(\tilde{w}_{(a_{2})},
\tilde{w}_{(a_{1})}, \tilde{w}_{(a_{3})}; x_{2}, x_{1})\rangle_{W^{a_4}}}\nno\\
&&=\iota_{21}\Phi(w'_{(a_{4})}, \tilde{w}_{(a_{1})},
\tilde{w}_{(a_{2})}, \tilde{w}_{(a_{3})}, [\mathcal{Z}]_P;z_1,z_2),\label{e1:36}
\end{eqnarray}
\begin{eqnarray}
\lefteqn{\langle w_{(a_4)}', (\tilde{\bf I}(\mathcal{ F}([\mathcal{Z}]_P)))(\tilde{w}_{(a_{1})}, \tilde{w}_{(a_{2})},
\tilde{w}_{(a_{3})}; x_{0}, x_{2})\rangle_{W^{a_4}}}\nno\\
&& =\iota_{20}\Phi(w'_{(a_{4})}, \tilde{w}_{(a_{1})},
\tilde{w}_{(a_{2})}, \tilde{w}_{(a_{3})},[\mathcal{Z}]_P;z_1,z_2).\label{e1:37}
\end{eqnarray}
So the preferred branches of $\Phi(w'_{(a_{4})}, \tilde{w}_{(a_{1})},
\tilde{w}_{(a_{2})}, \tilde{w}_{(a_{3})}, [\mathcal{Z}]_P;z_1,z_2)$ on $R_1$, $R_2$ and $R_3$
are
\begin{equation}\label{a13}
\langle w'_{(a_{4})}, (\tilde{\bf P}([\mathcal{Z}]_P))(\tilde{w}_{(a_{1})},
\tilde{w}_{(a_{2})}, \tilde{w}_{(a_{3})}; x_{1}, x_{2})\rangle_{W^{a_4}}
\mbar_{x_{1}^{n}=e^{n \log z_1},
x_{2}^{n}=e^{n\log z_2}},
\end{equation}
\begin{equation}\label{a36}
\langle w_{(a_4)}', (\tilde{\bf P}(\mathcal{B}([\mathcal{Z}]_P)))(\tilde{w}_{(a_{2})},
\tilde{w}_{(a_{1})}, \tilde{w}_{(a_{3})}; x_{2}, x_{1})\rangle_{W^{a_4}}\mbar_{x_{1}^{n}=e^{n \log z_1},
x_{2}^{n}=e^{n\log z_2}}
\end{equation}
and
\begin{equation}
\langle w_{(a_4)}', (\tilde{\bf I}(\mathcal{ F}([\mathcal{Z}]_P)))(\tilde{w}_{(a_{1})}, \tilde{w}_{(a_{2})},
\tilde{w}_{(a_{3})}; x_{0}, x_{2})\rangle_{W^{a_4}}\mbar_{x_{0}^{n}=e^{n \log (z_1-z_2)},
x_{2}^{n}=e^{n\log z_2}},
\end{equation}
respectively. And the multivalued analytic functions
\begin{equation}
\langle w'_{(a_{4})}, (\tilde{\bf P}([\mathcal{Z}]_P))(\tilde{w}_{(a_{1})},
\tilde{w}_{(a_{2})}, \tilde{w}_{(a_{3})}; x_{1}, x_{2})\rangle_{W^{a_4}}\mbar_{x_1=z_1,x_2=z_2},
\end{equation}
\begin{equation}\label{a17}
\langle w_{(a_4)}', (\tilde{\bf P}(\mathcal{B}([\mathcal{Z}]_P)))(\tilde{w}_{(a_{2})},
\tilde{w}_{(a_{1})}, \tilde{w}_{(a_{3})}; x_{2}, x_{1})\rangle_{W^{a_4}}\mbar_{x_1=z_1,x_2=z_2}
\end{equation}
and
\begin{equation}
\langle w_{(a_4)}', (\tilde{\bf I}(\mathcal{ F}([\mathcal{Z}]_P)))(\tilde{w}_{(a_{1})}, \tilde{w}_{(a_{2})},
\tilde{w}_{(a_{3})}; x_{0}, x_{2})\rangle_{W^{a_4}}\mbar_{x_0=z_1-z_2,x_2=z_2}
\end{equation}
are restrictions of the multivalued analytic function $\Phi(w'_{(a_{4})}, \tilde{w}_{(a_{1})},
\tilde{w}_{(a_{2})}, \tilde{w}_{(a_{3})},[\mathcal{Z}]_P;z_1,z_2)$ to their domains $|z_{1}|>|z_{2}|>0$, $|z_{2}|> |z_{1}|>0$ and $|z_{2}|>|z_{1}-z_{2}|>0$ respectively.

Then by the definition of the preferred branch of an element of $\mathbb{ G}^{a_{1}, a_{2}, a_{3}, a_{4}}$ on $R_3$, we can deduce that
\begin{eqnarray}
\lefteqn{ \langle w_{(a_4)}', (\tilde{\bf I}(\mathcal{ F}([\mathcal{Z}]_P)))(\tilde{w}_{(a_{1})}, \tilde{w}_{(a_{2})},
\tilde{w}_{(a_{3})}; x_{0}, x_{2})\rangle_{W^{a_4}}\mbar_{x_{0}^{n}=e^{n \log (z_1-z_2)},
x_{2}^{n}=e^{n\log z_2}}}\nno\\
&&= \langle w'_{(a_{4})}, (\tilde{\bf P}([\mathcal{Z}]_P))(\tilde{w}_{(a_{1})},
\tilde{w}_{(a_{2})}, \tilde{w}_{(a_{3})}; x_{1}, x_{2})\rangle_{W^{a_4}}
\mbar_{x_{1}^{n}=e^{n \log z_1},
x_{2}^{n}=e^{n\log z_2}} \label{a7}
\end{eqnarray}
on the region
\begin{equation*}
S_1=\{(z_1,z_2)\in \mathbb{C}^2 \mid \re z_1>\re z_2 > \re (z_1-z_2) >0,\ \im z_1 >  \im z_2 > \im (z_1-z_2) > 0\}.
\end{equation*}
Moreover, by skew-symmetry and (\ref{a7}), we see that
\begin{eqnarray}
\lefteqn{\langle w_{(a_4)}', (\tilde{\bf I}(\tilde{\Omega}^{(1)}(\mathcal{F}([\mathcal{Z}]_P)))) (\tilde{w}_{(a_{2})}, \tilde{w}_{(a_{1})},
\tilde{w}_{(a_{3})}; x_{0}, x_1)\rangle_{W^{a_4}}\mbar_{x_{0}^{n}=e^{n \log (z_2-z_1)},
x_{1}^{n}=e^{n\log z_1}}}\nno\\
&& =\langle w_{(a_4)}', (\tilde{\bf I}(\mathcal{F}([\mathcal{Z}]_P))) (\tilde{w}_{(a_{1})}, \tilde{w}_{(a_{2})},
\tilde{w}_{(a_{3})}; e^{-\pi i}x_{0}, x_2)\rangle_{W^{a_4}}\mbar_{x_{0}^{n}=e^{n \log (z_2-z_1)},
x_{2}^{n}=e^{n\log z_2}}\nno\\
&& =\langle w_{(a_4)}', (\tilde{\bf I}(\mathcal{F}([\mathcal{Z}]_P))) (\tilde{w}_{(a_{1})}, \tilde{w}_{(a_{2})},
\tilde{w}_{(a_{3})}; x_{0}, x_2)\rangle_{W^{a_4}}\mbar_{x_{0}^{n}=e^{n \log (z_1-z_2)},
x_{2}^{n}=e^{n\log z_2}}\nno\\
&& =\langle w'_{(a_{4})}, (\tilde{\bf P}([\mathcal{Z}]_P))(\tilde{w}_{(a_{1})},
\tilde{w}_{(a_{2})}, \tilde{w}_{(a_{3})}; x_{1}, x_{2})\rangle_{W^{a_4}}
\mbar_{x_{1}^{n}=e^{n \log z_1},
x_{2}^{n}=e^{n\log z_2}}
\end{eqnarray}
on the region $S_1$. So by the definition of the preferred branch of an element of $\mathbb{ G}^{a_{1}, a_{2}, a_{3}, a_{4}}$ on $R_4$, we deduce that the preferred branch on $R_4$ of $\Phi(w'_{(a_{4})}, \tilde{w}_{(a_{1})},
\tilde{w}_{(a_{2})}, \tilde{w}_{(a_{3})}, [\mathcal{Z}]_P;z_1,z_2)$ is equal to the single-valued analytic function
\begin{equation}\label{b9}
\langle w_{(a_4)}', (\tilde{\bf I}(\tilde{\Omega}^{(1)}(\mathcal{F}([\mathcal{Z}]_P)))) (\tilde{w}_{(a_{2})}, \tilde{w}_{(a_{1})},
\tilde{w}_{(a_{3})}; x_{0}, x_1)\rangle_{W^{a_4}}\mbar_{x_{0}^{n}=e^{n \log (z_2-z_1)},
x_{1}^{n}=e^{n\log z_1}}
\end{equation}
on the region $S_1$. Since the preferred branch on $R_4$ of $\Phi(w'_{(a_{4})}, \tilde{w}_{(a_{1})},
\tilde{w}_{(a_{2})}, \tilde{w}_{(a_{3})}, [\mathcal{Z}]_P;z_1,z_2)$ and the function (\ref{b9}) are both single-valued analytic functions on the domain $R_4$ which contains $S_1$, by the basic properties of analytic functions we conclude that they are equal on $R_4$; namely, the single-valued analytic function (\ref{b9})
defined on the region $R_4$ is the preferred branch of $\Phi(w'_{(a_{4})}, \tilde{w}_{(a_{1})},
\tilde{w}_{(a_{2})}, \tilde{w}_{(a_{3})}, [\mathcal{Z}]_P;z_1,z_2)$ on $R_4$.
Furthermore, by the definition of the preferred branch of an element of $\mathbb{ G}^{a_{1}, a_{2}, a_{3}, a_{4}}$ on $R_2$, we can conclude that
\begin{eqnarray}
\lefteqn{\langle w_{(a_4)}', (\tilde{\bf I}(\tilde{\Omega}^{(1)}(\mathcal{F}([\mathcal{Z}]_P)))) (\tilde{w}_{(a_{2})}, \tilde{w}_{(a_{1})},
\tilde{w}_{(a_{3})}; x_{0}, x_1)\rangle_{W^{a_4}}\mbar_{x_{0}^{n}=e^{n \log (z_2-z_1)},
x_{1}^{n}=e^{n\log z_1}}}\nno\\
&& =\langle w_{(a_4)}', (\tilde{\bf P}(\mathcal{B}([\mathcal{Z}]_P)))(\tilde{w}_{(a_{2})},
\tilde{w}_{(a_{1})}, \tilde{w}_{(a_{3})}; x_{2}, x_{1})\rangle_{W^{a_4}}\mbar_{x_{1}^{n}=e^{n \log z_1},
x_{2}^{n}=e^{n\log z_2}}\label{a9}
\end{eqnarray}
on the region
\begin{equation*}
S_2=\{(z_1,z_2)\in \mathbb{C}^2 \mid \re z_2>\re z_1 > \re (z_2-z_1) >0,\ \im z_2 >  \im z_1 > \im (z_2-z_1) > 0\}.
\end{equation*}
So this proposition holds.\epfv
\vspace{0.3cm}

\begin{thm}\label{l2}
Assume that the assumptions of Theorem \ref{t1} hold, then the Jacobi identity holds for the ordered triple
$(\tilde{w}_{(a_{2})}, \tilde{w}_{(a_{1})}, \tilde{w}_{(a_{3})})$.
\end{thm}

\pf
Consider any $a_{4}\in \mathcal{A}$, $w_{(a_4)}'\in (W^{a_4})'$ and $\mathcal{ Z}\in \coprod_{ a_{5}\in \mathcal{ A}}
\mathcal{ V}_{a_{1}a_{5}}^{a_{4}}\otimes
\mathcal{ V}_{a_{2}a_{3}}^{a_{5}}$. Since the assumptions of Theorem \ref{t1} contain the assumptions of Proposition \ref{l1}, we obtain a multivalued analytic function $\Phi(w'_{(a_{4})}, \tilde{w}_{(a_{1})},
\tilde{w}_{(a_{2})}, \tilde{w}_{(a_{3})}, [\mathcal{Z}]_P;z_1,z_2)$ (see (\ref{s7})) such that (\ref{s8})-(\ref{s11}) are its preferred branches on $R_1$, $R_2$, $R_3$ and $R_4$, respectively. Moreover, the preferred branches on $R_1$, $R_2$, $R_3$ and $R_4$ have relations (\ref{s12})-(\ref{s13}).
Interchanging $z_1$ and $z_2$ in the multivalued analytic function $\Phi(w'_{(a_{4})}, \tilde{w}_{(a_{1})},
\tilde{w}_{(a_{2})}, \tilde{w}_{(a_{3})}, [\mathcal{Z}]_P;z_1,z_2)$, we obtain another multivalued analytic function
\begin{equation}
\Phi(w'_{(a_{4})}, \tilde{w}_{(a_{1})},
\tilde{w}_{(a_{2})}, \tilde{w}_{(a_{3})}, [\mathcal{Z}]_P;z_2,z_1)
\end{equation}
 on $M^{2}=\{(z_{1},
z_{2})\in \mathbb{ C}^{2}\;|\;z_{1}, z_{2}\ne 0, z_{1}\ne
z_{2}\}$.
By interchanging $z_1$ and $z_2$ in (\ref{s9}), we see that
\begin{equation}\label{a12}
\langle w_{(a_4)}', (\tilde{\bf P}(\mathcal{B}([\mathcal{Z}]_P)))(\tilde{w}_{(a_{2})},
\tilde{w}_{(a_{1})}, \tilde{w}_{(a_{3})}; x_{1}, x_{2})\rangle_{W^{a_4}}\mbar_{x_{1}^{n}=e^{n \log z_1},
x_{2}^{n}=e^{n\log z_2}}
\end{equation}
is a branch of $\Phi(w'_{(a_{4})}, \tilde{w}_{(a_{1})},
\tilde{w}_{(a_{2})}, \tilde{w}_{(a_{3})}, [\mathcal{Z}]_P;z_2,z_1)$ on the region $R_1$.

Moreover, by interchanging $z_1$ and $z_2$ in (\ref{s13}), we get
\begin{eqnarray}
\lefteqn{\langle w_{(a_4)}', (\tilde{\bf I}(\tilde{\Omega}^{(1)}(\mathcal{F}([\mathcal{Z}]_P)))) (\tilde{w}_{(a_{2})}, \tilde{w}_{(a_{1})},
\tilde{w}_{(a_{3})}; x_{0}, x_2)\rangle_{W^{a_4}}\mbar_{x_{0}^{n}=e^{n \log (z_1-z_2)},
x_{2}^{n}=e^{n\log z_2}}}\nno\\
&& =\langle w_{(a_4)}', (\tilde{\bf P}(\mathcal{B}([\mathcal{Z}]_P)))(\tilde{w}_{(a_{2})},
\tilde{w}_{(a_{1})}, \tilde{w}_{(a_{3})}; x_{1}, x_{2})\rangle_{W^{a_4}}\mbar_{x_{1}^{n}=e^{n \log z_1},
x_{2}^{n}=e^{n\log z_2}}\label{a21}
\end{eqnarray}
on the region $S_1$. Since the single-valued analytic function
\begin{eqnarray}
\lefteqn{\langle w_{(a_4)}', (\tilde{\bf I}(\mathcal{F}(\mathcal{B}([\mathcal{Z}]_P))))(\tilde{w}_{(a_{2})},
\tilde{w}_{(a_{1})}, \tilde{w}_{(a_{3})}; x_{0}, x_{2})\rangle_{W^{a_4}}\mbar_{x_{0}^{n}=e^{n \log (z_1-z_2)},
x_{2}^{n}=e^{n\log z_2}}}\nno\\
&& =\langle w_{(a_4)}', (\tilde{\bf I}(\tilde{\Omega}^{(1)}(\mathcal{F}([\mathcal{Z}]_P)))) (\tilde{w}_{(a_{2})}, \tilde{w}_{(a_{1})},
\tilde{w}_{(a_{3})}; x_{0}, x_2)\rangle_{W^{a_4}}\mbar_{x_{0}^{n}=e^{n \log (z_1-z_2)},
x_{2}^{n}=e^{n\log z_2}}\qquad\label{a31}
\end{eqnarray}
on $S_1$ can be naturally analytically extended to the region $R_3$, we therefore get a branch of $\Phi(w'_{(a_{4})}, \tilde{w}_{(a_{1})},
\tilde{w}_{(a_{2})}, \tilde{w}_{(a_{3})}, [\mathcal{Z}]_P;z_2,z_1)$ on the region $R_3$.

By (\ref{a21}) and skew-symmetry, we have
\begin{eqnarray}
\lefteqn{ \langle w_{(a_4)}', (\tilde{\bf I}(\tilde{\Omega}^{(1)}\tilde{\Omega}^{(1)}(\mathcal{F}([\mathcal{Z}]_P))))
(\tilde{w}_{(a_{1})},
\tilde{w}_{(a_{2})}, \tilde{w}_{(a_{3})}; x_{0}, x_{1})\rangle_{W^{a_4}}\mbar_{x_{0}^{n}=e^{n \log (z_2-z_1)},
x_{1}^{n}=e^{n\log z_1}}}\nno\\
&& =\langle w_{(a_4)}', (\tilde{\bf I}(\tilde{\Omega}^{(1)}(\mathcal{F}([\mathcal{Z}]_P))))
(\tilde{w}_{(a_{2})},
\tilde{w}_{(a_{1})}, \tilde{w}_{(a_{3})}; e^{-\pi i}x_{0}, x_{2})\rangle_{W^{a_4}}\mbar_{x_{0}^{n}=e^{n \log (z_2-z_1)},
x_{2}^{n}=e^{n\log z_2}}\nno\\
&& =\langle w_{(a_4)}', (\tilde{\bf I}(\tilde{\Omega}^{(1)}(\mathcal{F}([\mathcal{Z}]_P))))
(\tilde{w}_{(a_{2})},
\tilde{w}_{(a_{1})}, \tilde{w}_{(a_{3})}; x_{0}, x_{2})\rangle_{W^{a_4}}\mbar_{x_{0}^{n}=e^{n \log (z_1-z_2)},
x_{2}^{n}=e^{n\log z_2}}\nno\\
&& =\langle w_{(a_4)}', (\tilde{\bf P}(\mathcal{B}([\mathcal{Z}]_P)))(\tilde{w}_{(a_{2})},
\tilde{w}_{(a_{1})}, \tilde{w}_{(a_{3})}; x_{1}, x_{2})\rangle_{W^{a_4}}\mbar_{x_{1}^{n}=e^{n \log z_1},
x_{2}^{n}=e^{n\log z_2}}\label{a22}
\end{eqnarray}
on the region $S_1$. Moreover, the first line of (\ref{a22}) on $S_1$ can be naturally analytically extended to the region $R_4$. So the single-valued analytic function
\begin{equation}\label{a30}
\langle w_{(a_4)}', (\tilde{\bf I}(\tilde{\Omega}^{(1)}\tilde{\Omega}^{(1)}(\mathcal{F}([\mathcal{Z}]_P))))
(\tilde{w}_{(a_{1})},
\tilde{w}_{(a_{2})}, \tilde{w}_{(a_{3})}; x_{0}, x_{1})\rangle_{W^{a_4}}\mbar_{x_{0}^{n}=e^{n \log (z_2-z_1)},
x_{1}^{n}=e^{n\log z_1}}
\end{equation}
defined on the region $R_4$ is a branch of $\Phi(w'_{(a_{4})}, \tilde{w}_{(a_{1})},
\tilde{w}_{(a_{2})}, \tilde{w}_{(a_{3})}, [\mathcal{Z}]_P;z_2,z_1)$ on the region $R_4$.

Observe that $[\mathcal{Z}]_P \in \pi_P(\coprod_{a_{5}\in \mathcal{ A}}
\mathcal{ V}_{a_{1}a_{5}}^{a_{4}}\otimes
\mathcal{ V}_{a_{2}a_{3}}^{a_{5}})$ implies
\begin{equation}
\mathcal{B}(\mathcal{B}([\mathcal{Z}]_P)) \in \pi_P(\coprod_{a_{5}\in \mathcal{ A}}
\mathcal{ V}_{a_{1}a_{5}}^{a_{4}}\otimes
\mathcal{ V}_{a_{2}a_{3}}^{a_{5}}).
\end{equation}
And since (\ref{s7})-(\ref{s13}) hold for any $[\mathcal{Z}]_P \in \pi_P(\coprod_{a_{5}\in \mathcal{ A}}
\mathcal{ V}_{a_{1}a_{5}}^{a_{4}}\otimes
\mathcal{ V}_{a_{2}a_{3}}^{a_{5}})$, they should hold with $[\mathcal{Z}]_P$ replaced by $\mathcal{B}(\mathcal{B}([\mathcal{Z}]_P))$ for any $[\mathcal{Z}]_P \in \pi_P(\coprod_{a_{5}\in \mathcal{ A}}
\mathcal{ V}_{a_{1}a_{5}}^{a_{4}}\otimes
\mathcal{ V}_{a_{2}a_{3}}^{a_{5}})$. So
replacing $[\mathcal{Z}]_P$ by $\mathcal{B}(\mathcal{B}([\mathcal{Z}]_P))$ in (\ref{s12}), we get
\begin{eqnarray}
\lefteqn{ \langle w'_{(a_{4})}, (\tilde{\bf P}(\mathcal{B}(\mathcal{B}([\mathcal{Z}]_P))))(\tilde{w}_{(a_{1})},
\tilde{w}_{(a_{2})}, \tilde{w}_{(a_{3})}; x_{1}, x_{2})\rangle_{W^{a_4}}
\mbar_{x_{1}^{n}=e^{n \log z_1},
x_{2}^{n}=e^{n\log z_2}}}\nno\\
&& =\langle w_{(a_4)}', (\tilde{\bf I}(\mathcal{ F}(\mathcal{B}(\mathcal{B}([\mathcal{Z}]_P)))))(\tilde{w}_{(a_{1})}, \tilde{w}_{(a_{2})},
\tilde{w}_{(a_{3})}; x_{0}, x_{2})\rangle_{W^{a_4}}\mbar_{x_{0}^{n}=e^{n \log (z_1-z_2)},
x_{2}^{n}=e^{n\log z_2}} \qquad\label{a11}
\end{eqnarray}
on the region $S_1$. Interchanging $z_1$ and $z_2$ in (\ref{a11}), we obtain
\begin{eqnarray}
\lefteqn{ \langle w'_{(a_{4})}, (\tilde{\bf P}(\mathcal{B}(\mathcal{B}([\mathcal{Z}]_P))))(\tilde{w}_{(a_{1})},
\tilde{w}_{(a_{2})}, \tilde{w}_{(a_{3})}; x_{2}, x_{1})\rangle_{W^{a_4}}
\mbar_{x_{1}^{n}=e^{n \log z_1},
x_{2}^{n}=e^{n\log z_2}}}\nno\\
&& =\langle w_{(a_4)}', (\tilde{\bf I}(\mathcal{ F}(\mathcal{B}(\mathcal{B}([\mathcal{Z}]_P)))))(\tilde{w}_{(a_{1})}, \tilde{w}_{(a_{2})},
\tilde{w}_{(a_{3})}; x_{0}, x_{1})\rangle_{W^{a_4}}\mbar_{x_{0}^{n}=e^{n \log (z_2-z_1)},
x_{1}^{n}=e^{n\log z_1}} \nno\\
&& =\langle w_{(a_4)}', (\tilde{\bf I}(\tilde{\Omega}^{(1)}\tilde{\Omega}^{(1)}\mathcal{F}([\mathcal{Z}]_P)))
(\tilde{w}_{(a_{1})},
\tilde{w}_{(a_{2})}, \tilde{w}_{(a_{3})}; x_{0}, x_{1})\rangle_{W^{a_4}}\mbar_{x_{0}^{n}=e^{n \log (z_2-z_1)},
x_{1}^{n}=e^{n\log z_1}}\qquad\ \label{a23}
\end{eqnarray}
on the region $S_2$. Moreover, the first line of (\ref{a23}) on $S_2$ can be naturally analytically extended to the region $R_2$. So the single-valued analytic function
\begin{equation}\label{a32}
\langle w'_{(a_{4})}, (\tilde{\bf P}(\mathcal{B}(\mathcal{B}([\mathcal{Z}]_P))))(\tilde{w}_{(a_{1})},
\tilde{w}_{(a_{2})}, \tilde{w}_{(a_{3})}; x_{2}, x_{1})\rangle_{W^{a_4}}
\mbar_{x_{1}^{n}=e^{n \log z_1},
x_{2}^{n}=e^{n\log z_2}}
\end{equation}
defined on the region $R_2$ is a branch of $\Phi(w'_{(a_{4})}, \tilde{w}_{(a_{1})},
\tilde{w}_{(a_{2})}, \tilde{w}_{(a_{3})}, [\mathcal{Z}]_P;z_2,z_1)$ on $R_2$.

From the above discussion (\ref{a12})-(\ref{a32}), we see that the multivalued analytic functions
\begin{equation}
\langle w_{(a_4)}', (\tilde{\bf P}(\mathcal{B}([\mathcal{Z}]_P)))(\tilde{w}_{(a_{2})},
\tilde{w}_{(a_{1})}, \tilde{w}_{(a_{3})}; x_{1}, x_{2})\rangle_{W^{a_4}}\mbar_{x_1=z_1,x_2=z_2},
\end{equation}
\begin{equation}
\langle w'_{(a_{4})}, (\tilde{\bf P}(\mathcal{B}(\mathcal{B}([\mathcal{Z}]_P))))(\tilde{w}_{(a_{1})},
\tilde{w}_{(a_{2})}, \tilde{w}_{(a_{3})}; x_{2}, x_{1})\rangle_{W^{a_4}}\mbar_{x_1=z_1,x_2=z_2}
\end{equation}
and
\begin{equation}
\langle w_{(a_4)}', (\tilde{\bf I}(\mathcal{F}(\mathcal{B}([\mathcal{Z}]_P))))(\tilde{w}_{(a_{2})},
\tilde{w}_{(a_{1})}, \tilde{w}_{(a_{3})}; x_{0}, x_{2})\rangle_{W^{a_4}}\mbar_{x_0=z_1-z_2,x_2=z_2}
\end{equation}
are restrictions of the multivalued analytic function $\Phi(w'_{(a_{4})}, \tilde{w}_{(a_{1})},
\tilde{w}_{(a_{2})}, \tilde{w}_{(a_{3})},[\mathcal{Z}]_P;z_2,z_1)$ to their domains $|z_{1}|>|z_{2}|>0$, $|z_{2}|> |z_{1}|>0$ and $|z_{2}|>|z_{1}-z_{2}|>0$ respectively.
So with the branch (\ref{a12}) chosen as the preferred branch on $R_1$, $\Phi(w'_{(a_{4})}, \tilde{w}_{(a_{1})},
\tilde{w}_{(a_{2})}, \tilde{w}_{(a_{3})},[\mathcal{Z}]_P;z_2,z_1)$ becomes an element of $\mathbb{ G}^{a_{2}, a_{1}, a_{3}, a_{4}}$. Moreover, by (\ref{a21}), (\ref{a31}), (\ref{a22}), (\ref{a30}), (\ref{a23}), (\ref{a32}), and by the definition of the preferred branches of an element of $\mathbb{ G}^{a_{2}, a_{1}, a_{3}, a_{4}}$ on $R_2$ and $R_3$, we see that
\begin{equation}
\langle w'_{(a_{4})}, (\tilde{\bf P}(\mathcal{B}(\mathcal{B}([\mathcal{Z}]_P))))(\tilde{w}_{(a_{1})},
\tilde{w}_{(a_{2})}, \tilde{w}_{(a_{3})}; x_{2}, x_{1})\rangle_{W^{a_4}}
\mbar_{x_{1}^{n}=e^{n \log z_1},
x_{2}^{n}=e^{n\log z_2}}
\end{equation}
and
\begin{equation}
\langle w_{(a_4)}', (\tilde{\bf I}(\mathcal{F}(\mathcal{B}([\mathcal{Z}]_P))))(\tilde{w}_{(a_{2})},
\tilde{w}_{(a_{1})}, \tilde{w}_{(a_{3})}; x_{0}, x_{2})\rangle_{W^{a_4}}\mbar_{x_{0}^{n}=e^{n \log (z_1-z_2)},
x_{2}^{n}=e^{n\log z_2}}
\end{equation}
are the preferred branches of $\Phi(w'_{(a_{4})}, \tilde{w}_{(a_{1})},
\tilde{w}_{(a_{2})}, \tilde{w}_{(a_{3})},[\mathcal{Z}]_P;z_2,z_1)$ on $R_2$ and $R_3$ respectively.
Therefore, we have
\begin{eqnarray}
\lefteqn{ \langle w_{(a_4)}', (\tilde{\bf P}(\mathcal{B}([\mathcal{Z}]_P)))(\tilde{w}_{(a_{2})},
\tilde{w}_{(a_{1})}, \tilde{w}_{(a_{3})}; x_1, x_2)\rangle_{W^{a_4}}}\nno\\
&& =\iota_{12}\Phi(w'_{(a_{4})}, \tilde{w}_{(a_{1})},
\tilde{w}_{(a_{2})}, \tilde{w}_{(a_{3})}, [\mathcal{Z}]_P;z_2,z_1),\label{e1:25}
\end{eqnarray}
\begin{eqnarray}
\lefteqn{ \langle w_{(a_4)}', (\tilde{\bf P}(\mathcal{B}(\mathcal{B}([\mathcal{Z}]_P))))(\tilde{w}_{(a_{1})},
\tilde{w}_{(a_{2})}, \tilde{w}_{(a_{3})}; x_{2}, x_{1}) \rangle_{W^{a_4}}}\nno\\
&& =\iota_{21}\Phi(w'_{(a_{4})}, \tilde{w}_{(a_{1})},
\tilde{w}_{(a_{2})}, \tilde{w}_{(a_{3})}, [\mathcal{Z}]_P;z_2,z_1),\label{e1:26}
\end{eqnarray}
\begin{eqnarray}
\lefteqn{\langle w_{(a_4)}', (\tilde{\bf I}(\mathcal{ F}(\mathcal{B}([\mathcal{Z}]_P))))(\tilde{w}_{(a_{2})}, \tilde{w}_{(a_{1})},
\tilde{w}_{(a_{3})}; x_{0}, x_{2})\rangle_{W^{a_4}}}\nno\\
&& =\iota_{20} \Phi(w'_{(a_{4})}, \tilde{w}_{(a_{1})},
\tilde{w}_{(a_{2})}, \tilde{w}_{(a_{3})}, [\mathcal{Z}]_P;z_2,z_1).\label{e1:27}
\end{eqnarray}
Let $\{e^{a_{2}, a_{1}, a_{3}, a_{4}}_{\alpha}\}_{\alpha\in
\mathbb{ A}(a_{2}, a_{1}, a_{3}, a_{4})}$ be a basis of
$\mathbb{ G}^{a_{2}, a_{1}, a_{3}, a_{4}}$ over the ring
\begin{equation}
\mathbb{ C}[x_{1}, x_{1}^{-1}, x_{2}, x_{2}^{-1}, (x_{1}-x_{2})^{-1}].
 \end{equation}
 Then there exists unique
\begin{equation}\label{e1:32}
G_{\alpha}(w'_{(a_{4})},\tilde{w}_{(a_{2})}, \tilde{w}_{(a_{1})},
\tilde{w}_{(a_{3})}, \mathcal{B}([\mathcal{Z}]_P); x_{1}, x_{2})\in
\mathbb{ C}[x_{1}, x_{1}^{-1}, x_{2}, x_{2}^{-1}, (x_{1}-x_{2})^{-1}]
\end{equation}
for $\alpha\in \mathbb{ A}(a_{2}, a_{1}, a_{3}, a_{4})$, such that only finitely many of them are nonzero and
\begin{eqnarray}
\lefteqn{ \Phi(w'_{(a_{4})}, \tilde{w}_{(a_{1})},
\tilde{w}_{(a_{2})}, \tilde{w}_{(a_{3})}, [\mathcal{Z}]_P;z_2,z_1)}\nno\\
&& =\sum_{\alpha\in \mathbb{ A}(a_{2}, a_{1}, a_{3}, a_{4})}G_{\alpha}(w'_{(a_{4})},\tilde{w}_{(a_{2})}, \tilde{w}_{(a_{1})},
\tilde{w}_{(a_{3})}, \mathcal{B}([\mathcal{Z}]_P); x_{1}, x_{2})e^{a_{2}, a_{1}, a_{3}, a_{4}}_{\alpha}.\quad \label{e1:28}
\end{eqnarray}
By (\ref{e3:1}) and (\ref{a18}), we see that $\mathcal{B}=\F^{-1}\tilde{\Omega}^{(1)}\mathcal{F}$ is an isomorphism and that
\begin{equation}\label{a19}
\mathcal{B}(\pi_P(\coprod_{a_{5}\in \mathcal{ A}}\mathcal{ V}_{a_{1}a_{5}}^{a_{4}}\otimes
\mathcal{ V}_{a_{2}a_{3}}^{a_{5}}))=\pi_P(\coprod_{a_{5}\in \mathcal{ A}}
\mathcal{ V}_{a_{2}a_{5}}^{a_{4}}\otimes
\mathcal{ V}_{a_{1}a_{3}}^{a_{5}})
\end{equation}
for any $a_1,\cdots,a_4\in\A$.
So we can define linear maps
\begin{eqnarray}
\lefteqn{f^{a_{2}, a_{1}, a_{3},
a_{4}}_{\alpha}(\tilde{w}_{(a_{2})}, \tilde{w}_{(a_{1})},
\tilde{w}_{(a_{3})}):}\nno\\
&&\hspace{4em}\pi_P(\coprod_{a_{5}\in \mathcal{ A}}
\mathcal{ V}_{a_{2}a_{5}}^{a_{4}}\otimes
\mathcal{ V}_{a_{1}a_{3}}^{a_{5}})\to  W^{a_{4}}[x_{1}, x_{1}^{-1}, x_{2},x_{2}^{-1}][[x_{2}/x_{1}]],
\end{eqnarray}
\begin{eqnarray}
\lefteqn{g^{a_{2}, a_{1}, a_{3},a_{4}}_{\alpha}(\tilde{w}_{(a_{2})}, \tilde{w}_{(a_{1})},\tilde{w}_{(a_{3})}):}\nno\\
&&\hspace{4em}
\pi_P(\coprod_{a_{5}\in \mathcal{ A}}
\mathcal{ V}_{a_{1}a_{5}}^{a_{4}}\otimes
\mathcal{ V}_{a_{2}a_{3}}^{a_{5}})  \to  W^{a_{4}}[x_{1}, x_{1}^{-1}, x_{2},
x_{2}^{-1}][[x_{1}/x_{2}]],
\end{eqnarray}
\begin{eqnarray}
\lefteqn{ h^{a_{2}, a_{1}, a_{3},
a_{4}}_{\alpha}(\tilde{w}_{(a_{2})}, \tilde{w}_{(a_{1})},
\tilde{w}_{(a_{3})}):}\nno\\
&&\hspace{4em}
\pi_I(\coprod_{a_{5}\in \mathcal{ A}}
\mathcal{ V}_{a_{2}a_{1}}^{a_{5}}\otimes
\mathcal{ V}_{a_{5}a_{3}}^{a_{4}}) \to W^{a_{4}}[x_{0},
x_{0}^{-1}, x_{2}, x_{2}^{-1}][[x_{0}/x_{2}]]
\end{eqnarray}
by
\begin{eqnarray}
\lefteqn{\langle w_{(a_4)}', f^{a_{2}, a_{1}, a_{3},
a_{4}}_{\alpha}(\tilde{w}_{(a_{2})}, \tilde{w}_{(a_{1})},
\tilde{w}_{(a_{3})}, \mathcal{B}([\mathcal{Z}]_P); x_{1}, x_{2})\rangle_{W^{a_4}}}\nno\\
&& =\iota_{12} G_{\alpha}(w'_{(a_{4})},\tilde{w}_{(a_{2})}, \tilde{w}_{(a_{1})},
\tilde{w}_{(a_{3})}, \mathcal{B}([\mathcal{Z}]_P); x_{1}, x_{2}),\label{e1:29}
\end{eqnarray}
\begin{eqnarray}
\lefteqn{\langle w_{(a_4)}', g^{a_{2}, a_{1}, a_{3},
a_{4}}_{\alpha}(\tilde{w}_{(a_{2})}, \tilde{w}_{(a_{1})},
\tilde{w}_{(a_{3})}, \mathcal{B}(\mathcal{B}([\mathcal{Z}]_P)); x_{1}, x_{2})\rangle_{W^{a_4}}}\nno\\
&& =\iota_{21} G_{\alpha}(w'_{(a_{4})},\tilde{w}_{(a_{2})}, \tilde{w}_{(a_{1})},
\tilde{w}_{(a_{3})}, \mathcal{B}([\mathcal{Z}]_P); x_{1}, x_{2}),\label{e1:30}
\end{eqnarray}
\begin{eqnarray}
\lefteqn{\langle w_{(a_4)}', h^{a_{2}, a_{1}, a_{3},
a_{4}}_{\alpha}(\tilde{w}_{(a_{2})}, \tilde{w}_{(a_{1})},
\tilde{w}_{(a_{3})}, \mathcal{F}(\mathcal{B}([\mathcal{Z}]_P)); x_0, x_2)\rangle_{W^{a_4}}} \nno\\
&& =\iota_{20} G_{\alpha}(w'_{(a_{4})},\tilde{w}_{(a_{2})}, \tilde{w}_{(a_{1})},
\tilde{w}_{(a_{3})}, \mathcal{B}([\mathcal{Z}]_P); x_2+x_0, x_{2})\label{e1:31}
\end{eqnarray}
for $w_{(a_4)}'\in (W^{a_4})'$, $\mathcal{ Z}\in \coprod_{ a_{5}\in \mathcal{ A}}
\mathcal{ V}_{a_{1}a_{5}}^{a_{4}}\otimes
\mathcal{ V}_{a_{2}a_{3}}^{a_{5}}$ and $\alpha\in \mathbb{ A}(a_{2}, a_{1}, a_{3}, a_{4})$. Then by (\ref{e1:25})-(\ref{e1:27}) and (\ref{e1:28}), we have
\begin{eqnarray}
\lefteqn{(\tilde{\bf P}(\mathcal{B}([\mathcal{Z}]_P)))(\tilde{w}_{(a_{2})},
\tilde{w}_{(a_{1})}, \tilde{w}_{(a_{3})}; x_1, x_2)}\nno\\
&& =\sum_{\alpha\in \mathbb{ A}(a_{2}, a_{1}, a_{3}, a_{4})} f^{a_{2}, a_{1}, a_{3},
a_{4}}_{\alpha}(\tilde{w}_{(a_{2})}, \tilde{w}_{(a_{1})},
\tilde{w}_{(a_{3})}, \mathcal{B}([\mathcal{Z}]_P); x_{1}, x_{2}) \iota_{12}(e^{a_{2}, a_{1}, a_{3}, a_{4}}_{\alpha}),\label{a37}
\end{eqnarray}
\begin{eqnarray}
\lefteqn{ (\tilde{\bf P}(\mathcal{B}(\mathcal{B}([\mathcal{Z}]_P))))(\tilde{w}_{(a_{1})},
\tilde{w}_{(a_{2})}, \tilde{w}_{(a_{3})}; x_2, x_1)} \nno\\
&& =\sum_{\alpha\in \mathbb{ A}(a_{2}, a_{1}, a_{3}, a_{4})} g^{a_{2}, a_{1}, a_{3},
a_{4}}_{\alpha}(\tilde{w}_{(a_{2})}, \tilde{w}_{(a_{1})},
\tilde{w}_{(a_{3})}, \mathcal{B}(\mathcal{B}([\mathcal{Z}]_P)); x_{1}, x_{2}) \iota_{21}(e^{a_{2}, a_{1}, a_{3}, a_{4}}_{\alpha}),\qquad
\end{eqnarray}
\begin{eqnarray}
\lefteqn{ (\tilde{\bf I}(\mathcal{ F}(\mathcal{B}([\mathcal{Z}]_P))))(\tilde{w}_{(a_{2})}, \tilde{w}_{(a_{1})},
\tilde{w}_{(a_{3})};  x_{0}, x_2)} \nno\\
&& =\sum_{\alpha\in \mathbb{ A}(a_{2}, a_{1}, a_{3}, a_{4})} h^{a_{2}, a_{1}, a_{3},
a_{4}}_{\alpha}(\tilde{w}_{(a_{2})}, \tilde{w}_{(a_{1})},
\tilde{w}_{(a_{3})}, \mathcal{F}(\mathcal{B}([\mathcal{Z}]_P)); x_0, x_2)
\iota_{20} (e^{a_{2}, a_{1}, a_{3}, a_{4}}_{\alpha}).\qquad
\end{eqnarray}
Moreover, by (\ref{e1:32}) and Proposition \ref{p3.1.1}, we have
 \begin{eqnarray}
\lefteqn{x_{0}^{-1}
\delta\left(\frac{x_{1}-x_{2}}{x_{0}}\right)
\iota_{12} G_{\alpha}(w'_{(a_{4})},\tilde{w}_{(a_{2})}, \tilde{w}_{(a_{1})},
\tilde{w}_{(a_{3})}, \mathcal{B}([\mathcal{Z}]_P); x_{1}, x_{2})}\nno\\
&&\quad -x_{0}^{-1}\delta\left(\frac{x_{2}-x_{1}}{-x_{0}}\right)
\iota_{21} G_{\alpha}(w'_{(a_{4})},\tilde{w}_{(a_{2})}, \tilde{w}_{(a_{1})},
\tilde{w}_{(a_{3})}, \mathcal{B}([\mathcal{Z}]_P); x_{1}, x_{2})
\nno\\
&&=x_{2}^{-1}\delta\left(\frac{x_{1}-x_{0}}{x_{2}}\right)
\iota_{20} G_{\alpha}(w'_{(a_{4})},\tilde{w}_{(a_{2})}, \tilde{w}_{(a_{1})},
\tilde{w}_{(a_{3})}, \mathcal{B}([\mathcal{Z}]_P); x_2+x_0, x_{2})\label{e1:33}
\end{eqnarray}
for $\alpha\in \mathbb{ A}(a_{2}, a_{1}, a_{3}, a_{4})$.
Since $w_{(a_4)}'\in (W^{a_4})'$ is arbitrary, by (\ref{e1:29})-(\ref{e1:31}) and (\ref{e1:33}), the Jacobi identity holds for the ordered triple $(\tilde{w}_{(a_{2})},\tilde{w}_{(a_{1})},
\tilde{w}_{(a_{3})})$:
\begin{eqnarray}
\lefteqn{x_{0}^{-1}
\delta\left(\frac{x_{1}-x_{2}}{x_{0}}\right)
f^{a_{2}, a_{1}, a_{3},
a_{4}}_{\alpha}(\tilde{w}_{(a_{2})},
\tilde{w}_{(a_{1})}, \tilde{w}_{(a_{3})}, \mathcal{B}([\mathcal{Z}]_P); x_{1}, x_{2})}\nno\\
&&\quad -x_{0}^{-1}\delta\left(\frac{x_{2}-x_{1}}{-x_{0}}\right)
g^{a_2, a_1, a_{3},
a_{4}}_{\alpha}(\tilde{w}_{(a_2)},
\tilde{w}_{(a_1)}, \tilde{w}_{(a_{3})}, \mathcal{ B}(\mathcal{B}([\mathcal{Z}]_P)); x_1, x_2)
\nno\\
&&=x_{2}^{-1}\delta\left(\frac{x_{1}-x_{0}}{x_{2}}\right)
h^{a_{2}, a_{1}, a_{3},
a_{4}}_{\alpha}(\tilde{w}_{(a_{2})},
\tilde{w}_{(a_{1})}, \tilde{w}_{(a_{3})}, \mathcal{ F}(\mathcal{B}([\mathcal{Z}]_P)); x_{0}, x_{2})\label{a20}
\end{eqnarray}
for $\mathcal{ Z}\in \coprod_{ a_{5}\in \mathcal{ A}}
\mathcal{ V}_{a_{1}a_{5}}^{a_{4}}\otimes
\mathcal{ V}_{a_{2}a_{3}}^{a_{5}}$ and $\alpha\in \mathbb{ A}(a_{2}, a_{1}, a_{3}, a_{4})$. Since $\mathcal{B}$ is isomorphic, by (\ref{a19}) we see that $\mathcal{B}([\mathcal{Z}]_P)$ in (\ref{a37})-(\ref{a20}) can be any element in $\pi_P(\coprod_{a_{5}\in \mathcal{ A}}
\mathcal{ V}_{a_{2}a_{5}}^{a_{4}}\otimes
\mathcal{ V}_{a_{1}a_{3}}^{a_{5}})$. So the Jacobi identity holds for the ordered triple
$(\tilde{w}_{(a_{2})}, \tilde{w}_{(a_{1})}, \tilde{w}_{(a_{3})})$.\epf
\vspace{0.2cm}

\begin{thm}\label{l3}
Assume that the assumptions of Theorem \ref{t1} hold, then the Jacobi identity holds for the ordered triple
$(\tilde{w}_{(a_{1})}, \tilde{w}_{(a_{3})}, \tilde{w}_{(a_{2})})$.
\end{thm}

\pf
Consider any $a_4\in \mathcal{ A}$, $w_{(a_4)}'\in (W^{a_4})'$ and $\mathcal{ Z}\in \coprod_{ a_{5}\in \mathcal{ A}}
\mathcal{ V}_{a_{1}a_{5}}^{a_{4}}\otimes
\mathcal{ V}_{a_{2}a_{3}}^{a_{5}}$.
Observe that
$$e^{-x_2 L(1)}w_{(a_4)}'\in (W^{a_4})'[x_2].$$
Then it can be easily derived that Proposition \ref{l1} holds with $w_{(a_4)}'$ replaced by $e^{-x_2 L(1)}w_{(a_4)}'$. In particular, replacing $w_{(a_4)}'$
by $e^{-x_2 L(1)}w_{(a_4)}'$ in (\ref{s7}), we get a multivalued analytic function
\begin{equation}\label{a35}
\Phi(e^{-x_2 L(1)}w_{(a_4)}', \tilde{w}_{(a_{1})},
\tilde{w}_{(a_{2})}, \tilde{w}_{(a_{3})}, [\mathcal{Z}]_P;z_1,z_2) \in \mathbb{ G}^{a_{1}, a_{2}, a_{3}, a_{4}}.
\end{equation}
Moreover, replacing the complex variables $(z_1,z_2)$ by $(z_1-z_2,-z_2)$ in (\ref{a35}), we get a multivalued analytic function on $M^{2}=\{(z_{1},
z_{2})\in \mathbb{ C}^{2}\;|\;z_{1}, z_{2}\ne 0, z_{1}\ne
z_{2}\}$, which shall simply be denoted by
\begin{equation}\label{b17}
\Phi(e^{-x_2 L(1)}w_{(a_4)}', \tilde{w}_{(a_{1})},
\tilde{w}_{(a_{2})}, \tilde{w}_{(a_{3})}, [\mathcal{Z}]_P;z_1-z_2,-z_2).
\end{equation}

Consider the simply connected region in $\mathbb{ C}^{2}$
obtained by cutting the region $|z_{1}-z_2|>|z_{2}|>0$ along the
intersection of this region with $$\{(z_{1}, z_{2})\in \C^{2}\;|\;z_{2}\in (-\infty,0]\}\cup \{(z_{1}, z_{2})\in \C^{2}\;|\;
z_{1}-z_{2}\in [0,+\infty)\}.$$ We denote it by $R_{5}$.
Replacing $w_{(a_4)}'$
by $e^{-x_2 L(1)}w_{(a_4)}'$, and then $(z_1,z_2)$ by $(z_1-z_2,-z_2)$ in (\ref{s8}), we see that
\begin{equation}\label{a33}
\langle e^{-x_2 L(1)}w_{(a_4)}', (\tilde{\bf P}([\mathcal{Z}]_P))(\tilde{w}_{(a_{1})},
\tilde{w}_{(a_{2})}, \tilde{w}_{(a_{3})}; x_0, x_{2}) \rangle_{W^{a_4}}\mbar_{x_{0}^{n}=e^{n \log (z_1-z_2)},x_{2}^{n}=e^{n \log (-z_2)}}
\end{equation}
is a branch of $\Phi(e^{-x_2 L(1)}w_{(a_4)}', \tilde{w}_{(a_{1})},
\tilde{w}_{(a_{2})}, \tilde{w}_{(a_{3})}, [\mathcal{Z}]_P;z_1-z_2,-z_2)
$ on the region $R_5$.
Moreover, the skew-symmetry isomorphism implies
\begin{eqnarray}
\lefteqn{ \langle w_{(a_4)}', (\tilde{\bf P}(\tilde{\Omega}^{(4)}([\mathcal{Z}]_P)))(\tilde{w}_{(a_{1})},
\tilde{w}_{(a_{3})}, \tilde{w}_{(a_{2})}; x_{1},x_{2}) \rangle_{W^{a_4}}\mbar_{x_{1}^{n}=e^{n \log z_1},x_{2}^{n}=e^{n \log z_2}}}\nno\\
&& =\langle e^{x_2 L(1)}w_{(a_4)}', (\tilde{\bf P}([\mathcal{Z}]_P))(\tilde{w}_{(a_{1})},
\tilde{w}_{(a_{2})}, \tilde{w}_{(a_{3})}; x_0, e^{-\pi i}x_{2}) \rangle_{W^{a_4}}\mbar_{x_{0}^{n}=e^{n \log (z_1-z_2)},x_{2}^{n}=e^{n \log z_2}}\nno\\
&& =\langle e^{-x_2 L(1)}w_{(a_4)}', (\tilde{\bf P}([\mathcal{Z}]_P))(\tilde{w}_{(a_{1})},
\tilde{w}_{(a_{2})}, \tilde{w}_{(a_{3})}; x_0, x_{2}) \rangle_{W^{a_4}}\mbar_{x_{0}^{n}=e^{n \log (z_1-z_2)},x_{2}^{n}=e^{n \log (-z_2)}}\qquad\quad
\end{eqnarray}
on the region $\{(z_{1}, z_{2})\in \C^{2}\;|\;
\re z_1>-\re z_2>0,
\im z_1>-\im z_2>0 \}$.
And observing that the single-valued analytic function
\begin{equation}\label{a24}
\langle w_{(a_4)}', (\tilde{\bf P}(\tilde{\Omega}^{(4)}([\mathcal{Z}]_P)))(\tilde{w}_{(a_{1})},
\tilde{w}_{(a_{3})}, \tilde{w}_{(a_{2})}; x_{1},x_{2}) \rangle_{W^{a_4}}\mbar_{x_{1}^{n}=e^{n \log z_1},x_{2}^{n}=e^{n \log z_2}}
\end{equation}
on the region $\{(z_{1}, z_{2})\in \C^{2}\;|\;
\re z_1>-\re z_2>0,
\im z_1>-\im z_2>0 \}$ can be naturally analytically extended to
the region $R_1$, we can conclude that the single-valued analytic function (\ref{a24}) on $R_1$ is a branch of $\Phi(e^{-x_2 L(1)}w_{(a_4)}', \tilde{w}_{(a_{1})},
\tilde{w}_{(a_{2})}, \tilde{w}_{(a_{3})}, [\mathcal{Z}]_P;z_1-z_2,-z_2)$ on $R_1$.
\vspace{0.2cm}

Let $(a_0,b_0)$, $(a_1,b_1)$, $(a_2,b_2)$ and $(a_3,b_3)$ be four pairs of fixed positive real numbers satisfying
\begin{eqnarray}
a_0>b_0>a_0-b_0>0,\ \
a_1>a_1-b_1>b_1>0,\nno\\
b_2>b_2-a_2>a_2>0,\ \ b_3>a_3>b_3-a_3>0.\label{b10}
\end{eqnarray}
Then we shall obtain branches of $\Phi(e^{-x_2 L(1)}w_{(a_4)}', \tilde{w}_{(a_{1})},
\tilde{w}_{(a_{2})}, \tilde{w}_{(a_{3})}, [\mathcal{Z}]_P;z_1-z_2,-z_2)$
by analytical continuations along curves.

First of all, we consider the simply connected region
\begin{eqnarray*}
\mathfrak{G}' &=&\mathbb{ C}^{2}\backslash \big(\{(z_{1}, z_{2})\in \C^{2}\;|\; z_1\in [0,+\infty)\}\cup \{(z_{1}, z_{2})\in \C^{2}\;|\; z_2\in [0,+\infty) \}\cup \\
 && \qquad\qquad \{(z_{1}, z_{2})\in \C^{2}\;|\;
z_{1}-z_{2}\in [0,+\infty)\}\big).
\end{eqnarray*}
Define a path $\gamma:\ [0,1]\rightarrow \mathfrak{G}'$ by
\begin{eqnarray}
\lefteqn{\gamma(t)=\left(\tilde{z}_1(t), \tilde{z}_2(t)\right)}\nno\\
&&=\left\{\begin{array}{ll}
\left( (a_0(1-7t)+7a_1t)e^{\frac{1}{4}\pi i}, \;\; (b_0(1-7t)+7b_1t)e^{\frac{1}{4}\pi i} \right)&
t\in[0,\frac{1}{7}],\\
\left( a_1 e^{\frac{1}{4}\pi i}, \;\; b_1 e^{\frac{1}{4}\pi i+(7t-1)\pi i} \right)&
t\in (\frac{1}{7},\frac{2}{7}],\\
\left( (a_1(3-7t)+a_2(7t-2)) e^{\frac{1}{4}\pi i}, \;\; (b_1(3-7t)+b_2(7t-2))e^{\frac{5}{4}\pi i} \right)&
t\in (\frac{2}{7},\frac{3}{7}],\qquad\\
\left( a_2 e^{\frac{1}{4}\pi i+(7t-3)\pi i}, \;\; b_2 e^{\frac{5}{4}\pi i} \right)&
t\in(\frac{3}{7},\frac{4}{7}],\\
\left( b_2 e^{\frac{5}{4}\pi i}+(b_2-a_2) e^{\frac{1}{4}\pi i+(7t-4)\pi i}, \;\; b_2 e^{\frac{5}{4}\pi i} \right)&
t\in(\frac{4}{7},\frac{5}{7}],\\
\left((2b_2-a_2)e^{\frac{5}{4}\pi i-(7t-5)\pi i}, \;\; b_2 e^{\frac{5}{4}\pi i-(7t-5)\pi i} \right)&
t\in(\frac{5}{7},\frac{6}{7}],\\
\Big( (7(2b_2-a_2)(1-t)+a_0(7t-6))e^{\frac{1}{4}\pi i}, \\
\quad\qquad (7b_2(1-t)+b_0(7t-6))e^{\frac{1}{4}\pi i} \Big)&
t\in(\frac{6}{7},1].
\end{array}\right.
\end{eqnarray}
\begin{figure}[!htb]
\caption[b]{$\gamma(t)$}
$\\$
\centering
\resizebox{11cm}{4.8cm}{\includegraphics{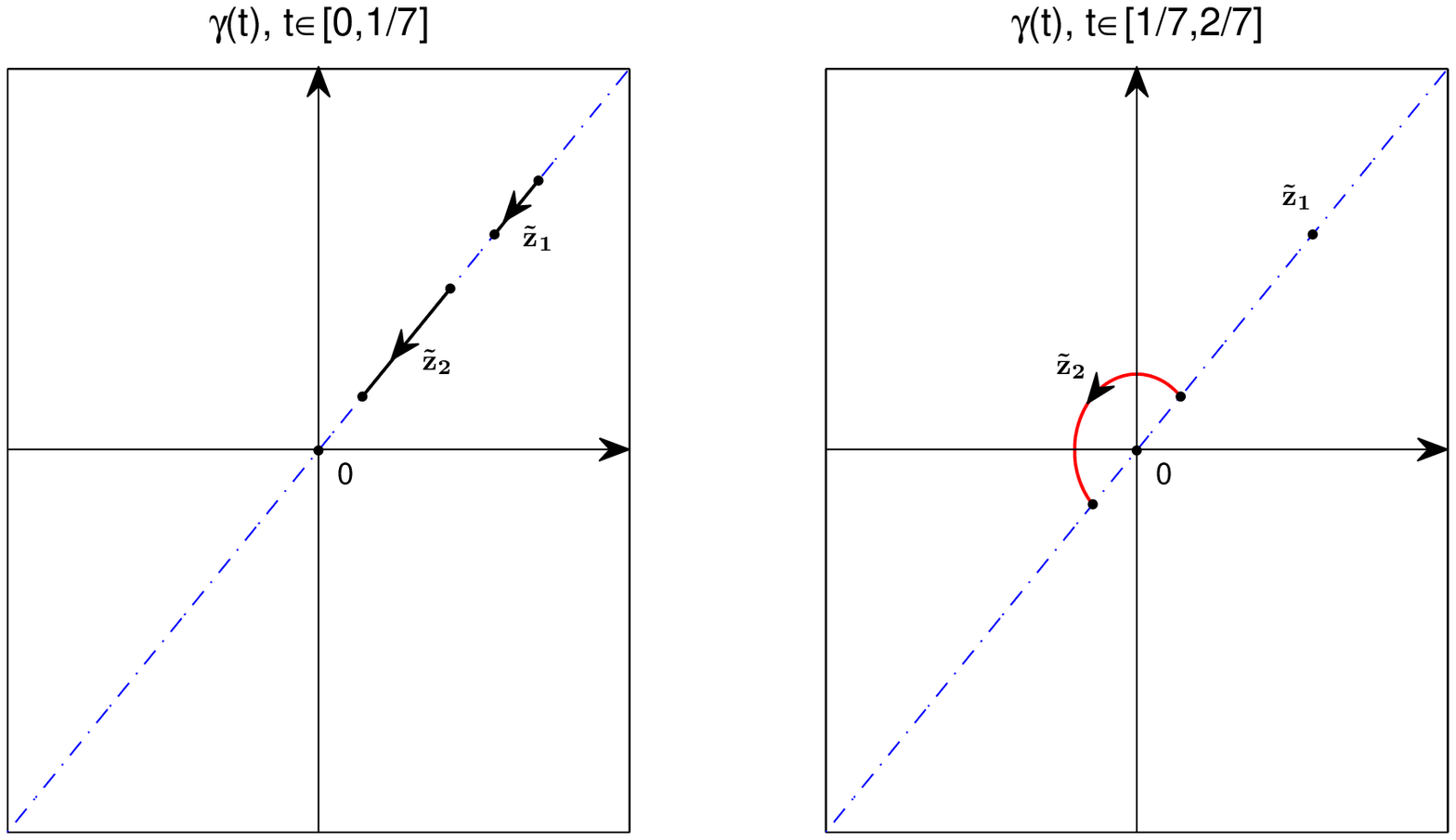}}

\resizebox{11cm}{4.8cm}{\includegraphics{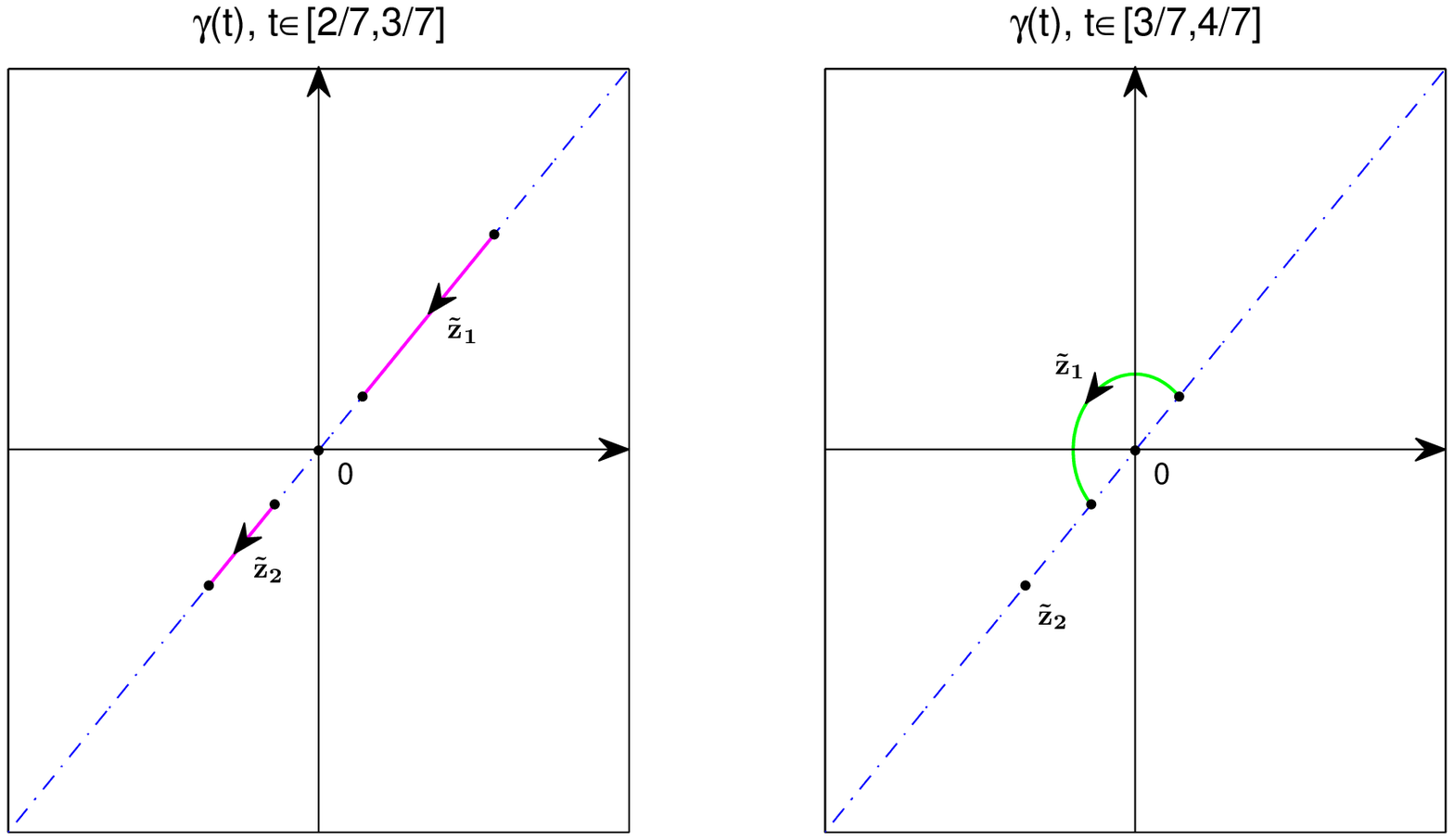}}

\resizebox{11cm}{4.8cm}{\includegraphics{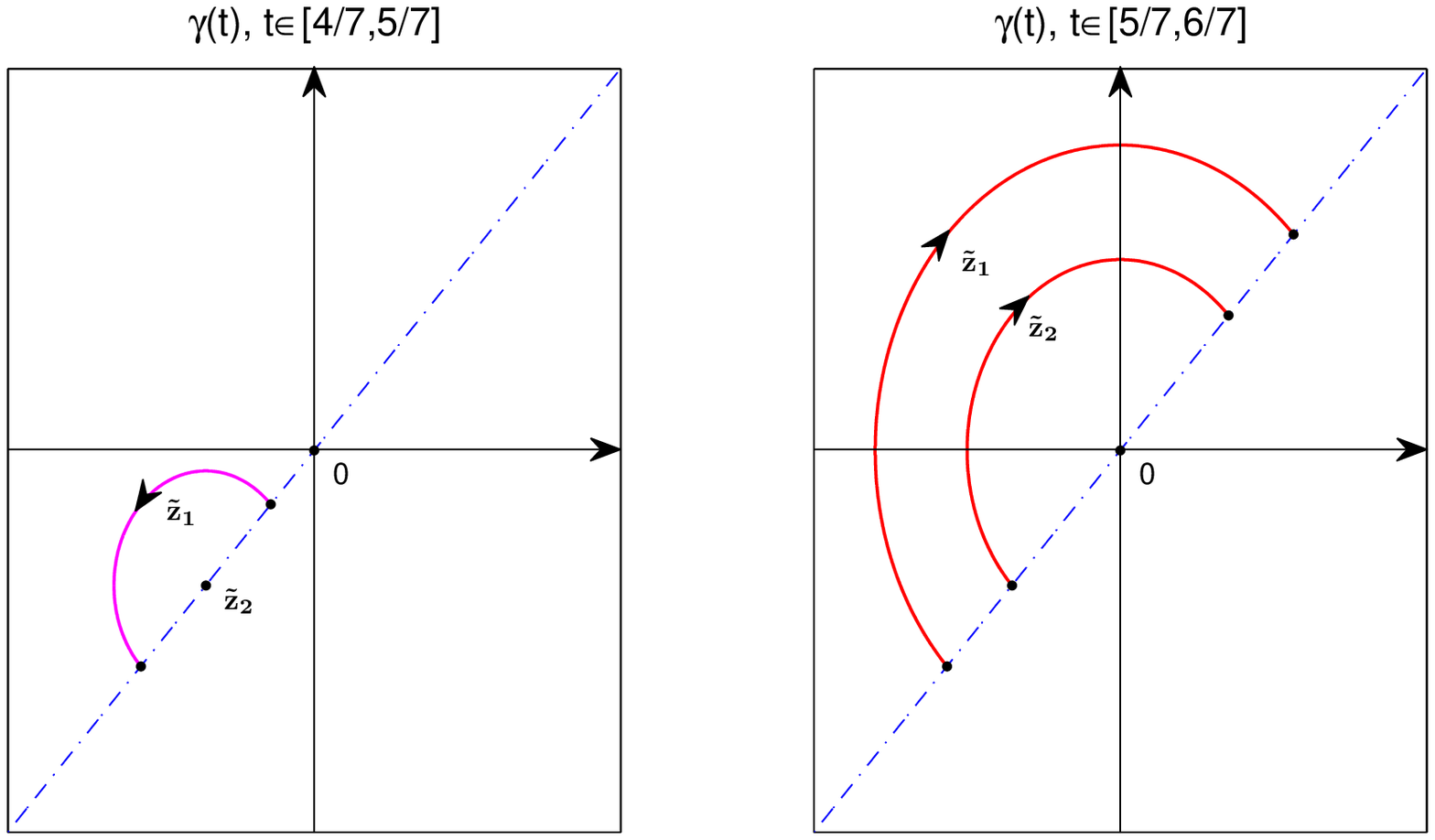}}

\resizebox{11cm}{4.8cm}{\includegraphics{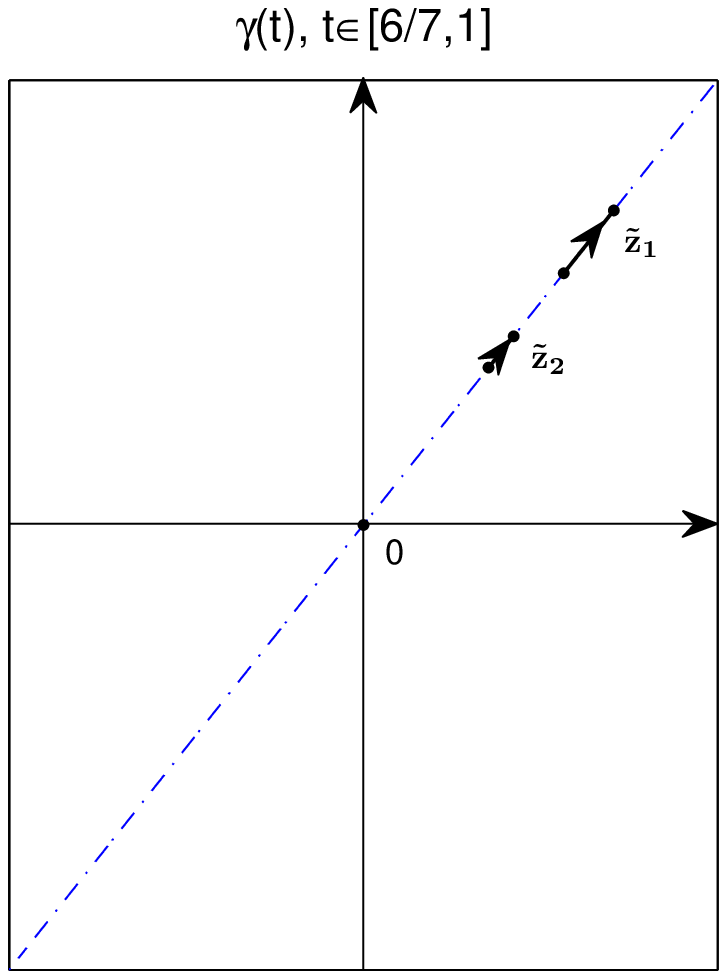}}
\end{figure}
See Figure 1 for an illustration. Then $\gamma(t)\subset \mathfrak{G}'$. We choose a simply connected region
\begin{equation}
D_t=\{(z_{1}, z_{2})\in \C^{2}\;|\; \textrm{max}(|z_1-\tilde{z}_1(t)|, |z_2-\tilde{z}_2(t)|)< \varepsilon_t\}
\end{equation}
for each $t\in [0,1]$, where $\varepsilon_t$ is a sufficiently small positive real number for each $t\in [0,1]$ such that
\begin{equation*}
D_0 \subset \mathfrak{G}'\cap  \{(z_1,z_2)\in \mathbb{C}^2 \mid \re z_1>\re z_2 > \re (z_1-z_2) >0,\ \im z_1 >  \im z_2 > \im (z_1-z_2) > 0\},
\end{equation*}
\begin{equation*}
D_t \subset \mathfrak{G}'\cap \{(z_{1}, z_{2})\in \C^{2}\;|\; |z_1|>|z_2|>0\} \ \ \textrm{ for } t\in (0,\frac{2}{7}),
\end{equation*}
\begin{equation*}
D_{\frac{2}{7}} \subset \mathfrak{G}'\cap  \{(z_1,z_2)\in \mathbb{C}^2 \mid \re z_1>-\re z_2 >0,\ \im z_1 > -\im z_2 > 0\},
\end{equation*}
\begin{equation*}
D_t \subset \mathfrak{G}'\cap \{(z_{1}, z_{2})\in \C^{2}\mid \re z_1>0>\re z_2,\ \im z_1 > 0>\im z_2\} \ \ \textrm{ for } t\in (\frac{2}{7},\frac{3}{7}),
\end{equation*}
\begin{equation*}
D_{\frac{3}{7}} \subset \mathfrak{G}'\cap  \{(z_1,z_2)\in \mathbb{C}^2 \mid -\re z_2 >\re z_1>0,\ -\im z_2 > \im z_1 > 0\},
\end{equation*}
\begin{equation*}
D_t \subset \mathfrak{G}'\cap \{(z_{1}, z_{2})\in \C^{2}\mid |z_2|>|z_1|>0\} \ \ \textrm{ for } t\in (\frac{3}{7},\frac{4}{7}),
\end{equation*}
\begin{equation*}
D_{\frac{4}{7}} \subset \mathfrak{G}'\cap  \{(z_1,z_2)\in \mathbb{C}^2 \mid \re z_2 <\re(z_2-z_1)<\re z_1<0,\ \im z_2 <\im (z_2-z_1)<\im z_1 < 0\},
\end{equation*}
\begin{equation*}
D_t \subset \mathfrak{G}'\cap \{(z_{1}, z_{2})\in \C^{2}\mid |z_2|>|z_1-z_2|>0\} \ \ \textrm{ for } t\in (\frac{4}{7},1),
\end{equation*}
\begin{equation*}
D_1 = D_0.
\end{equation*}
With some straightforward calculations, the existence of $\varepsilon_t$ can be easily verified. We omit the details here except that we shall write more about $\varepsilon_t$ for $t\in (\frac{3}{7},\frac{4}{7})$. Note that $|\tilde{z}_1(t)|=a_2$ and $|\tilde{z}_2(t)|=b_2$ for $t\in (\frac{3}{7},\frac{4}{7})$. So for each $t\in (\frac{3}{7},\frac{4}{7})$, to ensure that $D_t \subset \mathfrak{G}'$, we must have $\varepsilon_t<a_2$, which further implies $\varepsilon_t<\frac{1}{2}b_2$ by (\ref{b10}).
Thus $\re z_2<0$ and $\im z_2<0$ for any $(z_{1}, z_{2})\in D_t$ with $t\in (\frac{3}{7},\frac{4}{7})$. With these simply connected regions, we can see that
\begin{equation}
f_t=\langle w_{(a_4)}', (\tilde{\bf P}(\tilde{\Omega}^{(4)}([\mathcal{Z}]_P)))(\tilde{w}_{(a_{1})},
\tilde{w}_{(a_{3})}, \tilde{w}_{(a_{2})}; x_{1},x_{2}) \rangle_{W^{a_4}}\mbar_{x_{1}^{n}=e^{n \log z_1},x_{2}^{n}=e^{n \log z_2}}
\end{equation}
is a single-valued analytic function on the region $D_t$ for each $t\in [0,\frac{2}{7}]$;
\begin{equation}
f_t=\langle e^{-x_2 L(1)}w_{(a_4)}', (\tilde{\bf P}([\mathcal{Z}]_P))(\tilde{w}_{(a_{1})},
\tilde{w}_{(a_{2})}, \tilde{w}_{(a_{3})}; x_0, x_{2}) \rangle_{W^{a_4}}\mbar_{x_{0}^{n}=e^{n \log (z_1-z_2)},x_{2}^{n}=e^{n \log (-z_2)}}
\end{equation}
is a single-valued analytic function on the region $D_t$ for each $t\in (\frac{2}{7},\frac{3}{7}]$;
\begin{equation}
f_t=\langle e^{-x_2 L(1)}w_{(a_4)}', (\tilde{\bf I}(\mathcal{F}([\mathcal{Z}]_P)))(\tilde{w}_{(a_{1})},
\tilde{w}_{(a_{2})}, \tilde{w}_{(a_{3})}; x_1, x_2) \rangle_{W^{a_4}}\mbar_{x_{1}^{n}=e^{n \log z_1},x_{2}^{n}=e^{n \log (-z_2)}}
\end{equation}
is a single-valued analytic function on the region $D_t$ for each $t\in (\frac{3}{7},\frac{4}{7}]$; and
\begin{equation}
f_t=\langle w_{(a_4)}', (\tilde{\bf I}(\mathcal{F}(\tilde{\Omega}^{(4)}([\mathcal{Z}]_P))))(\tilde{w}_{(a_{1})},
\tilde{w}_{(a_{3})}, \tilde{w}_{(a_{2})}; x_0,x_{2}) \rangle_{W^{a_4}}\mbar_{x_{0}^{n}=e^{n \log (z_1-z_2)}, x_{2}^{n}=e^{n \log z_2}}
\end{equation}
is a single-valued analytic function on the region $D_t$ for each $t\in (\frac{4}{7},1]$.
Next, we shall show that $\{(f_t,D_t): 0\leq t \leq 1 \}$ is an analytic continuation along $\gamma$.

Firstly, it can be derived from the skew-symmetry property that on the region $D_{\frac{2}{7}}$,
\begin{eqnarray}
\lefteqn{ f_{\frac{2}{7}}=\langle w_{(a_4)}', (\tilde{\bf P}(\tilde{\Omega}^{(4)}([\mathcal{Z}]_P)))(\tilde{w}_{(a_{1})},
\tilde{w}_{(a_{3})}, \tilde{w}_{(a_{2})}; x_{1},x_{2}) \rangle_{W^{a_4}}\mbar_{x_{1}^{n}=e^{n \log z_1},x_{2}^{n}=e^{n \log z_2}}}\nno\\
&& =\langle e^{x_2 L(1)}w_{(a_4)}', (\tilde{\bf P}([\mathcal{Z}]_P))(\tilde{w}_{(a_{1})},
\tilde{w}_{(a_{2})}, \tilde{w}_{(a_{3})}; x_0, e^{-\pi i}x_{2}) \rangle_{W^{a_4}}\mbar_{x_{0}^{n}=e^{n \log (z_1-z_2)},x_{2}^{n}=e^{n \log z_2}}\nno\\
&& =\langle e^{-x_2 L(1)}w_{(a_4)}', (\tilde{\bf P}([\mathcal{Z}]_P))(\tilde{w}_{(a_{1})},
\tilde{w}_{(a_{2})}, \tilde{w}_{(a_{3})}; x_0, x_{2}) \rangle_{W^{a_4}}\mbar_{x_{0}^{n}=e^{n \log (z_1-z_2)},x_{2}^{n}=e^{n \log (-z_2)}}.\qquad\quad
\end{eqnarray}
Secondly, replacing $(w'_{(a_{4})},z_1,z_2)$ by $(e^{-x_2 L(1)}w_{(a_4)}',z_1-z_2,-z_2)$ in (\ref{s12}), we can derive that
\begin{eqnarray}
\lefteqn{f_{\frac{3}{7}}=\langle e^{-x_2 L(1)}w_{(a_4)}', (\tilde{\bf P}([\mathcal{Z}]_P))(\tilde{w}_{(a_{1})},
\tilde{w}_{(a_{2})}, \tilde{w}_{(a_{3})}; x_0, x_{2}) \rangle_{W^{a_4}}\mbar_{x_{0}^{n}=e^{n \log (z_1-z_2)},x_{2}^{n}=e^{n \log (-z_2)}}}\nno\\
&& =\langle e^{-x_2 L(1)}w_{(a_4)}', (\tilde{\bf I}(\mathcal{F}([\mathcal{Z}]_P)))(\tilde{w}_{(a_{1})},
\tilde{w}_{(a_{2})}, \tilde{w}_{(a_{3})}; x_1, x_2) \rangle_{W^{a_4}}\mbar_{x_{1}^{n}=e^{n \log z_1},x_{2}^{n}=e^{n \log (-z_2)}}\qquad\quad
\end{eqnarray}
on the region $D_{\frac{3}{7}}$. Thirdly, we shall prove that on the region $D_{\frac{4}{7}}$,
\begin{eqnarray}
\lefteqn{f_{\frac{4}{7}}=\langle e^{-x_2 L(1)}w_{(a_4)}', (\tilde{\bf I}(\mathcal{F}([\mathcal{Z}]_P)))(\tilde{w}_{(a_{1})},
\tilde{w}_{(a_{2})}, \tilde{w}_{(a_{3})}; x_1, x_2) \rangle_{W^{a_4}}\mbar_{x_{1}^{n}=e^{n \log z_1},x_{2}^{n}=e^{n \log (-z_2)}}}\nno\\
&& =\langle w_{(a_4)}', (\tilde{\bf I}(\mathcal{F}(\tilde{\Omega}^{(4)}([\mathcal{Z}]_P))))(\tilde{w}_{(a_{1})},
\tilde{w}_{(a_{3})}, \tilde{w}_{(a_{2})}; x_0,x_{2}) \rangle_{W^{a_4}}\mbar_{x_{0}^{n}=e^{n \log (z_1-z_2)}, x_{2}^{n}=e^{n \log z_2}}.\qquad\quad \label{b11}
\end{eqnarray}
By skew-symmetry, we have
\begin{eqnarray}
\lefteqn{\langle e^{-x_2 L(1)}w_{(a_4)}', (\tilde{\bf I}((\widetilde{\Omega^{-1}})^{(1)}\mathcal{F}([\mathcal{Z}]_P)))(\tilde{w}_{(a_{2})},
\tilde{w}_{(a_{1})}, \tilde{w}_{(a_{3})}; x_1, x_0) \rangle_{W^{a_4}}\mbar_{x_{1}^{n}=e^{n \log (-z_1)},x_{0}^{n}=e^{n \log (z_1-z_2)}}}\nno\\
&& =\langle e^{-x_2 L(1)}w_{(a_4)}', (\tilde{\bf I}(\mathcal{F}([\mathcal{Z}]_P)))(\tilde{w}_{(a_{1})},
\tilde{w}_{(a_{2})}, \tilde{w}_{(a_{3})}; e^{\pi i}x_1, x_2) \rangle_{W^{a_4}}\mbar_{x_{1}^{n}=e^{n \log (-z_1)},x_{2}^{n}=e^{n \log (-z_2)}}\nno\\
&& =\langle e^{-x_2 L(1)}w_{(a_4)}', (\tilde{\bf I}(\mathcal{F}([\mathcal{Z}]_P)))(\tilde{w}_{(a_{1})},
\tilde{w}_{(a_{2})}, \tilde{w}_{(a_{3})}; x_1, x_2) \rangle_{W^{a_4}}\mbar_{x_{1}^{n}=e^{n \log z_1},x_{2}^{n}=e^{n \log (-z_2)}}\qquad\quad\label{b12}
\end{eqnarray}
on the region $D_{\frac{4}{7}}$. Moreover, note that $[\mathcal{Z}]_P \in \pi_P(\coprod_{a_{5}\in \mathcal{ A}}
\mathcal{ V}_{a_{1}a_{5}}^{a_{4}}\otimes
\mathcal{ V}_{a_{2}a_{3}}^{a_{5}})$ implies
$$\mathcal{B}^{-1}(\mathcal{B}^{-1}([\mathcal{Z}]_P)) \in \pi_P(\coprod_{a_{5}\in \mathcal{ A}}
\mathcal{ V}_{a_{1}a_{5}}^{a_{4}}\otimes
\mathcal{ V}_{a_{2}a_{3}}^{a_{5}}).$$
So replacing $[\mathcal{Z}]_P$ by $\mathcal{B}^{-1}(\mathcal{B}^{-1}([\mathcal{Z}]_P))$, $(w'_{(a_{4})},z_1,z_2)$ by $(e^{-x_2 L(1)}w_{(a_4)}',z_1-z_2,-z_2)$ in (\ref{s13}),
 we obtain that
\begin{eqnarray}
\lefteqn{\langle e^{-x_2 L(1)}w_{(a_4)}', (\tilde{\bf P}(\mathcal{B}^{-1}([\mathcal{Z}]_P)))(\tilde{w}_{(a_{2})},
\tilde{w}_{(a_{1})}, \tilde{w}_{(a_{3})}; x_{2}, x_0) \rangle_{W^{a_4}}\lbar_{\substack{x_{0}^{n}=e^{n \log (z_1-z_2)}\\ x_{2}^{n}=e^{n \log (-z_2)}}}}\nno\\
&&=\langle e^{-x_2 L(1)}w_{(a_4)}', (\tilde{\bf I}((\widetilde{\Omega^{-1}})^{(1)}\mathcal{F}([\mathcal{Z}]_P)))(\tilde{w}_{(a_{2})},
\tilde{w}_{(a_{1})}, \tilde{w}_{(a_{3})}; x_1, x_0) \rangle_{W^{a_4}}\lbar_{\substack{x_{1}^{n}=e^{n \log (-z_1)}\\ x_{0}^{n}=e^{n \log (z_1-z_2)}}}\qquad\ \label{b13}
\end{eqnarray}
on the region $D_{\frac{4}{7}}$. Furthermore, by the Moore-Seiberg equations (\ref{hexagon1}) and
(\ref{hexagon2}), we have
\begin{equation}
\F \tilde{\Omega}^{(4)}
=\tilde{\Omega}^{(2)} \F^{-1} \F (\widetilde{\Omega^{-1}})^{(3)}\F \tilde{\Omega}^{(4)}
=\tilde{\Omega}^{(2)} \F^{-1} (\widetilde{\Omega^{-1}})^{(1)}\F (\widetilde{\Omega^{-1}})^{(4)} \tilde{\Omega}^{(4)}
=\tilde{\Omega}^{(2)}\mathcal{B}^{-1}.
\end{equation}
So this together with the skew-symmetry isomorphism implies
\begin{eqnarray}
\lefteqn{\langle w_{(a_4)}', (\tilde{\bf I}(\mathcal{F}(\tilde{\Omega}^{(4)}([\mathcal{Z}]_P))))(\tilde{w}_{(a_{1})},
\tilde{w}_{(a_{3})}, \tilde{w}_{(a_{2})}; x_0,x_{2}) \rangle_{W^{a_4}}\lbar_{\substack{x_{0}^{n}=e^{n \log (z_1-z_2)}\\ x_{2}^{n}=e^{n \log z_2}}}} \nno\\
&& =\langle w_{(a_4)}', (\tilde{\bf I}(\tilde{\Omega}^{(2)}\mathcal{B}^{-1}([\mathcal{Z}]_P)))(\tilde{w}_{(a_{1})},
\tilde{w}_{(a_{3})}, \tilde{w}_{(a_{2})}; x_0,x_{2}) \rangle_{W^{a_4}}\lbar_{\substack{x_{0}^{n}=e^{n\log (z_1-z_2)}\\ x_{2}^{n}=e^{n \log z_2}}} \nno\\
&& =\langle e^{x_2 L(1)}w_{(a_4)}', (\tilde{\bf P}(\mathcal{B}^{-1}([\mathcal{Z}]_P)))(\tilde{w}_{(a_{2})},
\tilde{w}_{(a_{1})}, \tilde{w}_{(a_{3})}; e^{-\pi i}x_{2}, x_0) \rangle_{W^{a_4}}\lbar_{\substack{x_{0}^{n}=e^{n \log (z_1-z_2)}\\ x_{2}^{n}=e^{n \log z_2}}}\nno\\
&& =\langle e^{-x_2 L(1)}w_{(a_4)}', (\tilde{\bf P}(\mathcal{B}^{-1}([\mathcal{Z}]_P)))(\tilde{w}_{(a_{2})},
\tilde{w}_{(a_{1})}, \tilde{w}_{(a_{3})}; x_{2}, x_0) \rangle_{W^{a_4}}\lbar_{\substack{x_{0}^{n}=e^{n \log (z_1-z_2)}\\ x_{2}^{n}=e^{n \log (-z_2)}}}\qquad\ \ \label{b14}
\end{eqnarray}
on the region $D_{\frac{4}{7}}$.
Therefore, (\ref{b11}) holds by (\ref{b12}), (\ref{b13}) and (\ref{b14}).

So to sum up, $\{(f_t,D_t): 0\leq t \leq 1 \}$ is an analytic continuation along $\gamma$.

Since $\mathfrak{G}'$ is simply connected, $\gamma\subset \mathfrak{G}'$ and $\gamma(0)=\gamma(1)$, we can derive that $f_0=f_1$ on the region $D_0\cap D_1=D_0$. Moreover, since $f_0$ and $f_1$ are both single-valued analytic functions on the domain $S_1$ which contains $D_0$, we deduce that $f_0=f_1$ on $S_1$.
Namely,
\begin{eqnarray}
\lefteqn{\langle w_{(a_4)}', (\tilde{\bf P}(\tilde{\Omega}^{(4)}([\mathcal{Z}]_P)))(\tilde{w}_{(a_{1})},
\tilde{w}_{(a_{3})}, \tilde{w}_{(a_{2})}; x_{1},x_{2}) \rangle_{W^{a_4}}\mbar_{x_{1}^{n}=e^{n \log z_1},x_{2}^{n}=e^{n \log z_2}}}\nno\\
&& =\langle w_{(a_4)}', (\tilde{\bf I}(\mathcal{F}(\tilde{\Omega}^{(4)}([\mathcal{Z}]_P))))(\tilde{w}_{(a_{1})},
\tilde{w}_{(a_{3})}, \tilde{w}_{(a_{2})}; x_0,x_{2}) \rangle_{W^{a_4}}\mbar_{x_{0}^{n}=e^{n \log (z_1-z_2)}, x_{2}^{n}=e^{n \log z_2}}\qquad\quad\label{b15}
\end{eqnarray}
on the region $S_1$. Furthermore, since the first line of (\ref{b15}) defined on $R_1$ is a branch of (\ref{b17}), and the second line of (\ref{b15}) on $S_1$ can be naturally analytically extended to the region $R_3$, we conclude that
\begin{equation}\label{b16}
\langle w_{(a_4)}', (\tilde{\bf I}(\mathcal{F}(\tilde{\Omega}^{(4)}([\mathcal{Z}]_P))))(\tilde{w}_{(a_{1})},
\tilde{w}_{(a_{3})}, \tilde{w}_{(a_{2})}; x_0,x_{2}) \rangle_{W^{a_4}}\mbar_{x_{0}^{n}=e^{n \log (z_1-z_2)}, x_{2}^{n}=e^{n \log z_2}}
\end{equation}
on $R_3$ is a branch of $\Phi(e^{-x_2 L(1)}w_{(a_4)}', \tilde{w}_{(a_{1})},
\tilde{w}_{(a_{2})}, \tilde{w}_{(a_{3})}, [\mathcal{Z}]_P;z_1-z_2,-z_2)$ on $R_3$.
\vspace{0.2cm}

Then, by skew-symmetry isomorphism and by (\ref{b15}), we can deduce that on the region $S_1$,
\begin{eqnarray}
\lefteqn{\langle w_{(a_4)}', (\tilde{\bf I}(\tilde{\Omega}^{(1)} \F (\tilde{\Omega}^{(4)}([\mathcal{Z}]_P))))(\tilde{w}_{(a_{3})},
\tilde{w}_{(a_{1})}, \tilde{w}_{(a_{2})};  x_{0},x_1) \rangle_{W^{a_4}}\mbar_{x_{0}^{n}=e^{n \log (z_2-z_1)},x_{1}^{n}=e^{n \log z_1}}}  \nno\\
&& =\langle w_{(a_4)}', (\tilde{\bf I}(\F (\tilde{\Omega}^{(4)}([\mathcal{Z}]_P))))(\tilde{w}_{(a_{1})},
\tilde{w}_{(a_{3})}, \tilde{w}_{(a_{2})};  e^{-\pi i}x_{0},x_2) \rangle_{W^{a_4}}\mbar_{x_{0}^{n}=e^{n \log (z_2-z_1)},x_{2}^{n}=e^{n \log z_2}}  \nno\\
&& =\langle w_{(a_4)}', (\tilde{\bf I}(\F (\tilde{\Omega}^{(4)}([\mathcal{Z}]_P))))(\tilde{w}_{(a_{1})},
\tilde{w}_{(a_{3})}, \tilde{w}_{(a_{2})};  x_{0},x_2) \rangle_{W^{a_4}}\mbar_{x_{0}^{n}=e^{n \log (z_1-z_2)},x_{2}^{n}=e^{n \log z_2}}\nno\\
&&=\langle w_{(a_4)}', (\tilde{\bf P}(\tilde{\Omega}^{(4)}([\mathcal{Z}]_P)))(\tilde{w}_{(a_{1})},
\tilde{w}_{(a_{3})}, \tilde{w}_{(a_{2})}; x_{1},x_{2}) \rangle_{W^{a_4}}\mbar_{x_{1}^{n}=e^{n \log z_1},x_{2}^{n}=e^{n \log z_2}}.\qquad\quad\label{b18}
\end{eqnarray}
Since the last line of (\ref{b18}) defined on $R_1$ is a branch of (\ref{b17}), and the first line of (\ref{b18}) on $S_1$ can be naturally analytically extended to the region $R_4$, we conclude that
\begin{equation}\label{b19}
\langle w_{(a_4)}', (\tilde{\bf I}(\tilde{\Omega}^{(1)} \F (\tilde{\Omega}^{(4)}([\mathcal{Z}]_P))))(\tilde{w}_{(a_{3})},
\tilde{w}_{(a_{1})}, \tilde{w}_{(a_{2})};  x_{0},x_1) \rangle_{W^{a_4}}\mbar_{x_{0}^{n}=e^{n \log (z_2-z_1)},x_{1}^{n}=e^{n \log z_1}}
\end{equation}
on $R_4$ is a branch of $\Phi(e^{-x_2 L(1)}w_{(a_4)}', \tilde{w}_{(a_{1})},
\tilde{w}_{(a_{2})}, \tilde{w}_{(a_{3})}, [\mathcal{Z}]_P;z_1-z_2,-z_2)$ on $R_4$.
\vspace{0.2cm}

Next, we consider the simply connected region
\begin{eqnarray*}
\mathfrak{G}'' &=&\mathbb{ C}^{2}\backslash \big(\{(z_{1}, z_{2})\in \C^{2}\;|\; z_1\in [0,+\infty)\} \cup \{(z_{1}, z_{2})\in \C^{2}\;|\; z_2\in [0,+\infty)\}\\
 && \qquad\qquad \cup \{(z_{1}, z_{2})\in \C^{2}\;|\;
z_{2}-z_{1}\in [0,+\infty)\}\big).
\end{eqnarray*}
Define a path $\sigma:\ [0,1]\rightarrow \mathfrak{G}''$ by
\begin{eqnarray}
\lefteqn{\sigma(t)=\left(\tilde{z}_1(t), \tilde{z}_2(t)\right)}\nno\\
&&=\left\{\begin{array}{ll}
\left( a_3e^{\frac{1}{4}\pi i}, \;\; a_3e^{\frac{1}{4}\pi i}+(b_3-a_3)e^{\frac{1}{4}\pi i+7t\pi i} \right)&
t\in[0,\frac{1}{7}],\\
\Big( (a_3(2-7t)+a_0(7t-1)) e^{\frac{1}{4}\pi i}, \\ \qquad\quad((2a_3-b_3)(2-7t)+b_0(7t-1))e^{\frac{1}{4}\pi i} \Big)&
t\in (\frac{1}{7},\frac{2}{7}],\\
\left( (a_0(3-7t)+a_1(7t-2)) e^{\frac{1}{4}\pi i}, \;\; (b_0(3-7t)+b_1(7t-2))e^{\frac{1}{4}\pi i} \right)&
t\in (\frac{2}{7},\frac{3}{7}],\qquad\\
\left( a_1 e^{\frac{1}{4}\pi i}, \;\; b_1 e^{\frac{1}{4}\pi i+(7t-3)\pi i} \right)&
t\in(\frac{3}{7},\frac{4}{7}],\\
\left( (a_1(5-7t)+a_2(7t-4)) e^{\frac{1}{4}\pi i}, \;\; (b_1(5-7t)+b_2(7t-4)) e^{\frac{5}{4}\pi i} \right)&
t\in(\frac{4}{7},\frac{5}{7}],\\
\left( a_2e^{\frac{1}{4}\pi i}, \;\; b_2 e^{\frac{5}{4}\pi i-(7t-5)\pi i} \right)&
t\in(\frac{5}{7},\frac{6}{7}],\\
\left( (7a_2(1-t)+a_3(7t-6))e^{\frac{1}{4}\pi i}, \;\;
 (7b_2(1-t)+b_3(7t-6))e^{\frac{1}{4}\pi i} \right)&
t\in(\frac{6}{7},1].
\end{array}\right.
\end{eqnarray}
\begin{figure}[!htb]
\caption[b]{$\sigma(t)$}
$\\$
\centering
\resizebox{11cm}{4.8cm}{\includegraphics{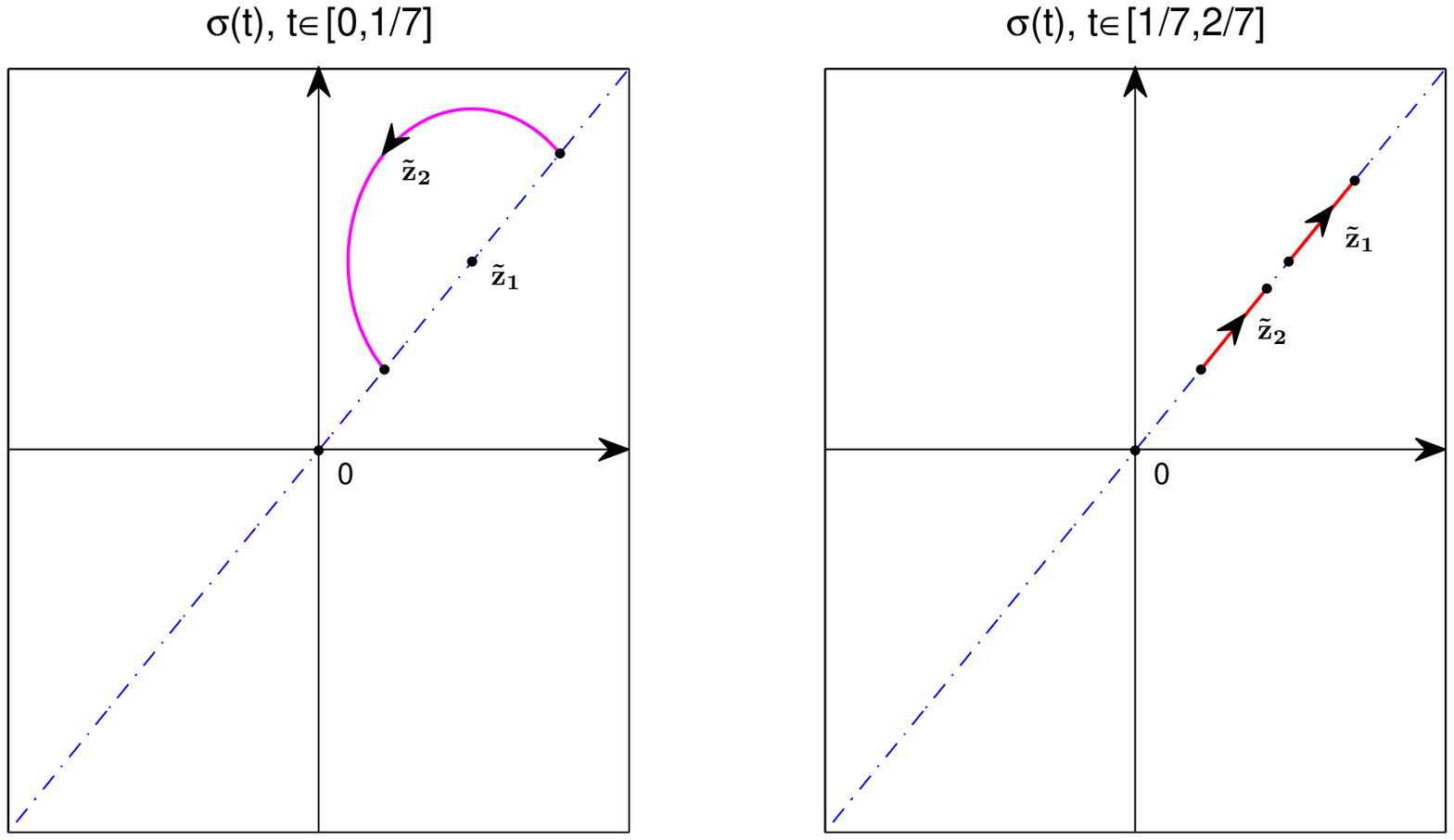}}

\resizebox{11cm}{4.8cm}{\includegraphics{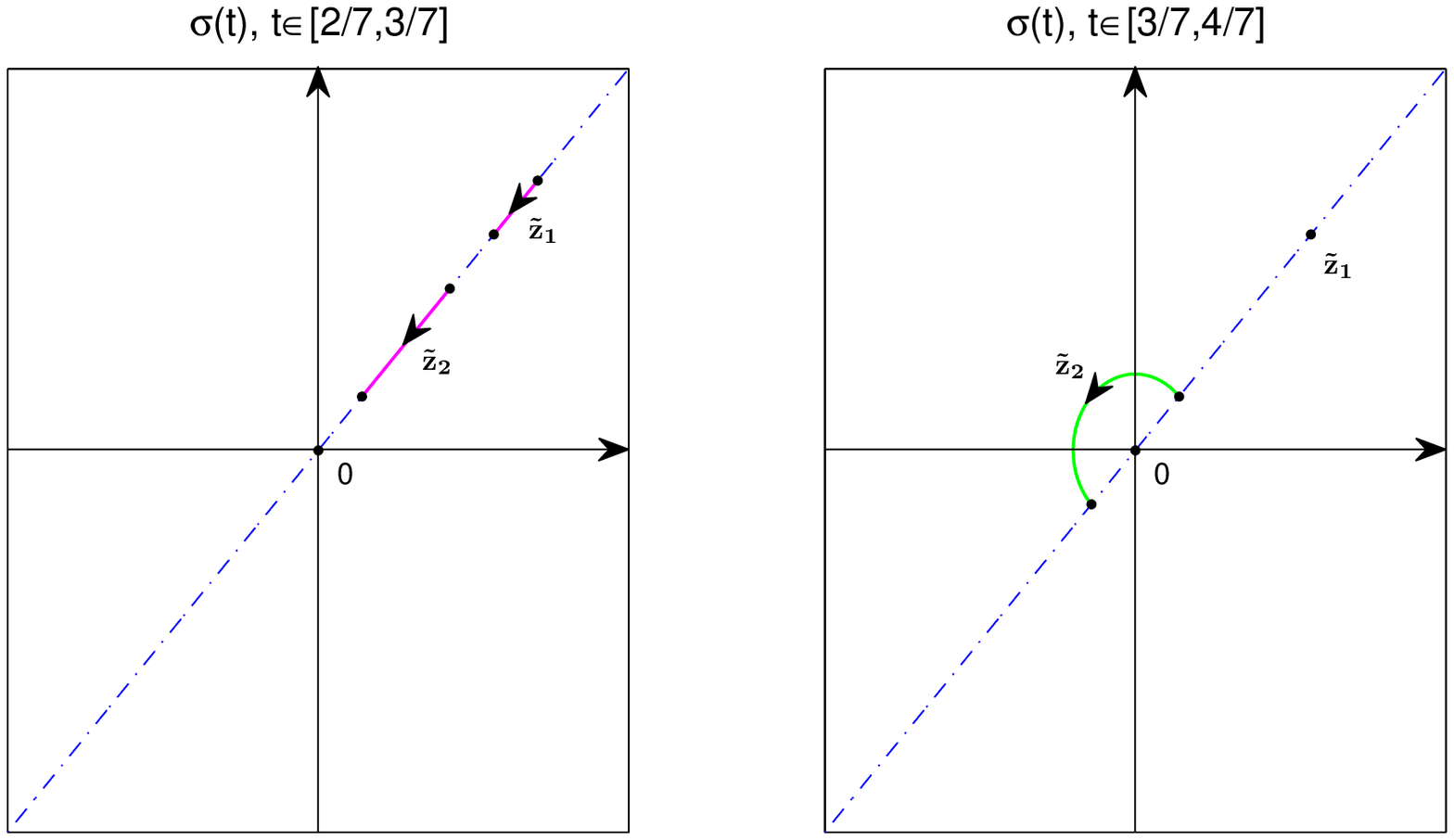}}

\resizebox{11cm}{4.8cm}{\includegraphics{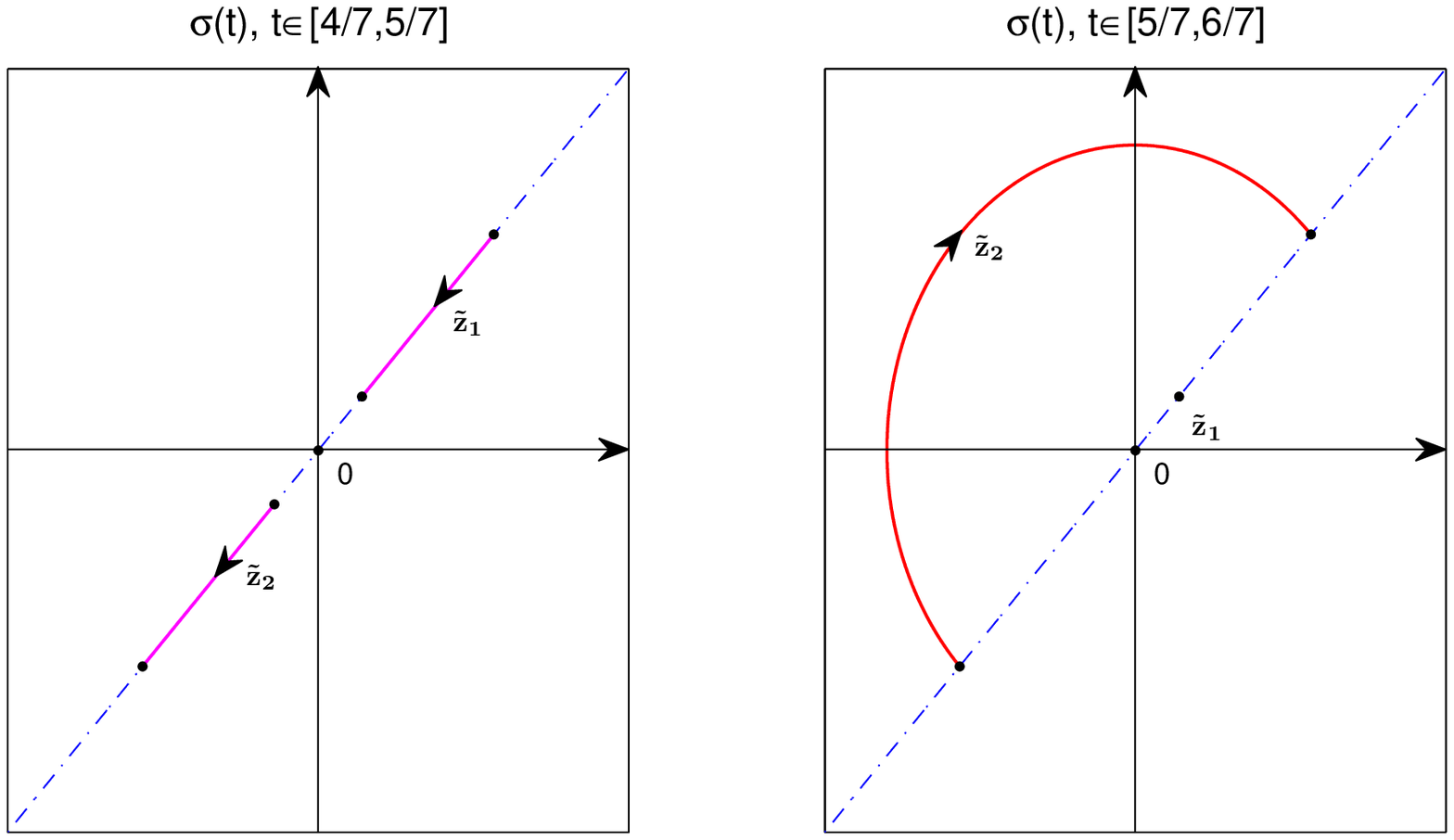}}

\resizebox{11cm}{4.8cm}{\includegraphics{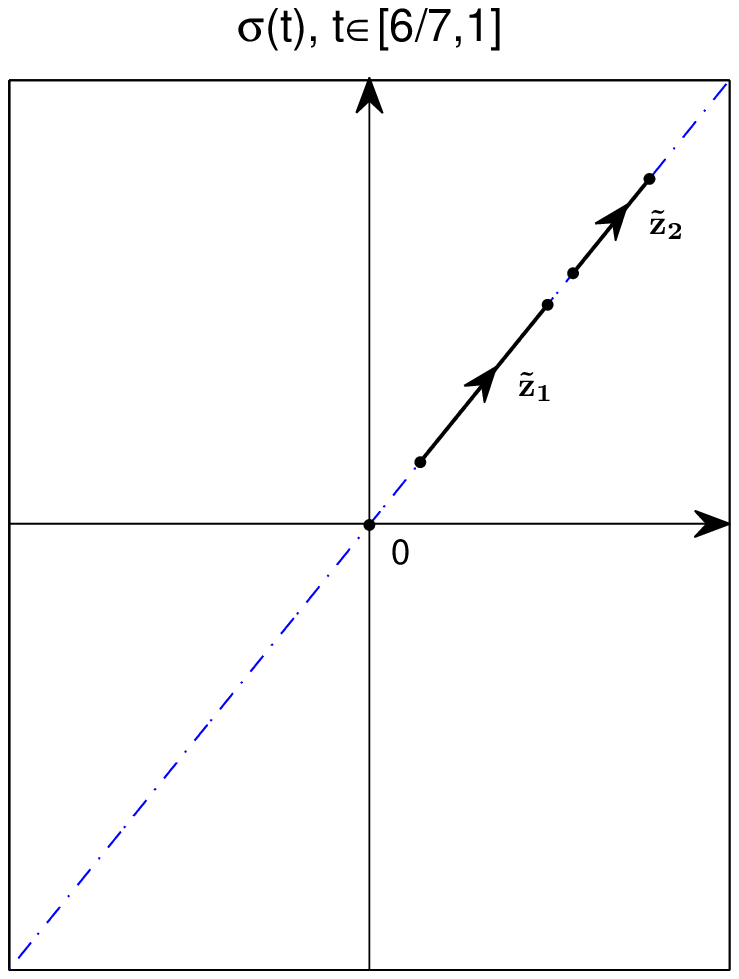}}
\end{figure}
See Figure 2 for an illustration. Then $\sigma(t)\subset \mathfrak{G}''$. We also choose a simply connected region
\begin{equation}
E_t=\{(z_{1}, z_{2})\in \C^{2}\;|\; \textrm{max}(|z_1-\tilde{z}_1(t)|, |z_2-\tilde{z}_2(t)|)< \epsilon_t\},
\end{equation}
for each $t\in [0,1]$, where $\epsilon_t$ is a sufficiently small positive real number for each $t\in [0,1]$ such that
\begin{equation*}
E_0 \subset \mathfrak{G}''\cap  \{(z_1,z_2)\in \mathbb{C}^2 \mid \re z_2>\re z_1 > \re (z_2-z_1) >0,\ \im z_2 >  \im z_1 > \im (z_2-z_1) > 0\},
\end{equation*}
\begin{equation*}
E_t \subset \mathfrak{G}''\cap \{(z_{1}, z_{2})\in \C^{2}\;|\; |z_1|>|z_1-z_2|>0\} \ \ \textrm{ for } t\in (0,\frac{2}{7}),
\end{equation*}
\begin{equation*}
E_{\frac{2}{7}} \subset \mathfrak{G}''\cap  S_1,
\end{equation*}
\begin{equation*}
E_t \subset \mathfrak{G}''\cap \{(z_{1}, z_{2})\in \C^{2}\mid |z_1|>|z_2|>0\} \ \ \textrm{ for } t\in (\frac{2}{7},\frac{4}{7}),
\end{equation*}
\begin{equation*}
E_{\frac{4}{7}} \subset \mathfrak{G}''\cap  \{(z_1,z_2)\in \mathbb{C}^2 \mid \re z_1>-\re z_2 >0,\  \im z_1 > -\im z_2 >0\},
\end{equation*}
\begin{equation*}
E_t \subset \mathfrak{G}''\cap \{(z_{1}, z_{2})\in \C^{2}\mid \re z_1 >0>\re z_2,\  \im z_1 >0> \im z_2 \} \ \ \textrm{ for } t\in (\frac{4}{7},\frac{5}{7}),
\end{equation*}
\begin{equation*}
E_{\frac{5}{7}} \subset \mathfrak{G}''\cap  \{(z_1,z_2)\in \mathbb{C}^2 \mid -\re z_2 >\re z_1>0,\ -\im z_2 > \im z_1 > 0\},
\end{equation*}
\begin{equation*}
E_t \subset \mathfrak{G}''\cap \{(z_{1}, z_{2})\in \C^{2}\mid |z_2|>|z_1|>0\} \ \ \textrm{ for } t\in (\frac{5}{7},1),
\end{equation*}
\begin{equation*}
E_1 = E_0.
\end{equation*}
Thus
\begin{equation}
g_t=\left.\langle w_{(a_4)}', (\tilde{\bf I}(\tilde{\Omega}^{(1)} \F (\tilde{\Omega}^{(4)}([\mathcal{Z}]_P))))(\tilde{w}_{(a_{3})},
\tilde{w}_{(a_{1})}, \tilde{w}_{(a_{2})};  x_{0},x_1) \rangle_{W^{a_4}}\right|_{\substack{x_{0}^{n}=e^{n \log (z_2-z_1)}\\x_{1}^{n}=e^{n \log z_1}\quad \ }}
\end{equation}
is a single-valued analytic function on the region $E_t$ for each $t\in [0,\frac{2}{7}]$;
\begin{equation}
g_t=\langle w_{(a_4)}', (\tilde{\bf P}(\tilde{\Omega}^{(4)}([\mathcal{Z}]_P)))(\tilde{w}_{(a_{1})},
\tilde{w}_{(a_{3})}, \tilde{w}_{(a_{2})}; x_{1},x_{2}) \rangle_{W^{a_4}}\mbar_{x_{1}^{n}=e^{n \log z_1},x_{2}^{n}=e^{n \log z_2}}
\end{equation}
is a single-valued analytic function on the region $E_t$ for each $t\in (\frac{2}{7},\frac{4}{7}]$;
\begin{equation}
g_t=\langle e^{-x_2 L(1)}w_{(a_4)}', (\tilde{\bf P}([\mathcal{Z}]_P))(\tilde{w}_{(a_{1})},
\tilde{w}_{(a_{2})}, \tilde{w}_{(a_{3})}; x_0, x_{2}) \rangle_{W^{a_4}}\mbar_{x_{0}^{n}=e^{n \log (z_1-z_2)},x_{2}^{n}=e^{n \log (-z_2)}}
\end{equation}
is a single-valued analytic function on the region $E_t$ for each $t\in (\frac{4}{7},\frac{5}{7}]$; and
\begin{equation}
g_t=\langle w_{(a_4)}', (\tilde{\bf P}(\mathcal{B}(\tilde{\Omega}^{(4)}([\mathcal{Z}]_P))))(\tilde{w}_{(a_{3})},
\tilde{w}_{(a_{1})}, \tilde{w}_{(a_{2})};  x_{2},x_1) \rangle_{W^{a_4}}\mbar_{x_{1}^{n}=e^{n \log z_1},x_{2}^{n}=e^{n \log z_2}}
\end{equation}
is a single-valued analytic function on the region $E_t$ for each $t\in (\frac{5}{7},1]$.
Next, we shall show that $\{(g_t,E_t): 0\leq t \leq 1 \}$ is an analytic continuation along $\gamma$.

Firstly, by (\ref{b18}) we have
\begin{eqnarray}
\lefteqn{g_{\frac{2}{7}}=\langle w_{(a_4)}', (\tilde{\bf I}(\tilde{\Omega}^{(1)} \F (\tilde{\Omega}^{(4)}([\mathcal{Z}]_P))))(\tilde{w}_{(a_{3})},
\tilde{w}_{(a_{1})}, \tilde{w}_{(a_{2})};  x_{0},x_1) \rangle_{W^{a_4}}\mbar_{x_{0}^{n}=e^{n \log (z_2-z_1)},x_{1}^{n}=e^{n \log z_1}}}  \nno\\
&& =\langle w_{(a_4)}', (\tilde{\bf P}(\tilde{\Omega}^{(4)}([\mathcal{Z}]_P)))(\tilde{w}_{(a_{1})},
\tilde{w}_{(a_{3})}, \tilde{w}_{(a_{2})}; x_{1},x_{2}) \rangle_{W^{a_4}}\mbar_{x_{1}^{n}=e^{n \log z_1},x_{2}^{n}=e^{n \log z_2}}\qquad\quad\label{b20}
\end{eqnarray}
on the region $E_{\frac{2}{7}}$. Secondly,
it can be derived from the skew-symmetry isomorphism that on the region $E_{\frac{4}{7}}$,
\begin{eqnarray}
\lefteqn{ g_{\frac{4}{7}}=\langle w_{(a_4)}', (\tilde{\bf P}(\tilde{\Omega}^{(4)}([\mathcal{Z}]_P)))(\tilde{w}_{(a_{1})},
\tilde{w}_{(a_{3})}, \tilde{w}_{(a_{2})}; x_{1},x_{2}) \rangle_{W^{a_4}}\mbar_{x_{1}^{n}=e^{n \log z_1},x_{2}^{n}=e^{n \log z_2}}}\nno\\
&& =\langle e^{x_2 L(1)}w_{(a_4)}', (\tilde{\bf P}([\mathcal{Z}]_P))(\tilde{w}_{(a_{1})},
\tilde{w}_{(a_{2})}, \tilde{w}_{(a_{3})}; x_0, e^{-\pi i}x_{2}) \rangle_{W^{a_4}}\mbar_{x_{0}^{n}=e^{n \log (z_1-z_2)},x_{2}^{n}=e^{n \log z_2}}\nno\\
&& =\langle e^{-x_2 L(1)}w_{(a_4)}', (\tilde{\bf P}([\mathcal{Z}]_P))(\tilde{w}_{(a_{1})},
\tilde{w}_{(a_{2})}, \tilde{w}_{(a_{3})}; x_0, x_{2}) \rangle_{W^{a_4}}\mbar_{x_{0}^{n}=e^{n \log (z_1-z_2)},x_{2}^{n}=e^{n \log (-z_2)}}.\qquad\quad
\end{eqnarray}
 Thirdly, we shall prove that on the region $E_{\frac{5}{7}}$,
\begin{eqnarray}
\lefteqn{g_{\frac{5}{7}}=\langle w_{(a_4)}', (\tilde{\bf P}(\mathcal{B}(\tilde{\Omega}^{(4)}([\mathcal{Z}]_P))))(\tilde{w}_{(a_{3})},
\tilde{w}_{(a_{1})}, \tilde{w}_{(a_{2})};  x_{2},x_1) \rangle_{W^{a_4}}\mbar_{x_{1}^{n}=e^{n \log z_1},x_{2}^{n}=e^{n \log z_2}}}  \nno\\
&& =\langle e^{-x_2 L(1)}w_{(a_4)}', (\tilde{\bf P}([\mathcal{Z}]_P))(\tilde{w}_{(a_{1})},
\tilde{w}_{(a_{2})}, \tilde{w}_{(a_{3})}; x_0, x_{2}) \rangle_{W^{a_4}}\mbar_{x_{0}^{n}=e^{n \log (z_1-z_2)},x_{2}^{n}=e^{n \log (-z_2)}}.\qquad\quad\label{b21}
\end{eqnarray}
Replacing $(w'_{(a_{4})},z_1,z_2)$ by $(e^{-x_2 L(1)}w_{(a_4)}',z_1-z_2,-z_2)$ in (\ref{s12}), we can derive that
\begin{eqnarray}
\lefteqn{\langle e^{-x_2 L(1)}w_{(a_4)}', (\tilde{\bf P}([\mathcal{Z}]_P))(\tilde{w}_{(a_{1})},
\tilde{w}_{(a_{2})}, \tilde{w}_{(a_{3})}; x_0, x_{2}) \rangle_{W^{a_4}}\mbar_{x_{0}^{n}=e^{n \log (z_1-z_2)},x_{2}^{n}=e^{n \log (-z_2)}}}\nno\\
&& =\langle e^{-x_2 L(1)}w_{(a_4)}', (\tilde{\bf I}(\mathcal{F}([\mathcal{Z}]_P)))(\tilde{w}_{(a_{1})},
\tilde{w}_{(a_{2})}, \tilde{w}_{(a_{3})}; x_1, x_2) \rangle_{W^{a_4}}\mbar_{x_{1}^{n}=e^{n \log z_1},x_{2}^{n}=e^{n \log (-z_2)}}\qquad\label{b22}
\end{eqnarray}
on the region $E_{\frac{5}{7}}$.
Moreover, by the Moore-Seiberg equations (\ref{hexagon1}) and
(\ref{hexagon2}), we have
\begin{equation}
\mathcal{B}\tilde{\Omega}^{(4)}
=\F^{-1}\tilde{\Omega}^{(1)}\F\tilde{\Omega}^{(4)}
=\tilde{\Omega}^{(3)}\F.
\end{equation}
So by skew-symmetry
we have
\begin{eqnarray}
\lefteqn{\langle w_{(a_4)}', (\tilde{\bf P}(\mathcal{B}(\tilde{\Omega}^{(4)}([\mathcal{Z}]_P))))(\tilde{w}_{(a_{3})},
\tilde{w}_{(a_{1})}, \tilde{w}_{(a_{2})};  x_{2},x_1) \rangle_{W^{a_4}}\mbar_{x_{1}^{n}=e^{n \log z_1},x_{2}^{n}=e^{n \log z_2}}}  \nno\\
&& =\langle w_{(a_4)}', (\tilde{\bf P}(\tilde{\Omega}^{(3)}\F([\mathcal{Z}]_P)))(\tilde{w}_{(a_{3})},
\tilde{w}_{(a_{1})}, \tilde{w}_{(a_{2})};  x_{2},x_1) \rangle_{W^{a_4}}\mbar_{x_{1}^{n}=e^{n \log z_1},x_{2}^{n}=e^{n \log z_2}}  \nno\\
&& =\langle e^{x_2 L(1)}w_{(a_4)}', (\tilde{\bf I}(\mathcal{F}([\mathcal{Z}]_P)))(\tilde{w}_{(a_{1})},
\tilde{w}_{(a_{2})}, \tilde{w}_{(a_{3})}; x_1, e^{-\pi i}x_2) \rangle_{W^{a_4}}\mbar_{x_{1}^{n}=e^{n \log z_1},x_{2}^{n}=e^{n \log z_2}}\nno\\
&& =\langle e^{-x_2 L(1)}w_{(a_4)}', (\tilde{\bf I}(\mathcal{F}([\mathcal{Z}]_P)))(\tilde{w}_{(a_{1})},
\tilde{w}_{(a_{2})}, \tilde{w}_{(a_{3})}; x_1, x_2) \rangle_{W^{a_4}}\mbar_{x_{1}^{n}=e^{n \log z_1},x_{2}^{n}=e^{n \log (-z_2)}}\qquad\quad\label{b23}
\end{eqnarray}
on the region $E_{\frac{5}{7}}$.
Therefore, (\ref{b21}) holds by (\ref{b22}) and (\ref{b23}).

So to sum up, $\{(g_t,E_t): 0\leq t \leq 1 \}$ is an analytic continuation along $\sigma$.

Since $\mathfrak{G}''$ is simply connected, $\sigma\subset \mathfrak{G}''$ and $\sigma(0)=\sigma(1)$, we can derive that $g_0=g_1$ on the region $E_0\cap E_1=E_0$. Moreover, since $g_0$ and $g_1$ are both single-valued analytic functions on the domain $S_2$ which contains $E_0$, we deduce that $g_0=g_1$ on $S_2$.
Namely,
\begin{eqnarray}
\lefteqn{\left.\langle w_{(a_4)}', (\tilde{\bf I}(\tilde{\Omega}^{(1)} \F (\tilde{\Omega}^{(4)}([\mathcal{Z}]_P))))(\tilde{w}_{(a_{3})},
\tilde{w}_{(a_{1})}, \tilde{w}_{(a_{2})};  x_{0},x_1) \rangle_{W^{a_4}}\right|_{x_{0}^{n}=e^{n \log (z_2-z_1)},x_{1}^{n}=e^{n \log z_1}}}\nno\\
&& =\langle w_{(a_4)}', (\tilde{\bf P}(\mathcal{B}(\tilde{\Omega}^{(4)}([\mathcal{Z}]_P))))(\tilde{w}_{(a_{3})},
\tilde{w}_{(a_{1})}, \tilde{w}_{(a_{2})};  x_{2},x_1) \rangle_{W^{a_4}}\mbar_{x_{1}^{n}=e^{n \log z_1},x_{2}^{n}=e^{n \log z_2}}\qquad\label{b24}
\end{eqnarray}
on the region $S_2$. Furthermore, since the first line of (\ref{b24}) defined on $R_4$ is a branch of (\ref{b17}), and the second line of (\ref{b24}) on $S_2$ can be naturally analytically extended to the region $R_2$, we conclude that
\begin{equation}\label{b25}
\langle w_{(a_4)}', (\tilde{\bf P}(\mathcal{B}(\tilde{\Omega}^{(4)}([\mathcal{Z}]_P))))(\tilde{w}_{(a_{3})},
\tilde{w}_{(a_{1})}, \tilde{w}_{(a_{2})};  x_{2},x_1) \rangle_{W^{a_4}}\mbar_{x_{1}^{n}=e^{n \log z_1},x_{2}^{n}=e^{n \log z_2}}
\end{equation}
on $R_2$ is a branch of $\Phi(e^{-x_2 L(1)}w_{(a_4)}', \tilde{w}_{(a_{1})},
\tilde{w}_{(a_{2})}, \tilde{w}_{(a_{3})}, [\mathcal{Z}]_P;z_1-z_2,-z_2)$ on $R_2$.
\vspace{0.2cm}

In conclusion of (\ref{b17})-(\ref{b25}), we see that the multivalued analytic functions
\begin{equation}
\langle w_{(a_4)}', (\tilde{\bf P}(\tilde{\Omega}^{(4)}([\mathcal{Z}]_P)))(\tilde{w}_{(a_{1})},
\tilde{w}_{(a_{3})}, \tilde{w}_{(a_{2})}; x_{1},x_{2}) \rangle_{W^{a_4}}\mbar_{x_1=z_1,x_2=z_2},
\end{equation}
\begin{equation}
\langle w_{(a_4)}', (\tilde{\bf P}(\mathcal{B}(\tilde{\Omega}^{(4)}([\mathcal{Z}]_P))))(\tilde{w}_{(a_{3})},
\tilde{w}_{(a_{1})}, \tilde{w}_{(a_{2})};  x_{2},x_1) \rangle_{W^{a_4}}\mbar_{x_1=z_1,x_2=z_2}
\end{equation}
and
\begin{equation}
\langle w_{(a_4)}', (\tilde{\bf I}(\mathcal{F}(\tilde{\Omega}^{(4)}([\mathcal{Z}]_P))))(\tilde{w}_{(a_{1})},
\tilde{w}_{(a_{3})}, \tilde{w}_{(a_{2})}; x_0,x_{2}) \rangle_{W^{a_4}}\mbar_{x_0=z_1-z_2,x_2=z_2}
\end{equation}
are restrictions of $\Phi(e^{-x_2 L(1)}w_{(a_4)}', \tilde{w}_{(a_{1})},
\tilde{w}_{(a_{2})}, \tilde{w}_{(a_{3})}, [\mathcal{Z}]_P;z_1-z_2,-z_2)$ to their domains $|z_{1}|> |z_{2}|>0$, $|z_{2}|>|z_{1}|>0$ and $|z_{2}|>|z_{1}-z_{2}|>0$ respectively.

So choosing (\ref{a24}) as the preferred branch of $\Phi(e^{-x_2 L(1)}w_{(a_4)}', \tilde{w}_{(a_{1})},
\tilde{w}_{(a_{2})}, \tilde{w}_{(a_{3})}, [\mathcal{Z}]_P;z_1-z_2,-z_2)$ on $R_1$, we see that $\Phi(e^{-x_2 L(1)}w_{(a_4)}', \tilde{w}_{(a_{1})},
\tilde{w}_{(a_{2})}, \tilde{w}_{(a_{3})}, [\mathcal{Z}]_P;z_1-z_2,-z_2)$ is an element of $\mathbb{ G}^{a_1, a_3, a_2, a_4}$. Moreover, by (\ref{b15}), (\ref{b16}), (\ref{b18}), (\ref{b19}), (\ref{b24}), (\ref{b25}),  and by the definition of the preferred branches of an element of $\mathbb{ G}^{a_1, a_3, a_2, a_4}$ on $R_2$ and $R_3$, we see that
\begin{equation}
\langle w_{(a_4)}', (\tilde{\bf P}(\mathcal{B}(\tilde{\Omega}^{(4)}([\mathcal{Z}]_P))))(\tilde{w}_{(a_{3})},
\tilde{w}_{(a_{1})}, \tilde{w}_{(a_{2})};  x_{2},x_1) \rangle_{W^{a_4}}\mbar_{x_{1}^{n}=e^{n \log z_1},x_{2}^{n}=e^{n \log z_2}}
\end{equation}
and
\begin{equation}
\langle w_{(a_4)}', (\tilde{\bf I}(\mathcal{F}(\tilde{\Omega}^{(4)}([\mathcal{Z}]_P))))(\tilde{w}_{(a_{1})},
\tilde{w}_{(a_{3})}, \tilde{w}_{(a_{2})}; x_0,x_{2}) \rangle_{W^{a_4}}\mbar_{x_{0}^{n}=e^{n \log (z_1-z_2)},x_{2}^{n}=e^{n \log z_2}}
\end{equation}
are the preferred branches of $\Phi(e^{-x_2 L(1)}w_{(a_4)}', \tilde{w}_{(a_{1})},
\tilde{w}_{(a_{2})}, \tilde{w}_{(a_{3})}, [\mathcal{Z}]_P;z_1-z_2,-z_2)$ on $R_2$ and $R_3$ respectively.

Therefore, we have
\begin{eqnarray}
\lefteqn{\langle w_{(a_4)}', (\tilde{\bf P}(\tilde{\Omega}^{(4)}([\mathcal{Z}]_P)))(\tilde{w}_{(a_{1})},
\tilde{w}_{(a_{3})}, \tilde{w}_{(a_{2})}; x_{1},x_{2}) \rangle_{W^{a_4}}}\nno\\
&& =\iota_{12}\Phi(e^{-x_2 L(1)}w_{(a_4)}', \tilde{w}_{(a_{1})},
\tilde{w}_{(a_{2})}, \tilde{w}_{(a_{3})}, [\mathcal{Z}]_P;z_1-z_2,-z_2),\label{e1:39}
\end{eqnarray}
\begin{eqnarray}
\lefteqn{\langle w_{(a_4)}', (\tilde{\bf P}(\mathcal{B}(\tilde{\Omega}^{(4)}([\mathcal{Z}]_P))))(\tilde{w}_{(a_{3})},
\tilde{w}_{(a_{1})}, \tilde{w}_{(a_{2})};  x_{2},x_1) \rangle_{W^{a_4}}}\nno\\
&& =\iota_{21}\Phi(e^{-x_2 L(1)}w_{(a_4)}', \tilde{w}_{(a_{1})},
\tilde{w}_{(a_{2})}, \tilde{w}_{(a_{3})}, [\mathcal{Z}]_P;z_1-z_2,-z_2),\label{e1:40}
\end{eqnarray}
\begin{eqnarray}
\lefteqn{\langle w_{(a_4)}', (\tilde{\bf I}(\mathcal{F}(\tilde{\Omega}^{(4)}([\mathcal{Z}]_P))))(\tilde{w}_{(a_{1})},
\tilde{w}_{(a_{3})}, \tilde{w}_{(a_{2})}; x_0,x_{2}) \rangle_{W^{a_4}}}\nno\\
&& =\iota_{20} \Phi(e^{-x_2 L(1)}w_{(a_4)}', \tilde{w}_{(a_{1})},
\tilde{w}_{(a_{2})}, \tilde{w}_{(a_{3})}, [\mathcal{Z}]_P;z_1-z_2,-z_2).\label{e1:41}
\end{eqnarray}
Let $\{e^{a_{1}, a_{3}, a_{2}, a_{4}}_{\alpha}\}_{\alpha\in
\mathbb{ A}(a_{1}, a_{3}, a_{2}, a_{4})}$ be a basis of
$\mathbb{ G}^{a_{1}, a_{3}, a_{2}, a_{4}}$ over the ring
\begin{equation}
\mathbb{ C}[x_{1}, x_{1}^{-1}, x_{2}, x_{2}^{-1}, (x_{1}-x_{2})^{-1}].
 \end{equation}
 Then there exist unique
\begin{eqnarray}
&H_{\alpha}(w'_{(a_{4})},\tilde{w}_{(a_{1})}, \tilde{w}_{(a_{3})},
\tilde{w}_{(a_{2})}, \tilde{\Omega}^{(4)}([\mathcal{Z}]_P); x_{1}, x_{2})&\nno\\
&  \in
\mathbb{ C}[x_{1}, x_{1}^{-1}, x_{2}, x_{2}^{-1}, (x_{1}-x_{2})^{-1}]&\label{e1:42}
\end{eqnarray}
for $\alpha\in \mathbb{ A}(a_{1}, a_{3}, a_{2}, a_{4})$, such that only finitely many of them are nonzero and
\begin{eqnarray}
\lefteqn{\Phi(e^{-x_2 L(1)}w_{(a_4)}', \tilde{w}_{(a_{1})},
\tilde{w}_{(a_{2})}, \tilde{w}_{(a_{3})}, [\mathcal{Z}]_P;z_1-z_2,-z_2)}\nno\\
&& =\sum_{\alpha\in \mathbb{ A}(a_{1}, a_{3}, a_{2}, a_{4})}H_{\alpha}(w'_{(a_{4})},\tilde{w}_{(a_{1})}, \tilde{w}_{(a_{3})},
\tilde{w}_{(a_{2})}, \tilde{\Omega}^{(4)}([\mathcal{Z}]_P); x_{1}, x_{2}) e^{a_{1}, a_{3}, a_{2}, a_{4}}_{\alpha}.\qquad\  \label{e1:38}
\end{eqnarray}
Recall that $\tilde{\Omega}^{(4)}$ is an isomorphism and that
\begin{equation}\label{a38}
\tilde{\Omega}^{(4)}\left(\pi_P(\coprod_{a_{5}\in \mathcal{ A}}
\mathcal{ V}_{a_{1}a_{5}}^{a_{4}}\otimes
\mathcal{ V}_{a_{2}a_{3}}^{a_{5}})\right)=\pi_P(\coprod_{a_{5}\in \mathcal{ A}}
\mathcal{ V}_{a_{1}a_{5}}^{a_{4}}\otimes
\mathcal{ V}_{a_{3}a_{2}}^{a_{5}})
\end{equation}
for any $a_1,\cdots,a_4\in\A$. So we can define linear maps
\begin{eqnarray}
\lefteqn{f^{a_{1}, a_{3}, a_{2},
a_{4}}_{\alpha}(\tilde{w}_{(a_{1})}, \tilde{w}_{(a_{3})},
\tilde{w}_{(a_{2})}):}\nno\\
&&\hspace{4em}
\pi_P(\coprod_{a_{5}\in \mathcal{ A}}
\mathcal{ V}_{a_{1}a_{5}}^{a_{4}}\otimes
\mathcal{ V}_{a_{3}a_{2}}^{a_{5}})\to  W^{a_{4}}[x_{1}, x_{1}^{-1}, x_{2},x_{2}^{-1}][[x_{2}/x_{1}]],
\end{eqnarray}
\begin{eqnarray}
\lefteqn{g^{a_{1}, a_{3}, a_{2},a_{4}}_{\alpha}(\tilde{w}_{(a_{1})}, \tilde{w}_{(a_{3})},\tilde{w}_{(a_{2})}):}\nno\\
&&\hspace{4em}
\pi_P(\coprod_{a_{5}\in \mathcal{ A}}
\mathcal{ V}_{a_{3}a_{5}}^{a_{4}}\otimes
\mathcal{ V}_{a_{1}a_{2}}^{a_{5}})  \to  W^{a_{4}}[x_{1}, x_{1}^{-1}, x_{2},
x_{2}^{-1}][[x_{1}/x_{2}]],
\end{eqnarray}
\begin{eqnarray}
\lefteqn{ h^{a_{1}, a_{3}, a_{2},
a_{4}}_{\alpha}(\tilde{w}_{(a_{1})}, \tilde{w}_{(a_{3})},
\tilde{w}_{(a_{2})}):}\nno\\
&&\hspace{4em}
\pi_I(\coprod_{a_{5}\in \mathcal{ A}}
\mathcal{ V}_{a_{1}a_{3}}^{a_{5}}\otimes
\mathcal{ V}_{a_{5}a_{2}}^{a_{4}}) \to W^{a_{4}}[x_{0},
x_{0}^{-1}, x_{2}, x_{2}^{-1}][[x_{0}/x_{2}]]
\end{eqnarray}
by
\begin{eqnarray}
\lefteqn{\langle w_{(a_4)}', f^{a_{1}, a_{3}, a_{2},
a_{4}}_{\alpha}(\tilde{w}_{(a_{1})}, \tilde{w}_{(a_{3})},
\tilde{w}_{(a_{2})}, \tilde{\Omega}^{(4)}([\mathcal{Z}]_P); x_{1}, x_{2})\rangle_{W^{a_4}}}\nno\\
&& =\iota_{12} H_{\alpha}(w'_{(a_{4})},\tilde{w}_{(a_{1})}, \tilde{w}_{(a_{3})},
\tilde{w}_{(a_{2})}, \tilde{\Omega}^{(4)}([\mathcal{Z}]_P); x_{1}, x_{2}),\label{e1:43}
\end{eqnarray}
\begin{eqnarray}
\lefteqn{\langle w_{(a_4)}', g^{a_{1}, a_{3}, a_{2},
a_{4}}_{\alpha}(\tilde{w}_{(a_{1})}, \tilde{w}_{(a_{3})},
\tilde{w}_{(a_{2})}, \mathcal{B}(\tilde{\Omega}^{(4)}([\mathcal{Z}]_P)); x_{1}, x_{2})\rangle_{W^{a_4}}}\nno\\
&& =\iota_{21} H_{\alpha}(w'_{(a_{4})},\tilde{w}_{(a_{1})}, \tilde{w}_{(a_{3})},
\tilde{w}_{(a_{2})}, \tilde{\Omega}^{(4)}([\mathcal{Z}]_P); x_{1}, x_{2}),\label{e1:44}
\end{eqnarray}
\begin{eqnarray}
\lefteqn{\langle w_{(a_4)}', h^{a_{1}, a_{3}, a_{2},
a_{4}}_{\alpha}(\tilde{w}_{(a_{1})}, \tilde{w}_{(a_{3})},
\tilde{w}_{(a_{2})}, \mathcal{F}(\tilde{\Omega}^{(4)}([\mathcal{Z}]_P)); x_0, x_2)\rangle_{W^{a_4}}} \nno\\
&& =\iota_{20} H_{\alpha}(w'_{(a_{4})},\tilde{w}_{(a_{1})}, \tilde{w}_{(a_{3})},
\tilde{w}_{(a_{2})}, \tilde{\Omega}^{(4)}([\mathcal{Z}]_P); x_2+x_0, x_{2})\label{e1:45}
\end{eqnarray}
for $w_{(a_4)}'\in (W^{a_4})'$, $\mathcal{ Z}\in \coprod_{ a_{5}\in \mathcal{ A}}
\mathcal{ V}_{a_{1}a_{5}}^{a_{4}}\otimes
\mathcal{ V}_{a_{2}a_{3}}^{a_{5}}$ and $\alpha\in \mathbb{ A}(a_{1}, a_{3}, a_{2}, a_{4})$. Then by (\ref{e1:39}), (\ref{e1:40}), (\ref{e1:41}) and (\ref{e1:38}),  we have
\begin{eqnarray}
\lefteqn{(\tilde{\bf P}(\tilde{\Omega}^{(4)}([\mathcal{Z}]_P)))(\tilde{w}_{(a_{1})},
\tilde{w}_{(a_{3})}, \tilde{w}_{(a_{2})}; x_{1},x_{2})}\nno\\
&& =\sum_{\alpha\in \mathbb{ A}(a_{1}, a_{3}, a_{2}, a_{4})} f^{a_{1}, a_{3}, a_{2},
a_{4}}_{\alpha}(\tilde{w}_{(a_{1})}, \tilde{w}_{(a_{3})},
\tilde{w}_{(a_{2})}, \tilde{\Omega}^{(4)}([\mathcal{Z}]_P); x_{1}, x_{2}) \iota_{12}(e^{a_{1}, a_{3}, a_{2}, a_{4}}_{\alpha}),\qquad \label{a39}
\end{eqnarray}
\begin{eqnarray}
\lefteqn{ (\tilde{\bf P}(\mathcal{B}(\tilde{\Omega}^{(4)}([\mathcal{Z}]_P))))(\tilde{w}_{(a_{3})},
\tilde{w}_{(a_{1})}, \tilde{w}_{(a_{2})};  x_{2},x_1)} \nno\\
&& =\sum_{\alpha\in \mathbb{ A}(a_{1}, a_{3}, a_{2}, a_{4})} g^{a_{1}, a_{3}, a_{2},
a_{4}}_{\alpha}(\tilde{w}_{(a_{1})}, \tilde{w}_{(a_{3})},
\tilde{w}_{(a_{2})}, \mathcal{B}(\tilde{\Omega}^{(4)}([\mathcal{Z}]_P)); x_{1}, x_{2}) \iota_{21}(e^{a_{1}, a_{3}, a_{2}, a_{4}}_{\alpha}),\qquad\quad
\end{eqnarray}
\begin{eqnarray}
\lefteqn{ (\tilde{\bf I}(\mathcal{F}(\tilde{\Omega}^{(4)}([\mathcal{Z}]_P))))(\tilde{w}_{(a_{1})},
\tilde{w}_{(a_{3})}, \tilde{w}_{(a_{2})}; x_0,x_{2})} \nno\\
&& =\sum_{\alpha\in \mathbb{ A}(a_{1}, a_{3}, a_{2}, a_{4})} h^{a_{1}, a_{3}, a_{2},
a_{4}}_{\alpha}(\tilde{w}_{(a_{1})}, \tilde{w}_{(a_{3})},
\tilde{w}_{(a_{2})}, \mathcal{F}(\tilde{\Omega}^{(4)}([\mathcal{Z}]_P)); x_0, x_2)
\iota_{20} (e^{a_{1}, a_{3}, a_{2}, a_{4}}_{\alpha}).\qquad\quad
\end{eqnarray}
Moreover, by (\ref{e1:42}) and Proposition \ref{p3.1.1}, we have
 \begin{eqnarray}
\lefteqn{x_{0}^{-1}
\delta\left(\frac{x_{1}-x_{2}}{x_{0}}\right)
\iota_{12}  H_{\alpha}(w'_{(a_{4})},\tilde{w}_{(a_{1})}, \tilde{w}_{(a_{3})},
\tilde{w}_{(a_{2})}, \tilde{\Omega}^{(4)}([\mathcal{Z}]_P); x_{1}, x_{2})}\nno\\
&& \quad -x_{0}^{-1}\delta\left(\frac{x_{2}-x_{1}}{-x_{0}}\right)
\iota_{21} H_{\alpha}(w'_{(a_{4})},\tilde{w}_{(a_{1})}, \tilde{w}_{(a_{3})},
\tilde{w}_{(a_{2})}, \tilde{\Omega}^{(4)}([\mathcal{Z}]_P); x_{1}, x_{2})
\nno\\
&& =x_{2}^{-1}\delta\left(\frac{x_{1}-x_{0}}{x_{2}}\right)
\iota_{20}  H_{\alpha}(w'_{(a_{4})},\tilde{w}_{(a_{1})}, \tilde{w}_{(a_{3})},
\tilde{w}_{(a_{2})}, \tilde{\Omega}^{(4)}([\mathcal{Z}]_P); x_2+x_0, x_{2})\qquad \label{e1:46}
\end{eqnarray}
for $\alpha\in \mathbb{ A}(a_{1}, a_{3}, a_{2}, a_{4})$.
Since $w_{(a_4)}'\in (W^{a_4})'$ is arbitrary, by (\ref{e1:43}), (\ref{e1:44}), (\ref{e1:45}) and (\ref{e1:46}), the Jacobi identity holds for the ordered triple $(\tilde{w}_{(a_{1})},\tilde{w}_{(a_{3})},
\tilde{w}_{(a_{2})})$:
\begin{eqnarray}
\lefteqn{x_{0}^{-1}
\delta\left(\frac{x_{1}-x_{2}}{x_{0}}\right)
f^{a_{1}, a_{3}, a_{2},
a_{4}}_{\alpha}(\tilde{w}_{(a_{1})}, \tilde{w}_{(a_{3})},
\tilde{w}_{(a_{2})}, \tilde{\Omega}^{(4)}([\mathcal{Z}]_P); x_{1}, x_{2})}\nno\\
&&\quad -x_{0}^{-1}\delta\left(\frac{x_{2}-x_{1}}{-x_{0}}\right)
g^{a_{1}, a_{3}, a_{2},
a_{4}}_{\alpha}(\tilde{w}_{(a_{1})}, \tilde{w}_{(a_{3})},
\tilde{w}_{(a_{2})}, \mathcal{ B}(\tilde{\Omega}^{(4)}([\mathcal{Z}]_P)); x_1, x_2)
\nno\\
&&=x_{2}^{-1}\delta\left(\frac{x_{1}-x_{0}}{x_{2}}\right)
h^{a_{1}, a_{3}, a_{2},
a_{4}}_{\alpha}(\tilde{w}_{(a_{1})}, \tilde{w}_{(a_{3})},
\tilde{w}_{(a_{2})}, \mathcal{ F}( \tilde{\Omega}^{(4)}([\mathcal{Z}]_P)); x_{0}, x_{2})\label{a40}
\end{eqnarray}
for $\mathcal{ Z}\in \coprod_{ a_{5}\in \mathcal{ A}}
\mathcal{ V}_{a_{1}a_{5}}^{a_{4}}\otimes
\mathcal{ V}_{a_{2}a_{3}}^{a_{5}}$ and $\alpha\in \mathbb{ A}(a_{1}, a_{3}, a_{2}, a_{4})$. By (\ref{a38}) we see that $\tilde{\Omega}^{(4)}([\mathcal{Z}]_P)$ in (\ref{a39})-(\ref{a40}) can be any element in $\pi_P(\coprod_{a_{5}\in \mathcal{ A}}
\mathcal{ V}_{a_{1}a_{5}}^{a_{4}}\otimes
\mathcal{ V}_{a_{3}a_{2}}^{a_{5}})$. So the Jacobi identity holds for the ordered triple $(\tilde{w}_{(a_{1})},\tilde{w}_{(a_{3})},
\tilde{w}_{(a_{2})})$.
\epfv
\vspace{0.2cm}

{\it Proof of Theorem \ref{t1}.}
 Since the permutations $(1\ 2)$ and $(2\ 3)$ generate the symmetric group $S_3$, in summary of Theorems \ref{l2} and \ref{l3}, we can conclude that the Jacobi identity holds for the triple $(\tilde{w}_{(a_{\tau(1)})},\tilde{w}_{(a_{\tau(2)})},
\tilde{w}_{(a_{\tau(3)})})$ for any $\tau\in S_3$. Thus Theorem \ref{t1} holds.
\epfv

\vspace{2em}

\noindent {\small \sc School of Mathematical Sciences,
University of Chinese Academy of Sciences, Beijing 100049, China}

\noindent {\em E-mail address}: chenling2013@ucas.ac.cn

\end{document}